\documentclass[11pt]{amsart}
\usepackage{amsmath, amsfonts, graphicx}
\newtheorem{thm}{Theorem}
\newtheorem{prop}{Proposition}
\newtheorem{lem}{Lemma}
\newtheorem{cor}{Corollary}
\newtheorem{defn}{Definition}
\newcommand{\C}{\mathbb{C}}
\newcommand{\h}{\mathbb{H}}
\newcommand{\lag}{\mathcal{L}}

\newcommand{\Z}{\mathbb{Z}}
\newcommand{\M}{\mathcal{M}}
\newcommand{\fp}{Fuk (W, U, a, \{\delta_i\})}
\newcommand{\fpt}{Fuk (W_t, U, a, \{\delta_{i, t}\})}
\newcommand{\R}{\mathbb{R}}
\newcommand{\pb}{{\bf P}}
\newcommand{\T}{\mathcal{T}}
\newcommand{\D}{\mathcal{D}^b}
\newcommand{\p}{\mathbb{P}}
\newcommand{\A}{{\bf a}}
\newcommand{\e}{\mathcal{E}}
\newcommand{\ti}{{\tilde i}}
\newcommand{\tV}{{\tilde V}}
\newcommand{\tW}{{\tilde W}}
\newcommand{\blo}{\mathcal{X}}
\newcommand{\gx}{\mathcal{X}}
\newcommand{\gm}{\mathbb{G}_m}
\newcommand{\gmc}{(\mathbb{C}^*)^n}
\newcommand{\nin}{\backslash }
\newcommand{\st}{\mathcal{O}}
\newcommand{\tr}{\triangle}
\newcommand{\Spec}{{\bf Spec}}
\newcommand{\aff}{\mathbb{A}}
\newcommand{\F}{\mathcal{F}}
\newcommand{\gy}{\mathcal{Y}}
\newcommand{\we}{{W^{-1}_\A (\varepsilon)}}
\newcommand{\wq}{{W^{-1}_\A (q)}}
\title{Weighted Blowups and Mirror Symmetry for Toric Surfaces}
\author{Gabriel Kerr}
\begin{document}

\maketitle
\begin{abstract} This paper explores homological mirror symmetry for weighted blowups of toric varieties. It will be shown that
both the A-model and B-model categories have natural semiorthogonal
decompositions. An explicit equivalence of the right orthogonal
categories will be shown for the case of toric surfaces.
\end{abstract}
\section{Introduction}

The homological mirror symmetry (HMS) conjecture was proposed by
Kontsevich \cite{Kontsevich} in 1994 as an attempt to gain a deeper
mathematical understanding of mirror symmetry. Since this time, many
papers have confirmed various versions of HMS. This paper explores
the relation between the B-model of a toric stack and the A-model of
its mirror Landau-Ginzburg model which is one version of HMS. Here
the B-model of the toric stack $\blo$ gives rise to the derived
category $\D (\blo )$ of coherent sheaves on the stack (or
equivariant sheaves on an atlas). The mirror is given by considering
the complex torus with a superpotential $W$ and constructing the
derived Fukaya category $\D (Fuk ((\C^*)^n, W))$ as in
\cite{Seidbook}. For this version of HMS, the conjecture is that
these two triangulated categories are equivalent.

This version of HMS has been confirmed for smooth Fano toric
surfaces \cite{Ueda} as well as weighted projective planes and
Hirzebruch surfaces \cite{AKO}. In this paper we consider the case
when $\blo$ is obtained by taking the weighted projective blowup of
a point on a toric variety $X$. For the standard blowup of a point,
it is well known that $\D (\blo )$ has a natural semi-orthogonal
decomposition in which one piece is independent of the original
variety $X$ \cite{Bondal}. This result will be extended to the case
of weighted blowups (Theorem $2$) where now the independent piece of
$\D (Coh (\blo ))$ also depends on the weights.

On the mirror side, there is a natural degeneration of the potential
$W_\blo$ mirror to $\blo$ into two potential functions, $W_X$ and
$\tilde W$. The first of these, $W_X$ is the potential mirror to the
original variety $X$ while $\tilde W$ depends only on the weights in
the weighted blowup. It will be shown that the derived Fukaya
category of $W_\blo$ admits a semi-orthogonal decomposition in which
the two pieces are the derived Fukaya categories of $W_X$ and
$\tilde W$ respectively (Theorem $6$). Thus the strategy employed is
to show the equivalence of these two categories with those given in
the decomposition of $\D (\blo )$. By itself, this is not enough to
prove the HMS conjecture for these Fano toric stacks as the
interaction between the two categories in the decomposition is
neglected. However, it gives strong evidence for the truth of the
conjecture.

We examine the two dimensional case in detail and show an explicit
equivalence between the triangulated category associated to the
weighted blow up and the derived Fukaya category associated to
$\tilde W$ (Theorem $8$). This result along with the results
established for Fano toric surfaces and weighted projective planes
yields a large class of toric stacks for which our strategy proves
successful.

In section $2$ we will define the weighted projective blowup of a
smooth variety. This procedure is analogous to the usual blowup,
where in this case the (reduced) exceptional divisor is a weighted
projective space. We then show that the derived category of the
blowup admits a semi-orthogonal decomposition. The part of this
decomposition corresponding to the exceptional divisor has an
exceptional collection whose quiver algebra is described explicitly.

Section $3$ is completely independent of the second section. Here we
address issues related to the derived Fukaya category of a
Landau-Ginzburg model. After a general introduction to the Fukaya
category of a Landau-Ginzburg model, we define the notion of a
partial Lefschetz fibration and its associated Fukaya category. The
advantage of this definition is that the derived Fukaya categories
of partial Lefschetz fibrations are invariant under perturbations of
the potential. We then address the case of a superpotential $W_\tr$
on the complex torus which is a Laurent polynomial with Newton
polytope $\tr$. In certain situations, there is a subdivision $\tr =
\tr_1 \cup \tr_2$ satisfying appropriate conditions and one can
define two associated potentials $W_{\tr_i}$. Using techniques from
\cite{GKZ}, we show that the derived Fukaya category associated to
$W$ then admits a semiorthogonal decomposition into the derived
Fukaya category associated to the Lefschetz fibration $W_{\tr_1}$
and the derived Fukaya category associated to the partial Lefschetz
fibration $W_{\tr_2}$.

In section $4$ we work through the example of a Fano stack which is
a weighted projective blowup of a smooth toric Fano surface. We give
an explicit isomorphism between the part of the derived Fukaya
category associated to one piece of the decomposed polytope as
detailed in section $3$ and the part of the derived category of
coherent sheaves of the toric stack associated to the weighted
blowup as detailed in section $2$.

{\it Acknowledgements}: I would like Paul Seidel for providing
insightful comments and suggestions on this topic. I would also like
 to thank Joseph Johns and Mohammed Abouzaid for valuable
conversations.

\section{Weighted Projective Blowups}
\subsection{Definition of the weighted projective blowup}
We first recall the notion of weighted projective space. Let $\A =
(a_0, a_1, \ldots , a_n) \in \Z_{>0}^{n + 1}$ and define an action
of $\C^*$ on $\C^{n + 1} \setminus \{0\}$ via $\lambda \cdot (z_0,
z_1, \ldots , z_n) = (\lambda^{a_0} z_0, \lambda^{a_1} z_1, \ldots,
\lambda^{a_n} z_n)$. We define $\pb (\A)$ to be the quotient of
$\C^{n + 1} \setminus \{0\}$ by this action. In general $\pb (\A)$
is not smooth as a scheme, so we would like to regard it as the
coarse space of a stack $\p (\A)$. There are several equivalent
definitions of a stack, however in this paper we will confine
ourselves to defining a particular atlas for the stacks under
consideration.

Recall that a groupoid atlas for a stack consists of two schemes
$M$, $O$ and morphisms $s$, $t$, $m$, $i$ and $e$. One should think
of the groupoid as a category with objects as points of $O$ and
morphisms as points of $M$. The morphisms $s$ and $t$ yield "source"
and "target" maps $M \hbox{}^\rightarrow_\rightarrow O$, $e: M \to
O$ the identity map, $m: M  \times_O M \to M$ composition and $i: O
\to M$ the inverse. To obtain an Artin stack, all of the expected
diagrams should commute, $s$ and $t$ should be flat morphisms and
$(s, t): M \to O \times O$ ought to be separated and quasi-compact.
By a coherent sheaf on the stack $M \hbox{ }
\hbox{}^\rightarrow_\rightarrow \hbox{ } O$ we will mean a coherent
sheaf $\F$ on $O$ together with a canonical isomorphism from $s^*
\F$ to $t^* \F$. For complete exposition on this subject see
\cite{Laumon}.

A common example of a stack is an algebraic group $G$ acting on a
scheme $X$. In this example we let $R = X \times G$ and $U = X$. The
morphisms are $s(x, g) = x$ and $t(x, g) = g\cdot x$ and composition
$m ((g x, h), (x, g)) = (x, hg)$ . This stack is often denoted $[X /
G]$ and called the quotient stack. In this case, coherent sheaves on
$[X / G]$ are simply $G$-equivariant coherent sheaves on $X$. Thus,
to consider $\pb (\A)$ as a stack, we define the stack $\p ( \A) =
[\C^{n + 1} \setminus \{0\} / \C^*]$ where the action is defined in
the first paragraph. Alternatively, one can define the graded ring
$R = \C [x_0, \ldots, x_n]$ where $\hbox{deg} (x_i) = a_i$. Of
course, this is simply the ring of functions on $\C^{n + 1}$ whose
grading reflects the characters of the group action. In more
generality, as discussed in \cite{AKO}, given a graded ring $R =
\oplus_{i \in \Z} R_i$ over a field $k$ and the ideal $I$ generated
by all elements of positive degree in $S$ one can give the following
definition.
\begin{defn} Define the
quotient stack  $\p {\bf roj} (R) = [(\Spec (R) \nin I) / \gm ]$
\end{defn}
Now suppose we have a smooth $(n+1)$-dimensional scheme $X$ with a
point $p$ as the origin in an affine chart $V = \Spec \C [x_0,
\ldots, x_n] \subset X$. The weighted blowup of a point is a
procedure which takes $(X, p, \A , V)$ and produces a stack $\blo$.
In the case of $\A = (1, 1, \ldots, 1)$, $\blo$ will be the usual
blowup of the scheme at the point $p$. Although there is a more
general procedure for weighted blowups of points which do not lie in
an affine chart, we will only need this case for this paper. We have
the following lemma whose proof occupies the rest of this
subsection.
\begin{lem} Given $(X, p, \A , V)$ as above, the weighted projective
blowup $\blo $ is a smooth stack which fits in a fiber square
\begin{equation}
\begin{array}{ccc}
\e & {\buildrel \ti \over \hookrightarrow} & \blo \\
\downarrow & & {\downarrow \scriptstyle \pi} \\ \{p\} &
\hookrightarrow & X \end{array} \end{equation} where $\e$ is the
exceptional divisor of the stack and $\e_{red}\simeq \p(\A)$.
Furthermore, $\pi: \blo \nin \e \to X \nin {p}$ is an equivalence of
schemes.
\end{lem}

Let ${\tilde V} = (\C^{n+1} \nin {0}) \times \C$ and ${\tilde W} =
(X \nin {p}) \times \C^*$; the atlas $\blo $ will be a quotient
stack with $M = (\tW \cup \tV) \times \C^*$ and $O = \tW \cup \tV$.
The action of $\C^*$ on $\tW$ will be multiplication on the second
factor. For $\tV$ we define $$ \lambda \cdot (x_0, \ldots, x_n, y) =
(\lambda^{a_0} x_0, \ldots, \lambda^{a_n} x_n, \lambda^{-1} y)
$$ The intersection of $\tV$ and $\tW$ is $(\C^{n +1} \nin {0})
\times \C^*$ with group action of $\C^*$ multiplication on the
second factor.  The equivariant inclusion maps are $$ j_\tV : (\C^{n
+1} \nin {0}) \times \C^* \to \tV $$
$$j_\tW: (\C^{n +1} \nin {0}) \times \C^* \to \tW$$ The inclusion map
$j_\tW$ is just the inclusion of $(V \nin {0}) \times \C^*$ in
$\tW$, while $j_\tV (s_0, \ldots, s_n, t) = (t^{a_0}s_0, \ldots,
t^{a_n}s_n, t^{-1})$. It is easily checked that $j_\tV$ is
equivariant with respect to the $\C^*$ action.
\begin{defn} Let $U$ be the equivariant pushout:
\begin{displaymath} \begin{array}{ccc}
{(\C^{n +1} \nin {0}) \times \C^*} & {\buildrel j_\tV \over
\hookrightarrow} & \tV \\ {{\scriptstyle j_\tW} \downarrow} & &
{\downarrow } \\ \tW & \hookrightarrow &
 U \end{array} \end{displaymath}
and define $\blo = [U / \C^*]$.
\end{defn}
As in the case of weighted projective spaces, one can view $\tV$ as
$\p {\bf roj} (S)$ where $S$ is the graded ring $\C [z_0, \ldots ,
z_n, y]$ with $\deg (z_i) = a_i$ and $\deg (y) = -1$. One sees that
the irrelevant ideal is $I = (z_0, \ldots, z_n) = S \cdot S_+$ so
that $\Spec (S) \nin I = (\C^{n + 1} \nin {0}) \times \C$ which is
just $\tV$. To define the map $\pi$ we simply take the quotient map
on $\tW$ while on $\tV$ we take the map induced by the homomorphism
of graded rings $\rho: \C [t_0, \ldots, t_n] \to \C [z_0, \ldots,
z_n, y]$ which sends $t_i$ to $z_i \cdot y^{a_i}$.

Thus the exceptional divisor $\e$ has ideal sheaf $J = (z_0 y^{a_0},
\ldots , z_n y^{a_n})$ and structure sheaf $\st_\e = S/ J$. If $b =
\max\{a_0, \ldots, a_n\}$ then the ideal $(y^b)$ clearly satisfies
$I \cdot (y^b) \subset J$ where $I$ is the irrelevant ideal of $S$.
So on $\e$, the ideal sheaf $(y^b)$ is supported on $V(I) \cap \e$
and therefore torsion. Thus,
$$\e = \p {\bf roj} \left( \frac{S}{J + (y^b)}\right)$$ From this we
see that $y \in Rad(S / [J + (y^b)])$ and
$$\e_{red} = \p {\bf roj} \left( \frac{S}{(y)}\right) = \p (\A)$$
Which yields the expected diagram (1). We also have the map $i : \p
(\A) = \e_{red} \to \e \to \tV$ which is induced by the homogeneous
quotient map of graded rings $\psi: S \to \frac{S}{(y)}$.

\subsection{Weighted projective blowups of toric varieties}
In this subsection we will see that the above definition has a
simple expression when $X$ is a smooth toric variety. Useful
references on the essentials of toric varieties are \cite{Fulton}
and \cite{Oda}; we will freely use results from \cite{Cox} as well.

We recall some standard notation. Let $N \simeq \Z^{n + 1}$ be a
lattice, $M = \hbox{Hom}(N, \Z)$, $N_\R := N \otimes \R$ and $M_\R =
M \otimes \R$. Let $\tr$ be a fan with cones $\sigma \subset N_\R$.
We will let $\tr (k)$ be the set of $k$-dimensional cones in $\tr$.
A smooth toric variety $X_\tr$ results from a fan $\tr$ such that
every $\sigma \in \tr (n + 1)$ can be represented as $\sigma =
\R_{\geq 0} \cdot v_0 + \cdots + \R_{\geq 0} \cdot v_n$ where
$\{v_0, \ldots, v_n\}$ generates $N$. For this subsection we will
choose such a $\sigma \in \tr (n + 1)$. This yields an affine chart
$U_\sigma = \Spec \C [\sigma^\vee \cap M] \simeq \Spec \C[v_0^\vee,
\ldots, v_n^\vee] \simeq \aff^{n + 1}$. Let $i: \aff^{n + 1}
\hookrightarrow X_\tr$ be the inclusion. If $p = i(0)$ then from the
correspondence between orbit closures and cones, we have $X_\tr \nin
{p} = X_{\tr \nin \sigma}$. Given $\A = (a_0, \ldots, a_n) \in \Z^{n
+ 1}_{>0}$ let $u = \sum_{i = 0}^n a_i v_i$. For any proper subset
$J \subset \{0, 1, \ldots, n\}$ let
$$\tau_J = \R_{\geq 0} \cdot u + \Sigma_{i \in J} \R_{\geq 0} \cdot
v_i$$ We define a new fan in the expected manner
$$\tr (\sigma, \A) = (\tr \nin \sigma) \cup \left( \cup_{J
} \tau_J\right)$$ Following the construction in subsection $1.1$, we
let $V = U_\sigma$ and $p = i(0)$. Let $\Sigma = \{\eta \in \tr
(\sigma, \A) | \eta \subset \sigma\}$. Then, as $X_{\tr (\sigma, \A)
\nin \cup \tau_J} = X_{\tr \nin \sigma}$ we need only examine
$X_{\Sigma}$ and the map induced by the identity $(N, \Sigma) \to
(N, \sigma)$. In order for the fan $\Sigma$ to yield a stack as
opposed to a singular scheme, one must regard $u$ as the element
generating the one dimensional cone $\{\R_{>0} \cdot u\} \in \Sigma
(1)$. For example, if $n = 1$ and $\A = (2, 2)$ we consider $u$ to
be the generating element of $\R_{> 0} \cdot u$ (as opposed to $v_1
+ v_2$). We then define the map $\alpha : \Z^{\Sigma (1)} \to N$ as
$\alpha (e_i ) = v_i$ for $0 \leq i \leq n$ and $\alpha (e_{n + 1} )
= u$ where $\{e_0, \ldots, e_{n + 1}\}$ is the standard basis . From
this we obtain the exact sequence:
$$ 0 \to M {\buildrel \alpha^\vee \over \longrightarrow} \hbox{ }
\Z^{\Sigma(1)} {\buildrel \beta \over \longrightarrow} A_n (X_\Sigma
) \to 0$$ where, using the standard basis, we have identified
$\Z^{\Sigma (1)}$ with it's dual.

As in \cite{Cox}, for each basis element $e_i \in \Z^{\Sigma (1)}$,
introduce the variable $z_i$ and consider the graded ring $$ S = \C
[z_0, \ldots , z_n, z_{n + 1}]$$ with the grading induced by $\deg
(z_i) = \beta (e_i)$. A quick check shows that $A_n (X_\Sigma )
\simeq \Z$ and that after a choice of generator, we have $\beta
(e_i) = a_i$ for $0 \leq i \leq n$ and $\beta (e_{n + 1} ) = -1$. In
the general case of a fan $F$, irrelevant ideal is the ideal
generated by all monomials $z_{\rho(1)} \cdots z_{\rho(k)}$ such
that $\{\rho (1), \ldots, \rho (k)\} = F (1) \nin \sigma (1)$ with
$\sigma \in F (n)$. Thus, in our situation the irrelevant ideal is
$I = (z_0, \ldots, z_n) = S \cdot S_+$ so as a stack we have
$X_\Sigma = [\Spec (S) \nin I , \C^*]$.

We recognize $\Spec (S) \nin I$ with the action of $\C^*$ given by
the grading as $\tV$. Likewise we have $V = U_\sigma = \Spec (\C
[x_0, \ldots , x_n])$ where we have introduced the variable $x_i$
for $v_i^\vee$. The map $\tV \to V$ is induced by the identity $(N,
\Sigma) \to (N, \sigma)$ and will be given by the map of homogeneous
coordinate rings $\rho :\C [x_0, \ldots , x_n] \to S $ where $\rho
(x_i) = z^{\alpha^\vee (v_i^\vee)} = z_i z_{n + 1}^{a_i}$ as can be
seen from examining the obvious diagram:
\begin{displaymath} \begin{array}{ccccccccc} 0 & \to & M &
{\buildrel = \over \longrightarrow} & \Z^{\sigma(1)} & {
\longrightarrow} & 0 & &
\\ && \parallel && {\downarrow \scriptstyle \alpha^\vee} &&&&
\\ 0 & \to & M & {\buildrel \alpha^\vee \over \longrightarrow} &
\Z^{\Sigma(1)} & {\buildrel \beta \over \longrightarrow} & A_n
(X_\Sigma ) & \to & 0
\end{array} \end{displaymath} For a more general procedure of
obtaining stacks from fans, see \cite{Borisov}.

\subsection{Sheaves on weighted projective blowups}
In this section we will extend some of the results from \cite{AKO}
and \cite{Bondal} to obtain a semi-orthogonal decomposition of the
derived category of coherent sheaves on $\blo$ in terms of the
derived categories of coherent sheaves on $X$ and $\p (\A)$. There
will be an explicit mirror decomposition on the symplectic side in
section $3$. We start with a theorem whose proof is given in either
\cite{AKO} or \cite{Cox}. Recall that $\tV = \p {\bf roj} (S)$ with
$S = \C [z_0, \ldots , z_n, y]$ and gradings and irrelevant ideal
$I$ given in section $2.1$. Define $gr (S)$ ($gr_f (S)$) to be the
category of (finitely generated) graded modules over $S$ and $Tor
(S)$ ($Tor_f (S)$)to be the full subcategory of modules $M$ such
that there exists a $k \in \Z_+$ with $I^k M = 0$. One sees that
$Tor (S)$ is a Serre subcategory of $gr (S)$.
\begin{thm} The category of quasi-coherent (coherent) sheaves on
$\tV$ is equivalent to the categorical quotient of $gr (S)$ ($gr_f
(S)$) by $Tor (S)$. \end{thm} The same theorem applies if one
replaces the ring $S$ by $R$ defined in section 1.1 and $\tV$ by $\p
(\A)$ (see \cite{AKO}). We recall some homological properties of
sheaves over $\p (\A)$. Given a graded module $M = \oplus_{i \in \Z}
M_i$ we let $M(k)$ be the graded module with $M(k)_i = M_{k + i}$.
Now, given $\A = (a_0, \ldots, a_n) \in \Z^{n +1}_{> 0}$ let $l =
\sum_{i = 0}^n a_i$. We recall from \cite{AKO} the following result
on the cohomology of the sheaf $\st (k)$ on $\p (\A)$ obtained from
the graded $R$-module $R(k)$.
\begin{prop} There are isomorphisms \begin{displaymath} H^p (\p
(\A), \st (k)) \simeq \left\{ \begin{array}{cc} R_k & for \hbox{ }
p=0\hbox{, } k \geq 0 \\ R^*_{-k-l} & for \hbox{ } p = n \hbox{, } k
\leq -l
\\ 0 & otherwise  \end{array} \right. \end{displaymath}
\end{prop} Our main concern in this section is with the
derived category of coherent sheaves on $\blo$. Given an Artin stack
$\gy$ we will denote the bounded derived category of coherent
sheaves on $\gy$ by $\D (\gy)$. We recall that this is the
triangulated category obtained from the abelian category $Coh(\gy)$
by localizing all quasi-isomorphisms. Given a proper (flat) map $f$
between two stacks, we will write $f_*$ ( $f^*$ ) for the derived
functors $Rf_*$ and $Lf^*$. Recall that if $A^\bullet$ is an object
in $\D (\gy)$ then $A[n]^\bullet$ is the translated object $A[n]^i =
A^{n + i}$.

From diagram (1) we have a map $i : \p (\A) \to \blo$. Utilizing
Theorem $1$ , we prove
\begin{prop} Given any $j$ and $k$, there are isomorphisms
\begin{displaymath} Ext^\bullet (i_* (\st (k)), i_* (\st (j)))
 \simeq Ext^\bullet (\st (k), \st (j)) \oplus Ext^\bullet (\st (k + 1), \st (j ))[-1]  \end{displaymath} \end{prop}
 \begin{proof} Since the image of $i$ is contained in $\tV$, it suffices to work
 locally with $i : \p (\A) \to \tV$ to compute the $Ext$ groups of the
pushforwards of $\st (p)$. Recall that the map $i : \p (\A) \to \tV$
is induced by the the homomorphism $\psi: S \to R$ with kernel
$(y)$. Observe that, when regarded as a module over $S$, the ideal
$(y)$ is isomorphic to $S (1)$. Thus there is the exact sequence of
$S$-modules $$ 0 \to S (1) {\buildrel m_y \over \longrightarrow} S
\to R \to 0$$ where $m_y$ is multiplication by $y$. For arbitrary
$k$ we of course have
$$0 \to S (k + 1) {\buildrel m_y \over \longrightarrow} S (k) \to
R (k) \to 0$$ Now, as an $S$-module, $R (k )$ is the module
corresponding to $i_* (\st (k ))$ and $i^* \circ i_*(\st (k ))$
corresponds to the module $R \otimes_S R(k)$. The first two modules
in the sequence above are free and therefore acyclic for the functor
$R \otimes_S \hbox{}_{\bf -}$. Thus in $gr (R)$ there is resolution
of $R \otimes_S R(k)$:
\begin{displaymath} \begin{array}{ccccccccc} 0 & \to &  R \otimes_S S (k + 1) &{\buildrel {id \otimes m_y} \over
\longrightarrow}& R \otimes_S S (k)& \to & R  \otimes_S R (k) & \to
& 0 \\ &&  \parallel & & \parallel &&\parallel && \\ 0 & \to & R (k
+ 1) &{\buildrel {0} \over \longrightarrow}& R (k)& \to & R
\otimes_S R (k) & \to & 0 \end{array} \end{displaymath} So one
obtains $R \otimes_S^L R(k) \simeq R(k) \oplus R(k + 1)[1]$ in the
bounded derived category $\D (gr (R))$. Thus we have the equivalence
$$i^* \circ i_*(\st (k )) \simeq \st (k) \oplus \st (k + 1) [1]$$
in $\D (\p (\A))$. Using this and the fact that $i_*$ is right
adjoint to $i^*$ yields the result.
\end{proof}

We will need the following explicit corollary.

\begin{cor} Given any $j$ and $k$ such that $0 \leq k - j \leq l - 2$ there
is an natural isomorphism compatible with composition:
\begin{displaymath} Ext^\bullet (i_* (\st (j)), i_* (\st (k)))
 \simeq R_{k - j} \oplus R_{k - j - 1} [-1]  \end{displaymath}
Where $R_{j -k}$ is placed in degree zero and $R_{k - j - 1} [-1]$
is placed in degree one. Furthermore, there is a natural
$\C^*$-action which is also compatible with composition.
\end{cor}
\begin{proof} The isomorphism is immediate from the previous
propositions. To see that it is compatible with composition we
return to the proof of Proposition $2$. As modules we have $i_*(\st
(j))$ is equivalent to $R(j)$ which is quasi-isomorphic to the
complex $0 \to S(j + 1) \to S(j) \to 0$ and likewise for $i_*(\st
(k))$. Following the proof of Proposition $2$ we see that any
monomial in $S_{k - j}$ which does not lie in the ideal $(y)$ gives
rise to a degree zero chain map. Such monomials correspond naturally
to the $R_{k - j}$ summand in the stated isomorphism. Similarly, the
degree one chain maps correspond to multiplication by monomials in
$S_{k - j - 1}$ not in $(y)$, or $R_{k - j - 1}$. It is clear that
composing two such chain maps corresponds to multiplication of two
such monomials (unless, of course, both are of degree $1$, in which
case the product is zero). The $\C^*$-action is induced by the
$\C^*$-action on $R$.
\end{proof}
 For basic definitions on triangulated categories, see
\cite{Bondal}. We will denote the full triangulated category in $\D
(\blo)$ generated by the object $i_* (\st (k))$ as $\T_k$. Then we
have the following proposition.

\begin{prop} The sequence $\left< \pi^* (\D ( X)) , \T_0 , \ldots , \T_{l - 2}
\right>$
 is an exceptional collection of triangulated subcategories in $\D (\blo )$. \end{prop}

\begin{proof} One can use standard methods to show that these are all
 admissible subcategories of $\D (\blo )$ (see \cite{Bondal}). To see that $\left< \T_0 , \ldots ,
\T_{l -2} \right>$ forms
 an exceptional collection we use Propositions $1$ and $2$. We note that it suffices to check
that $Ext^\bullet (\ti_* (\st (j)), \ti_* (\st (k))) = 0 $ for $0
\leq k < j \leq l - 2$, for this will imply that the categories
generated by these objects are semi-orthogonal. By Proposition $2$
and $1$ we have
\begin{displaymath} \begin{array}{ccc} Ext^\bullet (i_* (\st (j)), i_* (\st (k))) &
 = & Ext^\bullet (\st (j), \st (k)) \oplus Ext^\bullet (\st (j + 1), \st (k))[-1] \\ & =
 & Ext^\bullet (\st , \st (k - j)) \oplus Ext^\bullet
(\st, \st (k - j - 1))[-1] \\ & = & H^\bullet (\p (\A ), \st (k -
j)) \oplus H^{\bullet - 1}(\p (\A ), \st (k - j - 1))
\end{array}
\end{displaymath}
Since $0 \leq k < j \leq l - 2$ we have $2 - l \leq k - j < 0$ and
$1 - l \leq k - j - 1 < -1$. Applying Proposition 1 to each of these
cases we see that the right hand side of the above equation is zero.

Thus we need only show that $\left< \pi^* (\D ( X)) , \T_k \right>$
is a semi-orthogonal sequence for $0 \leq k \leq l - 2$. Again it
suffices to check this on a sheaf $\F \in Coh (X)$ and $\st (k)$,
i.e. we need to check that
\begin{displaymath} Ext^\bullet (i_* (\st (k)), \pi^* (\F )) = 0 \end{displaymath}
To do this we apply Serre Duality and the adjunction formula. Recall
from \cite{AKO} that the canonical sheaf of the stack $\p (\A )$ is
$\st (-l )$. Now, $i: \p (\A) \to \blo $ is a closed embedding, and
letting $Y = i(\p (\A ))$, we have that
 the sheaf $\st ([Y])$ on $\tV$ is equivalent to the module $S(-1)$ (as $(y) = S(1)$ is the ideal sheaf of the image).
 This implies $i^* (\st ([Y])) = \st (-1)$.
This fact and the adjunction formula yields
\begin{displaymath} \begin{array}{ccc} i^* (\omega_\gx ) & = & i^* (\omega_\gx ) \otimes
\st (-1 ) \otimes \st (1) \\ & = & i^* (\omega_\gx ) \otimes i^*
(\st ([Y])) \otimes \st (1)
 \\ & = & i^* (\omega_\gx \otimes \st ([Y])) \otimes \st (1) \\ & = & \omega_{\p (\A )} \otimes
 \st (1) \\ & = & \st (1 - l ) \end{array} \end{displaymath}
 Now we apply this with Serre Duality to see
\begin{displaymath} \begin{array}{ccc}  Ext^\bullet (i_* (\st (k)), \pi^* (\F )) & = &
 Ext^\bullet (\pi^* (\F ), \omega_\gx \otimes i_* (\st (k)))^\vee \\ & = & Ext^\bullet
 (\pi^* (\F ),i_*(i^*(\omega_\gx ) \otimes  \st (k)))^\vee \\ & = &
Ext^\bullet (\pi^* (\F ), i_* (\st (1 - l + k)))^\vee \\ & = &
Ext^\bullet (\F , (\pi \circ i )_* (\st (1 - l + k) ))^\vee \\ & = &
0 \end{array} \end{displaymath} Indeed, to verify the last line of
the equation observe that $\pi \circ i = inc \circ \rho$ where
$\rho$ is projection to a point and $inc$ is the inclusion of the
point to $p \in X$. Thus $(\pi \circ i)_*$ factors through the
derived global section functor $\rho_* =R\Gamma$. But using
Proposition $2$ and the fact that $1 - l \leq 1 - l + k \leq -1$ we
see that $R^\bullet \Gamma ( \st ( 1 - l + k)) = H^\bullet (\p (\A
), \st (1 - l + k) ) = 0$ yielding the asserted equality.
\end{proof}

Of course, we would also like to see whether or not this is a
complete exceptional collection; i.e., we would like to see if $\D
(\gx )$ is generated as a triangulated category by the subcategories
in the sequence. This is indeed the case, however owing to the fact
that the exceptional divisor $\e$ is not necessarily reduced, the
proof is fairly technical.
\begin{thm} The sequence $\left< \pi^* (\D ( X)) , \T_0 , \ldots , \T_{l - 2} \right>$
 gives a semi-orthogonal decomposition of $\D (\gx )$. \end{thm}
\begin{proof} By Proposition $3$ we need only check that the collection is complete. Let $\T$ be
the triangulated subcategory generated by $\left< \pi^* (\D ( X))
, \T_0 , \ldots , \T_{l - 2} \right>$. \\[8pt]
 {\it Claim} 1: If $i_* (\D (\p (\A )))$ is contained in $\T$ then $\T = \D ( \gx
 )$. \\[8pt]
To verify this claim let $Se$ be the Serre functor and suppose $A
\in {}^\perp \T$. Then for all $B \in \D (X)$ we have $Hom (\pi^*
(B), Se(A))^{\vee} = Hom (A, \pi^* (B)) = 0$ which implies $\pi_*
(Se (A)) = 0$. Thus $Se(A)$ and therefore $A$ has support on $Y =
i(\p (\A))$. If $A \not= 0$ then the identity morphism implies $ 0
\not= Hom (i^* (A), i^*(A)) = Hom (A, i_* (i^* (A)))$ which implies
$A \notin {}^\perp \T$ contradicting our assumption. Thus $A = 0$
and ${}^\perp \T = 0$. Since $\T$ is admissible, we have that $\T =
\D (\gx)$.

To complete the proof we need to establish that $i_* (\D (\p (\A
)))$ is a subcategory of $\T$. For this recall from \cite{AKO} that
$\left< \st , \ldots, \st (l - 1) \right>$ is a full exceptional
collection for $\p (\A )$. So if $i_* (\st (l - 1)) \in \T$ then
$\left< i_* (\st ) , \ldots, i_*(\st (l - 1)) \right> = i_* ( \left<
\st , \ldots, \st (l - 1) \right>) = i_* (\D (\p (\A )))$
is contained in $\T$. Thus we need only show \\[8pt] {\it Claim} 2: $i_* (\st
(l - 1)) \in \T$. \\[8pt]
Let us outline the proof of the claim. For each $0 \leq {m} \leq l$
we will introduce an object $K_m^\bullet \in \D ( \gx )$ as well as
a distinguished triangle $K_m^\bullet \to K_{m - 1}^\bullet \to L_{m
- 1}^\bullet$ where $L_{m -1}^\bullet$ is an object of $\T$ for $1
\leq m < l$, $K_0^\bullet \in \pi^* (\D ( X))$, $K_l^\bullet = 0$
and $L_{l - 1}^\bullet$ is a direct sum of objects in $\T$ and shift
of $i_* (\st (l - 1 ))$. An easy induction argument then shows that
$L_{l - 1}^\bullet$ is in $\T$ and therefore $i_* (\st (l - 1)) \in
 \T$. We will now fill in the details.

 We require some notation to verify this claim. Let $0 \leq m \leq
l$ and $c_m = (c_{m, 0}, c_{m, 1}, \ldots , c_{m , n})$ be the
unique element in $\Z_{\geq 0}^{n + 1}$ satisfying
\begin{displaymath} \begin{array}{cc} (1) & \Sigma_{j =0}^n c_{m, j} = m \\
(2) & \Sigma_{j =0}^k c_{m, j} = \Sigma_{j =0}^k a_j \hbox{ if }
\Sigma_{j =0}^k a_j \leq m
\\ (3) & c_{m, k} = 0 \hbox{ if } \Sigma_{j =0}^{k - 1} a_j \geq m \end{array}
\end{displaymath}
Alternatively, one can see that $c_m$ is the $m$-th term in the
lexigraphically ordered sequence in $\Z_{\geq 0}^{n + 1}$ with each
$j$-coordinate between $0$ and $a_j$ and starting at $(0, \ldots,
0)$.

By Theorem $2$, we can regard sheaves over $\tV$ as graded modules
over $S = \C [ z_0, \ldots, z_n, y]$. For $0 \leq m \leq l$, let
\begin{displaymath} N_m = \bigoplus_{j = 0}^n S(c_{m, j})
\end{displaymath}
and define the degree zero element $$s_m = \oplus_{j = 0}^n z_j
y^{a_j - c_{m,j}} \in N_m$$ We now define the complex $K^\bullet_m$
to be the Koszul complex, \begin{displaymath} 0  \to   S {\buildrel
{\wedge s_m} \over \longrightarrow} N_m  {\buildrel {\wedge s_m}
\over \longrightarrow}  \wedge^2 N_m  \to  \cdots  {\buildrel
{\wedge s_m} \over \longrightarrow}  \wedge^{n + 1} N_m  \to  0
\end{displaymath} One should note that $s_m$ does not generally come from a
regular sequence for a given $m$. Following the above outline, we
would like to define a chain map $f_m^\bullet : K^\bullet_m \to
K^\bullet_{m - 1}$. Observe that $c_m - c_{m - 1} = (0, \ldots, 0,
1, 0, \ldots , 0)$ where $1$ is in the $i$-th coordinate for some
$i$. We can define ${\tilde f_m} : N_m \to N_{m - 1}$ by $\left(
\oplus_{j \not= i} id \right) \oplus m_y$ where $m_y$ is
multiplication by $y$ in the $i$-th summand. One easily shows that
${\tilde f_m} (s_m ) = s_{m - 1}$ which implies that ${\tilde f_m}$
induces an injective chain map $f_m^\bullet : K^\bullet_m \to
K^\bullet_{m - 1}$ (by injective we will mean that $f^j_m$ is
injective for every $j$ ). By the proof of Proposition $2$, we see
that $coker({\tilde f_m}) = R(c_{m-1 , i})$ and extrapolating this
to the chain complex one can show
$$coker (f_m^k) \simeq \bigoplus_{ i \in J, |J| = k} R ( \Sigma_{j
\in J} c_{m -1, j} ) =: L^k_{m - 1}$$ where the direct sum is taken
over all $J \subset \{0, \ldots, n\}$. As ${\tilde f_m} (s_m ) =
s_{m - 1}$, the differential in the cokernel $L^\bullet_{m - 1}$
must be zero. Thus, on the level of chain complexes, we have an
exact sequence $$0 \to K^\bullet_m \to K^\bullet_{m - 1} \to
L^\bullet_{m -1} \to 0$$ where $L^\bullet_{m - 1}$ is a direct sum
of the modules given above plus shifts. In the derived category,
this gives us a distinguished triangle, $K^\bullet_m \to
K^\bullet_{m - 1} \to L^\bullet_{m -1}$. Now, for $1 \leq m < l$ and
any subset $J \subset \{0, \ldots, n\}$, condition ($1$) implies
that $$0 \leq \sum_{j \in J} c_{m - 1, j} \leq m - 1 < l - 1$$
Therefore, for $1 \leq m < l$, all of the summands in $L^\bullet_{m
-1}$ are equal to $R(k)$ for $0 \leq k \leq l - 2$. But these
modules are equivalent to the $i_* (\st( k ))$ which appear in our
exceptional collection. Thus $L^\bullet_{m -1} \in \T$ for all $1
\leq m < l$. For $m = l$ we have $c_{m - 1} = (a_0, \ldots, a_n -
1)$ and we observe that for all subsets $J \not= \{0, \ldots, n\}$,
the sum $\Sigma_{j \in J} c_{m - 1, j} < l - 1$ implying that for
each such $J$, the corresponding summand in $L^\bullet_{l - 1}$ is
contained in $\T$. On the other hand, for the summand corresponding
to $J = \{0, \ldots, n\}$ we have $\Sigma_{j \in J} c_{m - 1, j} = l
- 1$. Thus $L^\bullet_{l - 1}$ is a direct sum of an object in $\T$
and a shift of $i_* ( \st (l - 1))$ as was desired.

Finally, we must examine the objects, $K^\bullet_0$ and
$K^\bullet_l$. For $K^\bullet_0$ observe that $c_0 = (0, \ldots, 0)$
so that $N_0 = S^{n + 1}$. Also, the element $s_0 = (z_0 y^{a_0},
\ldots, z_n y^{a_n})$. Recall that the map $\pi : \tV \to \C^{n +
1}$ is induced by the ring homomorphism $\rho: \C [x_0, \ldots, x_n]
\to \C [z_0, \ldots, z_n, y]$ which sends $x_i$ to $z_i \cdot
y^{a_i}$. Thus $K^\bullet_0$ is the pullback of the Koszul complex
over $\C [x_0, \ldots, x_n]$ generated by the regular sequence
$\left< x_0, \ldots, x_n \right>$. One thus identifies $K^\bullet_0$
as the pullback via $\pi$ of the skyscraper sheaf on $\C^{n + 1}
\subset X$ at zero. In particular $K^\bullet_0 \in \T$.

To see that $K^\bullet_l = 0$ we observe that $c_l = (a_0, \ldots,
a_n)$ and $s_l = (z_0, \ldots, z_n)$. This shows that $K^\bullet_l$
is the Koszul complex associated to the regular sequence $\left<
z_0, \ldots, z_n \right>$. But this is just the regular sequence of
the irrelevant ideal $I$ in $S$ so, modulo torsion, the complex
$K^\bullet_l$ is exact and thus equivalent to zero in the derived
category.

The induction argument then goes as follows. For $m = 0$, we have
seen that $K^\bullet_0 \in \T$. Assume $K^\bullet_m \in \T$ for $0
\leq m < l - 1$, then this object fits into the distinguished
triangle $K^\bullet_{m + 1} \to K^\bullet_m \to L^\bullet_m$. For
all such $m$, we have seen that $L^\bullet_m \in \T$ implying that
$K^\bullet_{m + 1} \in \T$. By induction, we have $K^\bullet_{l - 1}
\in \T$ and from the above observations $0 = K^\bullet_l \in \T$.
The distinguished triangle $K^\bullet_{l} \to K^\bullet_{l - 1} \to
L^\bullet_{l - 1}$ implies $L^\bullet_{l - 1} \in \T$. But as was
observed, $L^\bullet_{l-1} \simeq A \oplus i_*(\st (l - 1) )[-(n +
1)]$ where $A \in \T$. Thus the distinguished triangle, $A[n + 1]
\to L^\bullet_{l - 1}[n + 1] \to i_*(\st (l - 1 ))$ shows that $i_*
(\st (l - 1 )) \in \T$ verifying the second claim and proving the
theorem.

\end{proof}

\subsection{Koszul Duality for Weighted Projective Blowups}
For the purpose of mirror symmetry, one often finds that the natural
exceptional collection to use is the Koszul dual of one formed from
line bundles (for examples see \cite{AKO}, \cite{SeidVC}). The same
is true in the case of the weighted projective blowup. We will need
to calculate the $Ext$ groups and their compositions of the dual
collection in order to see that the mirror derived directed Fukaya
category forms an equivalent category. For more on Koszul duality,
see \cite{Beil1}.

Given any exceptional collection $E = \left< E_1, \ldots, E_r
\right>$ generating the derived category $\T$ , one can regard the
dual collection in the following way. Homological perturbation
theory asserts that $\T$ is equivalent to the bounded derived
category of graded right modules of the quiver algebra $Q_E$
associated to $E$ if $Q_E$ is formal as a dga (see, for example,
\cite{SeidThom}, \cite{Keller}).
 Recall $$ Q_E = \bigoplus_{1 \leq i \leq j \leq r}
Ext^\bullet (E_i, E_j)$$

The exceptional collection in $\D (grMod_r(Q_E))$ corresponding to
$E$ is simply the collection of projective objects $P_i = e_i \cdot
Q_E$ where $e_i$ is the identity morphism in $Ext^\bullet (E_i,
E_i)$. In this situation, the Koszul dual of this collection
consists of the simple objects $S_i = e_{i} \cdot Q_E \cdot e_{i}$
(neglecting shifts). The collection $E^\prime =\left< S_r , \ldots,
S_1 \right>$ is full and exceptional for the same category and
yields a dual quiver algebra $Q_{E^\prime}$. It is this algebra
which we will compute for the exceptional collection $E = \left<
i_*(\st(0)), \ldots, i_*(\st(l - 2)) \right>$ corresponding to the
piece of $\D (\blo )$  from the weighted blowup $\p(\A)$.

Corollary $1$ gives the structure of $Q_E$ as a dga with zero
differential.  In order to apply duality in the straightforward way
outlined above, it must be shown that $Q_E$ is formal. I.e. $Q_E$
carries an $A^\infty$-structure, unique up to quasi-isomorphism,
which comes from the underlying dga via the Homological Perturbation
Lemma. This structure is quasi-isomorphic to one in which higher
products vanish if $HH^k(Q_E, Q_E)$ has no elements of degree $2 -
k$ for $k > 2$ \cite{SeidThom}. In this case, $Q_E$ is intrinsically
formal and we can pursue the above approach to Koszul duality.

In fact, we can restrict attention to an equivariant version of the
Homological Perturbation Lemma. As was shown in Corollary $1$, $Q_E$
has a natural $\C^*$-action, i.e. the morphisms in $Ext^\bullet
(i_*(\st(k)), i_*(\st(j))) = = R_{k - j} \oplus R_{k - j + 1} [-1]$
have weight in $\Z$ associated with the weights in $R = \C [x_0,
\ldots, x_n]$ and composition respects these weights (we will use
the term weight to distinguish from the degree given by the
grading). It follows the $A^\infty$-structure must also respect this
action and the intrinsic formality argument adapts using equivariant
Hochschild cohomology. Furthermore, we take Hochschild cohomology
over the semi-simple base $\oplus_{i = 0}^{l - 2} \C \cdot
1_{\st(i)} = T$.

\begin{prop} The equivariant version of $HH^k(Q_E, Q_E)$ over $T$
is zero in degree $2 - k$ for $k > 2$. \end{prop}
\begin{proof} Suppose $f: Q_E \otimes \cdots \otimes Q_E \to Q_E$ is
a non-zero cocycle representing an element of degree $d$ in
$HH^k(Q_E, Q_E)$. For $1 \leq i \leq k$, let $r_i \in Ext^\bullet
(i_*(\st(a_i)), i_*(\st(b_i)))$ be elements homogeneous in weight
and degree such that $f(r_1 \otimes \cdots \otimes r_k) = r_{k + 1}
\not= 0$ and define $n_i = b_i - a_i$. Then, as we are working over
$T$, we must have $\sum_{i = 1}^k n_i = n_{k + 1}$ and as $f$ is
equivariant with respect to weight we have $\sum_{i = 1}^k
\hbox{wt}(r_i) = \hbox{wt}(r_{k + 1})$. By the description of
$Ext^\bullet (i_*(\st(a_i)), i_*(\st(b_i)))$ in Corollary $1$ we
have that $n_i + \deg (r_i ) = \hbox{wt}(r_i)$ which implies
$\sum_{i = 1}^k \deg (r_i) = \deg(r_{k + 1})$. Thus $d = \deg (f) =
\deg(r_{k + 1}) - \sum_{i = 1}^k \deg (r_i) = 0$ confirming the
claim.
\end{proof}

Now we will apply the above strategy and represent the dual quiver
algebra to $Q_E$. To compute the $Ext$ groups between the simple
modules of $Q_E$, we must find a projective resolution for each such
object. In our case, the additional graded factor in the $Ext$
groups of $E$ poses a technical challenge in this computation. This
is overcome by constructing a double complex of projective objects
whose total complex resolves $S_i$ and whose spectral sequence
converges at the third stage.

We start by resolving $S_i$ by modules over $Q_E$ which are modules
over the degree zero piece of $Q_E$. By Corollary $1$, we have that
$Q_E = Q \oplus I$ where $Q$ consists of the morphisms graded at
zero and $I$ is the ideal whose grading is $1$. Given a graded
module $M$ over $Q_E$, we let ${\tilde M} = M / MI$. Again by
Corollary $1$, we have
$${\tilde P_k} = \bigoplus_{ j = 0}^k Ext^\bullet (i_*(\st (j)),
i_*(\st (k)))
 \approx \bigoplus_{j = 0}^k R_j$$ Recall that in $R$, given an
$\alpha \in \Z_{\geq 0}^{n + 1}$, one has $\deg (x^\alpha ) =
\Sigma_{i = 0}^n \alpha_i \cdot a_i$ where $\A = (a_0, \ldots, a_n)$
is the weight for $\p (\A)$. Thus ${\tilde P_k}$ is a vector space
generated by monomials $\{x^\alpha | \deg (x^\alpha ) \leq k\}$. One
can also see that the $Q_E$-module morphisms $Hom ( {\tilde P_j},
{\tilde P_k})$ are generated by maps sending $e_j$ to $x^\beta$ with
$\deg (x^\beta ) = k - j$. We also note that $S_k = e_k \cdot Q_E
\cdot e_k$ is the quotient module of $\tilde P_k$ by the submodule
generated by all non-trivial monomials. To resolve $S_k$ by modules
$\tilde P_j$, we will exploit the interplay between these modules
and sheaves on $\p (\A )$.

Now, let $[n] = \{0, \ldots, n\}$ and for any subset $J \subset [n]$
let $a_J = \Sigma_{i \in J} a_i$. We have the following exact Koszul
resolution $K^\bullet $ of sheaves on $\p (\A )$ generated by the
regular sequence $\left< x_0, \ldots, x_n \right>$ in $R$
 \begin{displaymath}  0 \to \oplus_{|J| = n + 1} \st ( - a_J ) \to   \oplus_{|J| = n} \st ( - a_J ) \to \cdots
  \to \oplus_{|J| = 1} \st ( -a_J) \to \st  \to 0
 \end{displaymath}
Here we will place $\st$ in the zeroth position. One sees that this
sequence remains exact after tensoring with any line bundle $\st
(i)$ so that the hypercohomology $\h^\bullet (K^\bullet \otimes \st
(i)) = 0$. Furthermore, given $1 \leq i $, $0 < r$ and any $j$,
Proposition $1$ tells us that $H^r (\p (\A), K^j \otimes \st (i)) =
0$ . Thus, for $i \geq 0$, the first term in the spectral sequence
\begin{displaymath} E_{p, q}^1 =  H^q (\p (\A), K^p \otimes
\st(i))) \Rightarrow \h^\bullet (K^\bullet \otimes \st(i)) = 0
\end{displaymath}
is concentrated on the $q = 0$ axis and yields an exact sequence. On
the other hand, the cohomology of $H^0(\p(\A ), K^\bullet)$ itself
is clearly just $\C$ concentrated in degree zero.

As a convention, we will take $P_i = 0$ and $R_i = 0$ if $i < 0$.
Define
$${M_{k, j}} = \bigoplus_{J \subset [n], |J| = j} {
P_{k - a_J}}$$ With these preliminaries in mind, we can state:

\begin{lem} The sequence $H^0(\p (\A ), K^\bullet \otimes \left(
\oplus_{j = 0}^k \st (j) \right) )$ is naturally identified with a
resolution of $S_k$ by modules $\tilde M_{k, j}$.\end{lem}
\begin{proof} By the above definitions and Proposition $1$ we have that \begin{displaymath} \begin{array}{ccc} H^0(\p (\A ), K^i \otimes \left(
\oplus_{j = 0}^k \st (j) \right) ) & = & \oplus_{j = 0}^k H^0(\p (\A
), K^i \otimes \st (j)  ) \\  & = & \oplus_{j = 0}^k H^0(\p (\A ),
\oplus_{|J| = -i} \st(-a_J) \otimes \st (j)  ) \\ & = & \oplus_{j =
0}^k \oplus_{|J| = -i} H^0(\p (\A ), \st (j - a_J)  ) \\ & = &
\oplus_{|J| = -i} \left( \oplus_{j = 0}^k R_{j - a_J} \right) \\ & =
& \oplus_{|J| = - i} {\tilde P_{k - a_J}} \\ & = & {\tilde M_{k,
-i}}
\end{array} \end{displaymath} Furthermore, all maps in the Koszul
sequence $K^\bullet$ are induced from multiplication by $\pm x_i$
and this descends to cohomology. As was pointed out, these are
precisely the morphisms between the modules $\tilde P_i$, thus the
maps in the sequence are $Q_E$-module morphisms. Examining the last
map $$ \oplus_{i = 0}^n {\tilde P_{k - a_i}} \to {\tilde P_k} $$ we
see that this is simply the map that takes $e_i$ in each summand
where $k \geq a_i$ to $x_i$ in ${\tilde P_k}$. So the cokernel of
this map is just $S_k$. Rewriting the sequence in the language of
$Q_E$-modules and observing that for all $0 \leq k \leq l - 2$,
${\tilde M_{k, n + 1}} = 0$ we have the sequence $C_k^\bullet$
\begin{displaymath} 0 \to {\tilde M_{k, n}} \to {\tilde M_{k, n - 1}} \to \cdots
  \to {\tilde M_{k , 1}} \to {\tilde M_{k, 0}}  \to 0
 \end{displaymath} which is quasi-isomorphic to $S_k$ concentrated
 in degree zero.
\end{proof}

We now find truly projective resolutions for the modules $\tilde
M_{k, j}$ and extend $C_k^\bullet$ to a double complex whose total
complex resolves $S_k$. For this we start by resolving $\tilde P_k$
by projective modules. This is actually quite simple. By Corollary
$1$, we have
\begin{displaymath} \begin{array}{ccc} {P_k} & = & \bigoplus_{ j =
0}^k \left( Ext^\bullet (i_*(\st (j)), i_*(\st (k))) \oplus
Ext^\bullet (i_*(\st (j + 1)), i_*(\st (k)))[-1] \right) \\ &
 \approx & \bigoplus_{j = 0}^k \left( R_j  \oplus R_{j-1}[-1] \right)
 \end{array} \end{displaymath}
The summand of $P_k$ in degree $1$ is generated by the identity
element in $R_0[-1]$ which we will denote $e_k[-1]$
. There is a degree $0$ morphism $t_{k - 1}: P_{k -1 }[-1] \to P_k$
which sends $e_{k-1}$ to $e_k[-1]$. We see that the kernel of this
morphism consists of all degree $2$ elements of $P_{k-1}[-1]$ and
the image consists of all degree $1$ elements of $P_k$. Thus the
cokernel is $\tilde P_k$ and one sees easily that the following is a
projective resolution of $\tilde P_k$:
\begin{displaymath} 0 \to P_0[-k] {\buildrel t_0 \over \longrightarrow}
P_1[-k + 1] {\buildrel t_1 \over \longrightarrow} \cdots {\buildrel
t_{k - 1} \over \longrightarrow} P_k \to {\tilde P_k} \to 0
\end{displaymath}
Taking direct sums of these sequences yields a resolution
\begin{displaymath} 0 \to M_{0, j}[-k] { \longrightarrow}
M_{1, j}[-k +1] {\longrightarrow} \cdots {\longrightarrow} M_{k, j}
\to {\tilde M_{k,j}} \to 0
\end{displaymath} Using these resolutions, we extend $C_k^\bullet$ to
a third quadrant double complex $C_k^{\bullet \bullet}$ with
\begin{displaymath} C_k^{p, q} =  \left\{ \begin{array}{cc} M_{k + q , - p}[q] & for \hbox{ }
q \leq 0 \\ 0 & otherwise \end{array} \right. \end{displaymath}
Using the spectral sequence for $C_k^{\bullet \bullet}$ one sees
that the total complex is a projective resolution for $S_k$. Define
the $Q_E$ module \begin{equation} N_{k, j} = \bigoplus_{p + q = -j}
C^{p, q} = \bigoplus_{i = 0}^j M_{k - i, j - i}[-i] = \bigoplus_{J
\subset [n], |J| \leq j} P_{k - j + |J| - a_J}[|J| - j]
\end{equation} Then we have proved
\begin{prop} $S_k$ has a projective resolution $D_k^\bullet$ \begin{displaymath}
0 \to N_{k, n} { \longrightarrow} N_{k, n-1} {\longrightarrow}
\cdots {\longrightarrow} N_{k, 0} \to S_k \to 0
\end{displaymath}
\end{prop}
Now, one easily sees that $Hom (P_i, S_j) = \delta_{ij} \cdot \C$.
Indeed if $i \not= j$ and $f \in Hom (P_i, S_j)$ then $0 = f(e_i
\cdot e_j) = f(e_i) \cdot e_j = f(e_i)$. Using this and the above
resolution, we can calculate the morphisms from $S_k$ to $S_i$. Let
$V$ be a vector space over $\C$ with basis $\{e_0, \ldots, e_n\}$.
On the exterior algebra $\bigwedge^\bullet V$  we have two gradings
generated by $\deg (e_i) = 1$ and $\deg_\A (e_i) = a_i$. Thus we
denote $\bigwedge^{r, s} V = \{w \in \bigwedge^\bullet V | \deg (w)
= r , \deg_\A (w) = s \}$.
\begin{prop} There is a natural equivalence $$Ext^r (S_k, S_i
) = \bigoplus_{s \leq k - i} {\bigwedge }^{r, s} V $$ which is
compatible with composition. \end{prop}
\begin{proof} To calculate $Ext^\bullet (S_k, S_i)$ we first examine
the resolution $D_k^\bullet$. Given a $J \subset [n]$ one sees that
$P_{k - j + |J| - a_J}[|J| - j]$ occurs as a summand in
$D_k^\bullet$ for some $j$ if and only if $k \geq a_J$. Indeed, as
$j \geq |J|$ in each such summand, if $k < a_J$ then $k - j + |J| -
a_J < 0$ implying the associated module is zero. Conversely, if $k
\geq a_J$ then for each $j$ satisfying $|J| \leq j \leq k - a_J +
|J|$ one such summand occurs and these modules are $\{P_0[a_J - k],
\ldots , P_{k - a_J}[0]\}$. For such a $J$, $P_i$ occurs in this
list if and only if $ i \leq k - a_J$. Thus, given any $J \subset
[n]$, the projective module $P_i$ occurs once as a direct summand in
$D^\bullet_k$ of the form $P_{k - j + |J| - a_J}[|J|-j]$ if and only
if $a_J \leq k - i$. Such a summand occurs in $D_k^{-j}$ so that, as
a graded module, the total grading on the summand is $P_i[|J| -
j][j] = P_i[|J|]$. Now, as each map in $D^\bullet_k$ restricted to
any summand either sends $e_p$ to $e_p \cdot x_q$ or to $e_p[-1]$,
we have that the complex $Hom (D_k^\bullet, S_i)$ has a zero
differential. This implies that $Ext^\bullet (S_k, S_i) = Hom
(D_k^\bullet , S_i)$ and recalling the above note that $Hom (P_j,
S_i) = \delta_{ij} \C$ we have
\begin{displaymath}
\begin{array}{ccc} Ext^\bullet (S_k , S_i) & = & \bigoplus_{J \subset [n], a_J \leq k - i}
Hom (P_i[|J|], S_i) \\ & = & \bigoplus_{J \subset [n], a_J \leq k -
i} \C [-|J|]
 \end{array} \end{displaymath} The identification is now clear. For
 any $J = \{j_1, \ldots, j_{|J|}\}$ with $a_J \leq k - i$ assign the
 projection map in $Hom (P_i[|J|], S_i) = \C [-|J|]$ to the
 element $v_J = e_{j_1}\wedge \cdots \wedge e_{j_{|J|}} \in \bigwedge
 V$. Then we see that $\deg (v_J)= |J|$
 so that the grading of $v_J \in \C [-|J|]$ is $ \deg (v_J)$ as claimed.

 To see that these isomorphisms are compatible with composition, one observes that any
 element of $Ext^c (S_k, S_i)$ is a morphism from $D_k^\bullet$ to
 $D_i^{\bullet + c}$. Since the horizontal differential of $C^{\bullet \bullet}$
 takes identity elements to degree $0$ elements and the vertical differential takes identities to elements of degree $1$, such maps naturally extend to
 maps of the double complexes giving each such map a bi-degree.
 Working out the degrees, one sees that the element corresponding to
 $v_J$ has bi-degree $( - \deg (v_J), \deg_\A (v_J) + i - k)$. In
 particular, for $J = \emptyset$ we have the map which is simply a
 vertical shift plus projection of $D_k^{\bullet \bullet}$ to
$ D_i^{\bullet \bullet}$. One easily sees that this implies that
$v_\emptyset$ acts as an identity element under compositions (after
one has made the identification between the $Ext$ groups and the
exterior algebra). Thus in examining relations, one needs only
consider elements in $Ext^\bullet (S_k, S_i)$ of bi-degree $(d, 0)$
which correspond to elements for which $a_J = k - i$. As all such
elements are of degree zero, one can quotient the complex by the
weight $-1$ ideal in $Q_E$. Viewing the double complex over $Q$
gives a zero vertical differential, so examining the relations with
respect to the complexes $C^\bullet$ over $Q$ is sufficient. These
relations have been worked out in \cite{AKO}, section 2.6, where it
was found that for elements of bi-degree $(d, 0)$, there is a
natural equivalence
$$Ext^\bullet_{Q} (S_k, S_i) = \bigwedge^{r, k - i} V [-r]$$
compatible with composition.
 \end{proof}

 {\it Note}: There are two alternative and less computational ways of approaching this problem.
  One could approach the above proof from the perspective of Koszul algebras. Indeed,
 letting $k =  \oplus_{i = 0}^{l - 2} \C \cdot e_i$ and viewing $M$ as a module generated by
 monomials $x_i: P_k \to P_j$ and $e_i[-1]$ for all $i$, $j$ and
 $k$, $Q_E$ can be viewed as a quadratic algebra $T(M)/ I$ where $I$
 is generated by elements in $M \otimes_k M$. If one could show that
 this algebra is in fact Koszul, then the dual quadratic algebra
 is isomorphic to the opposite algebra of the Koszul dual. This algebra is easily seen to be the one
 exhibited in the previous proposition. For more on quadratic
 algebras over semi-simple rings, see \cite{Beil2}.

 Alternatively, one could view $\T$ as a subcategory of graded modules over
 the super symmetric algebra $T = Sym (V_0 \oplus V_1)$ where the even part $V_0 = V$
 is the weighted vector space as above and and the odd $V_1 = \C$ encodes the infinitesimal information of
 the exceptional divisor. One readily sees that $i_*(\st(k))$ could
 be identified with the module $T_k$ by relating $Hom (T_k, T_j)$ to
 $Ext^\bullet (i_*(\st(k)), i_*(\st(j)))$ via Corollary $1$. Then
 the Koszul dual collection is found by examining the dual collection
 to the dual algebra $T^\vee = Sym (V[1])$ which gives a Grassmanian in $n + 1$ variables
 tensored with a symmetric algebra in one variable. One can see that
 $S_k$ can be associated to $T^\vee_{-k}$ by rewriting Proposition
 $5$ in terms of homogeneous elements of $T^\vee$.

 \section{Degenerations of Lefschetz Pencils}
\subsection{The Directed Fukaya Category}

This preliminary subsection will outline the definition of the
directed Fukaya category corresponding to the Landau-Ginzburg model.
There will be no proofs, but all of the stated results can be found
in the literature, see for example \cite{SeidVC}, \cite{AKO},
\cite{Fukaya}. We will impose some technical assumptions which
greatly simplify the construction. These will always be satisfied in
our applications. First, assume $(M, \omega )$ is an exact K\"ahler
manifold. The condition of exactness is that there exists a $1$-form
$\theta $ such that $\omega = d \theta $. Now let $W: M \to \C$ be a
holomorphic function with isolated Morse singularities. We will
assume that there is at most one critical point in any given fiber
of $W$ and will call $W$ the potential. Given such a potential
function and a regular point $p \in M$, one can use the symplectic
form to define a splitting of $T_p M = T_p (W^{-1}(p)) \oplus U_p$.
Here, $U_p$ is the symplectic orthogonal to the tangent space of the
fiber. One sees that $U_p$ is mapped isomorphically onto the
$T_{W(p)} \C$ via the differential of $W$. So, given a tangent
vector at $W(p)$, one can lift this to a vector field on $W^{-1}(p)$
in $M$ (at a critical point on a fiber, we take the vector field to
equal zero).

Now suppose $\{q_0, \ldots, q_m\}$ are the critical values of $W$
and choose a regular value $q$. Let $\gamma_i :[0, 1] \to \C$ be a
set of paths such that $\gamma_i(0) = q$, $\gamma_i(1) = q_i$,
$\gamma_i((0, 1]) \cap \gamma_j ((0, 1]) = \emptyset$ and the paths
are oriented counter-clockwise around $q$. We will refer to such a
collection $\{\gamma_i\}$ as a distinguished basis. Along each such
path, one can lift the vector field $d/dt$ to a vector field on
$W^{-1}(\gamma_i ([0, 1])$ which will define a parallel transport.
We will assume for the moment that this vector field is integrable
on $M$ and let $D_i$ be the set of points flowing into $\tilde q_i$,
the critical point which is mapped to $q_i$. We let $L_i = D_i \cap
W^{-1}(q)$ be its boundary. $D_i$ is often called a Lefschetz
thimble and $L_i$ the vanishing cycle. One can show that the
vanishing cycles $L_i$ are actually exact lagrangian sub-manifolds
in the fiber $W^{-1}(q)$. Such a submanifold is one for which
$\theta|_{L_i}$ is an exact one form on $L_i$. Assume that the
vanishing cycles intersect transversely (indeed, one can always
perturb them slightly to accomplish transversality). As the
potential $W$ is not proper, the parallel transport vector field may
not be integrable in general. To get around this one must
occasionally perturb the parallel transport vector field. So long as
such a perturbation results in an exact isotopy of the $L_i$, the
theory will remain unaffected. For more on vanishing cycles, see
\cite{Donaldson}, \cite{Milnor}.

Now, one can form the Lagrangian Grassmanian bundle $\Lambda$ of the
tangent bundle $T (W^{-1}(q))$ on the fiber $W^{-1}(q)$. Assuming
certain obstructions vanish, we also have the fiberwise universal
covering $\eta : {\tilde \Lambda} \to \Lambda$. For each vanishing
cycle $L_i$, one has a natural lift
\begin{displaymath}
\begin{array}{ccc} L_i & {\buildrel {\phi }
\over \dashrightarrow} & \Lambda \\ & \searrow & \downarrow \\  & &
W^{-1}(q)
\end{array} \end{displaymath}
where $\phi$ sends a point of $L_i$ to its tangent space. A grading
on $L_i$ is a lift ${\tilde \phi}:L_i \to {\tilde \Lambda}$
satisfying $\eta \circ {\tilde \phi} = \phi$. As in \cite{Seidbook},
we define a graded lagrangian to be a lagrangian submanifold with
such a lift $\tilde \phi$. For more details on gradings see Section
$4.4$.

The objects of the directed Fukaya category for $(W ,\{\gamma_i \})$
are graded vanishing cycles $[L_i, {\tilde \phi}_i]$ where we have
chosen one grading for each vanishing cycle. Given a point $p \in
L_i \cap L_j$ one can choose a path $\delta_p:[0, 1] \to \tilde
\Lambda_p$ defined by setting $\delta_p (0) = {\tilde \phi_i}(p)$
and $\delta_p (1) = {\tilde \phi_i}(p)$; as this path is in the
simply connected universal cover of the lagrangian grassmanian, it
is unique up to homotopy. Taking the Maslov index of $\eta \circ
\delta_p$ gives an integer which we will call $\deg (p)$. We can now
define the Hom sets with their associated grading.

\begin{displaymath} Hom^\bullet (L_i, L_j) = \left\{ \begin{array}{cc} \oplus_{p \in L_i \cap L_j} \C \left< p \right> [- \deg(p)] & for \hbox{ }
i < j  \\ \C \left< e_i \right> & for \hbox{ } i = j
\\ 0 & otherwise  \end{array} \right.
\end{displaymath}

In order to proceed to products, one must confront the fact that the
directed Fukaya category is actually an $A^\infty$-category. Thus,
for every $k \in \Z_{>0}$ one has a higher product $m_k$ of degree
$2 - k$. To define these maps, one considers the moduli space of
holomorphic maps from the disc with marked boundary points to the
fiber $W^{-1}(q)$ with certain boundary conditions. More
specifically,  Let $\M_{k} (D)$ be the moduli space of the unit disc
$(D,
\partial D)$ with $k + 1$ distinct marked points $\{z_0, \ldots, z_{k}\}$ oriented clockwise
along the boundary. Now, assume $i_0 < i_2 < \cdots < i_{k}$ and
$p_{i_j} \in L_{i_j} \cap L_{i_{j + 1}}$ then we take $\M (p_{i_0},
\ldots, p_{i_k})$ to be the moduli space of all holomorphic maps $u:
D_{\bf z} \to W^{-1}(q)$ such that $D_{\bf z} \in \M_k (D)$, $u(z_j)
= p_{i_j}$ and the arc along the boundary from $z_j$ to $z_{j + 1}$
is sent to $L_{i_{j + 1}}$. In the above notation we take $j$ mod
$k+1$ and when $k = 1$, we need to quotient this space by the action
of $\R$ on such maps. In the situation described above, this moduli
space has a natural compactification which is a manifold with
corners. We will denote $\M (p_{i_0}, \ldots, p_{i_k})$ by $\M_0
(p_{i_0}, \ldots, p_{i_k})$ if it has dimension zero, otherwise
$\M_0 (p_{i_0}, \ldots, p_{i_k})$ will be the empty set. This space
has a natural orientation class $sgn$ and defines our products.
Namely, using the above notation one defines
$$m_k: Hom^\bullet (L_{i_{k - 1}}, L_{i_{k}}) \otimes \cdots \otimes
Hom^\bullet (L_{i_{0}}, L_{i_{1}}) \to Hom^\bullet (L_{i_{0}},
L_{i_{k}})[2 - k]$$ via
$$m_k (p_{i_{k - 1}} \otimes \cdots \otimes p_{i_0}) = \Sigma_{r \in
L_{i_0} \cap L_{i_k}} \left( \Sigma_{u \in \M_0 (p_{i_0}, \ldots,
p_{i_{k - 1}}, r)} sgn (u) \right) r$$ A Maslov index calculation
shows that $m_k$ has the indicated degree. As was noted in
\cite{AKO}, one normally weights this sum by the exponential of the
symplectic area of $u(D)$ in $W^{-1}(q)$. However, since the
lagrangians are exact, this is not necessary.

The version of homological mirror symmetry examined in this paper
proposes an equivalence between the derived category of coherent
sheaves on a toric variety and the derived Fukaya category of the
mirror Landau-Ginzburg model. To pass from the $A^\infty$-category
described above to its derived category, one must take twisted
complexes of formal sums and shifts of the above objects, their
idempotent splittings and formal inverses of quasi-isomorphisms.
This procedure was invented by Kontsevich \cite{Kontsevich} and is
fully explained in Seidel's recent book \cite{Seidbook}. After
forming such a construction, the objects $[L_i, {\tilde \phi_i}]$
form a complete exceptional collection in the derived category $\D
(Fuk (W, \{\gamma_i \})$. A result of Seidel is that any other
choice of paths $\{\gamma_i^\prime \}$ yields objects $[L_i^\prime ,
{\tilde \phi_i}^\prime ]$ that can be obtained from $[L_i, {\tilde
\phi_i}]$ in $\D (Fuk(W, \{\gamma_i \})$ by a sequence of mutations
\cite{SeidVC}. As a consequence, the derived Fukaya category is
independent of the choice of paths and one can simply write $\D (Fuk
(W))$. We will need the following theorem for what follows which is
a standard fact whose proof can be found in Seidel's book in the
section on directed $A^\infty$-categories \cite{Seidbook}.

\begin{thm} (Seidel) $\D (Fuk (W, \{\gamma_i \})$ is invariant under
exact perturbations of $\omega$. Furthermore, if $[L_i^\prime ,
{\tilde \phi_i}^\prime ]$ are Hamiltonian isotopic to $[L_i, {\tilde
\phi_i}]$, then $\D (Fuk (W, \{\gamma_i \})$ is equivalent to the
derived Fukaya category generated by the graded lagrangians
$[L_i^\prime , {\tilde \phi_i}^\prime ]$. \end{thm}

This theorem is not only true on the derived level, but also up to
$A^\infty$-quasi-isomorphism. By an exact perturbation we mean a
perturbation of $\theta$ through one forms whose exterior derivative
is non-degenerate.

In our case we will see that all $m_k$ vanish except $m_2$. This
implies that the $A^\infty$-category is indeed a category. In such a
situation, it is known that $\D (Fuk (W, \{\gamma_i \})$ is
equivalent as a triangulated category to the bounded derived
category of graded modules over the quiver algebra defined by the
collection $[L_i, {\tilde \phi_i}]$. All of this will be explored in
more detail in Section 3. \\[4pt]

{\it Note}: In other accounts of this subject the distinguished
basis of paths $\{\gamma_i\}$ is ordered clockwise and the moduli
space $\M (D)$ has curves with marked points oriented
counterclockwise. Our formulation is equivalent and more
advantageous for the examples worked through in Section 3; however,
as will be seen, it will affect the definition of the Maslov index
as well.

\subsection{Partial Lefschetz Fibrations} Let  $W$ be a potential function
from a K\"ahler manifold $M$ to $\C$. In many cases of interest
symplectic parallel transport is not well defined with respect to
$W$ for all paths in $\C$. These cases arise when certain fibers of
$W$ do not transversely intersect the divisor at infinity of a
suitable compactification of $M$. This fact can make it difficult to
define the directed Fukaya category for the pair $(M, W)$. However,
when confronted with such a potential $W$, one can often find a
connected open set $U \subset \C$ for which parallel transport can
be defined. In the ideal case, all critical values of $W$ are
contained in $U$ and $U$ is simply connected; this case yields the
usual concept of a Landau-Ginzburg model explored in the last
subsection. In the less ideal case, one can capture part of the
Landau-Ginzburg model by examining non-simply connected domains $U$
containing only some of the critical values of $W$. The motivation
for this procedure is that this part of the model will remain
invariant under perturbation of $W$ yielding a semi-orthogonal piece
of the directed Fukaya category corresponding to the perturbed
model.

We now fix a Stein manifold $M$ with an exhaustive plurisubharmonic
function $\rho$ and the associated (exact) symplectic structure
$\omega$. We also take a potential function $W$ from $M$ to $\C$
associated to the Landau-Ginzburg model on $M$.

\begin{defn} A partial Lefschetz fibration on $(M, \rho )$ consists of
the data $(W, U, a)$ where $W$ is a holomorphic function on $M$, $U$
is a connected open subset of $\C$ with smooth boundary and $a \in
\R$ such that
\\ (i) $a$ is a regular value of $\rho$
\\ (ii) $W^{-1}(p)$ intersects $\rho^{-1}(a)$ transversely for each $p
\in {\overline U}$  \\ (iii) The critical points $p_i$ such that
$W(p_i) \in U$ are
Morse and no two such points lie on the same fiber \\ (iv) The set of critical values of $W$ does not intersect $\partial U$ \\
(v) $\rho (p) < a$ for all critical points $p \in M$ with $W(p) \in U$  \\
Given a partial Lefschetz fibration, let $M_{(U, a)} =
W^{-1}({\overline U}) \cap \rho^{-1}((-\infty, a])$.
\end{defn}

The word "partial" refers to the part of $M$ corresponding to
$M_{(U, a)}$ which is a K\"ahler manifold with boundary and
codimension $2$ corners. On this part of $M$, the fibers of $W$ are
K\"ahler manifolds with contact type boundary and the setup looks
similar to the one considered in \cite{SeidVC}. A major flaw however
is the possible failure of $\rho^{-1}(a_i) \cap W^{-1}({\overline
U})$ being horizontal with respect to the symplectic orthogonal
connection on $M$. In other words, attempting parallel transport of
a fiber $W$ in $M_{(U, a)}$ along a path in $U$ may not be possible
as one may be led outside the subspace $M_{(U, a)}$. On can remedy
this situation by adding the Liouville vector field to any given
parallel transport \cite{SeidSmith}. While flowing downward along
the Liouville vector field does not preserve the symplectic
structure, it will be an exact isotopy for exact Lagrangian
submanifolds. These are the manifolds of interest when defining the
directed Fukaya category.

Let us start to pursue this line of reasoning by recalling some
basic definitions and results. The symplectic form on $M$ is $\omega
= -d (J^* d \rho )$ where $J$ is the complex structure on $M$. In
the notation of the previous subsection $\theta = - J^* d \rho$. The
Liouville vector field $X$ is defined by the equality
$$\omega (X , Y) = \theta (Y)$$ for all vector fields $Y$. We have
the classical result that $L_X \theta = \theta$.

Alternatively, one can view $X$ as the gradiant vector field of
$\rho$ with respect to the K\"ahler metric. In particular, $X$ is
normal to the real hypersurface $\rho^{-1}(a)$ for all regular
values $a$. Now suppose $W$ is a potential function on $M$. On each
fiber $F_p$, the restriction of the K\"ahler form $\omega$ to $F_p$
is generated by the restriction of the K\"ahler potential $\rho$ to
the fiber. Thus, for each fiber we can define a Liouville vector
field which taken together on all of $M$ we will call the fiberwise
Liouville vector field $X_f$. Although one no longer has $L_{X_f}
\theta = \theta$ for $M$, this equality is certainly true fiberwise.
One can see also that as $X_f$ is the gradiant vector field on $F_p$
to $\rho$ restricted to the fiber, it will be normal to
$\rho^{-1}(a)$ in $M$ for those points at which $F_p$ transversely
intersects $\rho^{-1}(a)$.

Now assume $(W, U, a)$ is a partial Lefschetz fibration with
critical values $\{q_0, \ldots, q_m\} \subset U$ and corresponding
critical points $\{p_0, \ldots, p_m \}$. Choose a regular point $q
\in U$ of $W$ and choose a distinguished basis of paths $\{
\delta_0, \ldots, \delta_m \}$ from $[0,1]$ to $U$ satisfying the
conditions of the previous subsection. We will call such a choice of
paths a distinguished basis for the partial Lefschetz fibration $(W,
U, a)$. For each such path we have a parallel transport vector field
$Y_i$ on $W^{-1}(\delta_i([0,1]))$. As was discussed above, for $p
\in \rho^{-1}(a)$, one may have that $Y_i$ at $p$ enters the space
$M_{(U,a)}$, i.e. $d\rho (Y) < 0$. However, there is a sufficiently
large constant $C_i$ such that $d\rho (Y_i + C_i X_f) > 0$ for every
point in the compact space $\rho^{-1}(a) \cap M_{(U, a)}$ (this
follows from the assumption of transversality in the definition of
partial Lefschetz fibrations). It is necessary to add the fiberwise
Liouville vector field $C_i X_f$ to $Y_i$ as opposed to the
Liouville vector field as $X$ is not generally defined on the
tangent space of $W^{-1}(\delta_i([0,1]))$, i.e. certain tangent
vectors may lead one off of this subspace. Let $C$ be the maximum of
such $C_i$ and $Z_i = Y_i + C X_f$. We let $D_i^{lv}$ be the space
of all points in $W^{-1}(\delta_i([0,1]))$ which flow to the
critical point $p_i$ via $Z_i$ and $L_i^{lv} =
\partial D_i^{lv} \subset W^{-1}(q)$. Let $D_i$ and $L_i$ be defined
as in the previous subsection. For the following proposition, we
will assume that the vector fields $Y_i$ are integrable.

\begin{prop} For any partial Lefschetz fibration $(W, U, a)$ on $M$, one has $D^{lv}_i
\subset M_{(U, a)}$ for every $i$. Furthermore, there is an exact
isotopy from $L_i^{lv}$ to $L_i$ in $W^{-1}(q)$ for each $i$.

\end{prop}
\begin{proof} To see that $D^{lv}_i \subset M_{(U, a)}$, observe
that if $f:[0, \infty) \to W^{-1}(\delta_i([0,1]))$ is a flow line
for which $\lim_{t \to \infty} f(t) = p_i$ and $ \rho( f(0)) \geq a$
then as $p_i \in M_{(U, a)}$ there is some $p \in \rho^{-1}(a) \cap
M_{(U, a)}$ and some $t_0$ for which $f(t_0) = p$ and $d/dt (\rho
\circ f)_{t_0} \leq 0$. But then $$0 \geq d/dt (\rho \circ f)_{t_0}
= d \rho ( (df/dt)_{t_0}) = d \rho (Z_i(p)) > 0$$ Therefore, $\rho
(f(0)) < a$ and $f(0) \in M_{(U, a)}$.

The second assertion follows from the fact that $L_{X_f} \theta =
\theta$ on the fibers of $W$. Thus, flowing along $X_f$ simply
rescales $\theta|_{L_i}$ so that exactness is never violated.
\end{proof}

Now, given the collection $\{L_0^{lv}, \ldots, L_m^{lv} \}$ in
$W^{-1}(q)$ we can define a directed $A^\infty$-category following
the procedure outlined in the last subsection where all products are
defined using moduli spaces $\M_0 (p_{i_0}, \ldots, p_{i_k})$ where
the target is $W^{-1}(q)$. As $\rho$ is a subharmonic function on
$W^{-1}(q)$, the maximum principle shows that any holomorphic map
from the unit disc into $W^{-1}(q)$ such that $\rho(\partial D)
\subset M_{(U, a)}$ must have $\rho (D) \subset M_{(U, a)}$. Thus
the space $\M_0 (p_{i_0}, \ldots, p_{i_k})$ only contains
holomorphic maps whose image lies in the subspace $M_{(U, a)}$. With
this in mind we define
\begin{defn} Let $(W, U, \rho)$ be a partial Lefschetz fibration and
$\{\delta_0, \ldots, \delta_m\}$ be paths as above. Define $\fp$ to
be the $A^\infty$-category generated by the lagrangians $\{L_0^{lv},
\ldots, L_m^{lv} \}$ with higher products defined in 3.1 \end{defn}

{\it Note}: Unlike the situation for ordinary Landau-Ginzburg
models, the derived version of this category will depend on the
choice of paths. This situation arises when one takes $U$ to be a
non-simply connected domain. Then the monodromy prevents one from
concluding that different choices correspond to mutations in the
derived category.

From the arguments given above and Theorem $3$ we have

\begin{cor} Let $W$ be a Lefschetz fibration on $M$. If $(W, U, a)$
is a partial Lefschetz fibration and $\{\delta_i\}$ a distinguished
basis for $(W, U, a)$. Adding paths $\{\tau_j\}$ to obtain a
distinguished basis for $W$, let $Fuk(W|_U)$ be the
$A^\infty$-category contained in $Fuk(W)$ generated by the vanishing
cycles associated to $\{\delta_i\}$. Then $\D (Fuk(W|_U))$ is
equivalent to $\D (\fp)$.
\end{cor}

This corollary shows that $\D (\fp )$ can be considered as a
subcategory of $\D (Fuk (W))$ when $W$ is a Lefschetz fibration.
Furthermore, if there is a distinguished basis of $W$ of the form
$\{\tau_1, \ldots, \tau_r, \delta_0, \ldots, \delta_m\}$ we see that
 the triangulated category generated by the first $r$ vanishing
cycles and $\D (\fp)$ yields a semiorthogonal decomposition of $\D
(Fuk (W))$.

\subsection{Deformations of the potential}

The motivation behind the preceding definitions is not simply to
take a piece of a single Landau-Ginzburg model $(M, W)$, but rather
to take a piece that will remain somewhat stable under a
perturbation of $W$. Of course, slight perturbations of the
potential $W$ may introduce new critical points and new topology
into the fibers. The simple example of $W_\epsilon (z) = \epsilon
z^3 + z^2$ on $\C$ demonstrates such a situation. For $\epsilon =
0$, we have one critical point and the regular fibers consist of $2$
points. After slightly perturbing, a critical point enters in from
infinity, and the generic fibers become $3$ points. This basic
example shows that the full derived Fukaya category of $M$ with a
varying potential $W_\epsilon$ is not generally invariant. In this
paper, we will be considering such unruly types of perturbations, so
it will be desirable to pick out a piece of the Fukaya category that
does remain invariant. Of course, this piece will be precisely the
one corresponding to a partial Lefschetz fibration. This motivates
the following definition.

\begin{defn} Suppose $(W, U, a)$ is a partial Lefschetz fibration on
$(M, \rho)$. Letting  $D_\varepsilon$ be the $\varepsilon$ disc in
$\C$, a deformation of $(W, U, a)$ consists of a holomorphic map
$W_{-}: D_\varepsilon \times M \to \C$ such that $W_0 = W$ and
$(W_t, U, a)$ is a partial Lefschetz fibration for all $t \in
D_\varepsilon$. Given such a deformation, we let $M_{(U, a, t)} =
W_t^{-1}(U) \cap \rho^{-1}((-\infty, a])$.
\end{defn}

Observe that a consequence of this definition is that no critical
values or points will enter into $U$ or $M_{(U, a)}$ as we perturb
$W$. Furthermore, the transversality condition insures that
$M_{(U,a)}$ and its fibers are topologically unchanged. Moreover, we
have the following proposition.

\begin{prop} Let $W_{-}$ be a deformation of a partial Lefschetz
fibration $(W, U, a)$. If $q_t \in U$ is a smoothly varying regular
value of $W_t$, then $W_t^{-1}(q_t) \cap M_{(U, a, t)}$ contains a
subspace $N_{(U, a, t)}$ that is an exact perturbation of
$W^{-1}(q_0) \cap M_{(U, a)}$ up to a rescaling of $\omega$.
\end{prop}
\begin{proof} The proof of this theorem follows the spirit of the
proof of Proposition 7. We let $F: D_\varepsilon \times M \to
D_\varepsilon \times \C$ be defined as $F = (\pi_1 , W_{-})$.  Given
a path $\gamma(t) = (\gamma_1(t) , \gamma_2(t)) $ in $D_\varepsilon
\times \C$, we can form a parallel transport vector field in
$D_\varepsilon \times M$ along the fibers $W^{-1}_{\gamma_1
(t)}(\gamma_2(t))$. Although this vector field will generally not be
integrable, adding a fiberwise Liouville flow will force
integrability in the subspaces $M_{(U, a, t)}$. We have that
parallel transport is an exact isotopy and the Liouville flow simply
rescales the exact form $\theta$, implying the result. The details
of this argument follow the same line of reasoning as in Proposition
$7$.
\end{proof}

Now suppose $\{\delta_{i, t}\}$ is a smoothly varying set of
distinguished bases for $W_t$ and $\{L_{i, t}^{lv}\}$ are their
associated vanishing cycles in $W^{-1}_t(q_t)$. Using the previous
proposition, there is an embedding $j_{0, t}$ from $W^{-1}(q_0) \cap
M_{(U, a, 0)}$ into $W^{-1}(q_t) \cap M_{(U, a, t)}$ which is an
exact symplectic deformation. Thus we have the exact Lagrangians,
$\{j_{0, t}(L_{i, 0}^{lv})\}$ contained in $W^{-1}(q_t) \cap M_{(U,
a, t)}$ as well.

\begin{prop} With the notation above, $\{j_{0, t}(L_{i, 0}^{lv})\}$ is
exact isotopic in $M_{(U, a, t)}$ to $\{L_{i, t}^{lv}\}$ for every
$i$.
 \end{prop}
 \begin{proof} This result follows from considering the definitions
 of both $j$ and $L_{i, t}^{lv}$. Both such objects were constructed
 by adding a sufficiently large multiple of the fiberwise Liouville
 vector field and integrating the new parallel transport map.
 Indeed, given a smoothly varying critical value $q_{i, t}$ one can consider the path
 $\gamma (s) = (\gamma_1(s), \gamma_2(s))$ such that $\gamma_1(s) =
 0$ and $\gamma_2(s) = \delta_{i, 0}(2s)$  for $0 \leq s \leq 1/2$.
 For $1/2 \leq s \leq 1$ let $\gamma_1(s)$ be the path from $0$ to
 $t$ and $\gamma_2(s) = q_{\gamma_1(s)}$. One may slightly
 perturb $\gamma$ to make it a smooth map and add a
 sufficiently large multiple of the Liouville vector field to the
 parallel transport field. Observe that integrating the resulting
 field from the critical point $p_{i, o}$ simply gives $j_{0,
 t}(L_{i, 0}^{lv})$. On the other hand, one can deform $\gamma$
 through smooth paths $\gamma_r$ to $\gamma^\prime(s) = (\gamma^\prime_1(s), \gamma^\prime_2(s))$
 such that for $0 \leq s \leq 1 / 2$ one has $\gamma^\prime_1(s) =
 \gamma_1(1/2 + s)$ and $\gamma^\prime_2(s) = \delta_{i,
 \gamma^\prime_1(s)}(0)$. For $1/2 \leq s \leq 1$ let
 $\gamma^\prime_1(s) = t$ and $\gamma^\prime_2(s) = \delta_{i,
 t}(2s - 1)$. One sees that integrating the resulting vector field
 along this path from the point $p_{i, 0}$ yields $L_{i, t}^{lv}$.
 For each path in the deformation from $\gamma$ to $\gamma^\prime$
 we can avoid critical values of $F$ and obtain a smooth family of
 exact lagrangians in $W^{-1}(q_t) \cap M_{(U, a, t)}$. This gives
 us the required isotopy.
\end{proof}

Using this and Theorem $3$, we have the following corollary which
also holds on the level of quasi-isomorphic $A^\infty$-categories.

\begin{cor} Let $W_{-}$ be a deformation of a partial Lefschetz
fibration $(W, U, a)$ and $\{\delta_{i, t}\}$ be a smoothly varying
family of distinguished bases for $(W_t, U, a)$. Then $\D (\fp )$ is
equivalent to $\D (\fpt )$ for every $t \in D_\varepsilon$.
 \end{cor}
\begin{proof} Using Theorem $3$ and Proposition $8$, we have that the derived category
generated by $L_{i, 0}^{lv}$ in $M_{(U, a)}$ is equivalent to that
generated by $j_{(0, t)} (L_{i, 0}^{lv})$ in the image of $N_{(U,
a,t)}$ in $M_{(U, a, t)}$. By the previous proposition, these
lagrangians are exact isotopic to $\{L_{i, t}^{lv}\}$ in $M_{(U, a,
t)}$ so that again by Theorem $3$, their derived categories are
equivalent in $M_{(U, a, t)}$. To be precise, one should use the
proposition which can be found in \cite{Seidbook} under "Properties
of the Fukaya Category" to see that all moduli spaces defining the
higher products are equivalent (i.e. to see that the image of a
holomorphic curve in $M_{(U, a, t)}$ with boundary in $N_{(U, a,
t)}$ actually lies in $N_{(U, a, t)}$).
\end{proof}

Utilizing this corollary and the results of previous subsections, we
come to the main result of this section. To understand the statement
we recall that the Fukaya category has an obvious exceptional
collection corresponding to the ordered distinguished basis of
paths. Given a connected open subset $U$ in $\C$ and a set of
critical values of some potential $W$, we can choose a basis of
paths from a regular point in $U$ to the critical points in $U$ such
that each path is contained in $U$. We can then extend this set to a
distinguished basis of $W$ so that all the paths of the basis which
tend to points outside of $U$ occur before those in our set. Doing
so yields a semi-orthogonal decomposition of the Fukaya category
into two subcategories, the left piece being the subcategory
generated by the vanishing cycles associated to points outside $U$
and the right piece the subcategory generated by vanishing cycles
associated to the critical points in $U$.

\begin{thm} Let $W_{-}$ be a deformation of a partial Lefschetz
fibration $(W, U, a)$ and $\{\delta_{i, t}\}$ be a smoothly varying
family of distinguished basis' for $(W_t, U, a)$. Suppose $W_{t_0}$
is a Lefschetz fibration for some $t_0$. Then there is a
semi-orthogonal decomposition
$$\D (Fuk (W_{t_0})) = \left< \T, \D (\fp )\right>$$ where the
triangulated category $\T$ is generated by the vanishing cycles
associated to critical values of $W_{t_0}$ in $\C - U$.

\end{thm}
\begin{proof} We start by choosing a distinguished basis $\{\gamma_1, \ldots, \gamma_k, \delta_{0,
t} \ldots, \delta_{m, t}\}$ for $W_t$. Here, the paths $\gamma_j$
are from $q_t \in U$ to the critical values of $W_t$ in $\C - U$. We
have that $\D (Fuk (W_t))$ decomposes into the category $\T$
generated by the vanishing cycles associated to the $\gamma_j$ and
into the category generated by the vanishing cycles associated to
$\delta_{i, t}$. By Corollary 2, the latter category is equivalent
to $\fpt$, and by the Corollary 3, this category is equivalent to
$\fp$.

\end{proof}

\subsection{The case $M = (\C^*)^n$}
In this section we will be concerned with Laurent polynomials
$W(z_1, \ldots, z_n) = \sum_{\alpha \in A} c_\alpha z^\alpha$ from
$\gmc$ to $\C$. Here $\alpha = (\alpha_1, \ldots, \alpha_n)$ and
$z^\alpha = z_1^{\alpha_1} \cdots z_n^{\alpha_n}$. Also $A \subset
\Z^n$ is the support of $W$, i.e. $c_\alpha \not= 0$ if and only if
$\alpha \in A$. In certain situations, we will be able to decompose
$A$ as the union of two subsets $A_1$ and $A_2$. Each subset comes
with it's own potential $W_i$ on $\gmc$ which can be taken to be $W$
restricted to $A_i$, i.e. we set all coefficients with subscripts
not in $A_i$ equal to zero. The idea is to apply Theorem $4$ to this
situation and obtain a semi-orthogonal decomposition of $\D (Fuk
(W))$ into the categories $\D (Fuk (W_1))$ and $\D (Fuk (W_2))$. One
of these categories, however, will require use of a partial
Lefschetz fibration, as the associated symplectic vector field will
not be integrable for all paths in $\C$. In carrying out this plan,
we will make heavy use of techniques from \cite{GKZ}.

We start by establishing some notation. Throughout this paper, we
will take $\rho (z_1, \ldots , z_n) = \Sigma_{i = 1}^n {1/2} (\log
|z_i|^2)^2$ to be the K\"ahler potential on $\gmc$. One can see that
any $a \not= 0$ is a regular value for $\rho$ and that
$\rho^{-1}((-\infty, a])$ is compact for all such $a$, i.e. $\rho$
is an exhaustive plurisubharmonic function on $(\C^*)^n$. We also
observe that $\omega = \Sigma_{i = 1}^n |z_i|^{-2} dz_i \wedge
d{\overline z_i}$. Now, given any Laurent polynomial $f(z_1, \ldots,
z_n) = \sum_{\alpha \in \Z^n} c_\alpha z^\alpha$ we will call $A =
\{\alpha \in \Z^n | c_\alpha \not= 0\}$ the support of $W$ and $Newt
(W) \subset \R^n$ the Newton polytope of $W$ where $Newt(W)$ is the
convex hull $Conv(A)$ of $A$.

Given any subset $A \subset \Z^n$ we can consider the space $\C^A$
consisting of all Laurent polynomials whose support is contained in
$A$. Recall from \cite{GKZ} that there is a polynomial $E_A : \C^A
\to \C$ called the principal $A$-determinant whose zero locus
consists generically of those polynomials $f \in \C^A$ which have a
solution to the system of equations \begin{displaymath}
\begin{array}{ccc} f(z) & =
& 0 \\ z_1 \frac{\partial}{\partial z_1} f(z) & = & 0 \\ & \vdots & \\
z_n \frac{\partial}{\partial z_n} f(z)&  = & 0 \end{array}
\end{displaymath}

Now, if $0 \in A$ and $f = \sum_{\alpha \in A} c_\alpha z^\alpha \in
\C^A$ then we let $f_s = f - s \in \C^A$ where $s$ is a constant.
This gives a $1$-dimensional subspace $V_f$ parameterized by $s$ of
$\C^A$ and restricting $E_A$ to this subspace gives a one variable
polynomial $E_A(f_s)$. One can see that the set $Z_f = \{s \in \C |
E_A (f_s) = 0\}$ contains the set of critical values of $f$. We will
see that picking sufficiently generic coefficients for $f$ and
assuming some basic properties for $A$, $Z_f$ will equal the set of
critical values of $f$.

A crucial result in \cite{GKZ} is the description of the Newton
polytope of $E_A$ as the secondary polytope $\Sigma (A)$ of $A$. For
what follows we will require some notation and results on $\Sigma
(A)$ which are taken directly from \cite{GKZ}.

\begin{defn} (i) A marked polytope $(Q, A)$ is a subset $A \in \Z^n$ and
a polytope $Q \subset \R^n$ such that all vertices of $Q$ lie in $A$
and $A \subset Q$ \\ (ii) A subdivision of a marked polytope $(Q,
A)$ consists of marked polytopes $\{(Q_i, A_i)\}$ such that the
dimension of $Q_i$ is $n$ for all $i$, the intersection of any $Q_i$
and $Q_j$ gives a
face of each, $A_i \cap Q_i \cap Q_j = A_j \cap Q_i \cap Q_j$ and the union of all $Q_i$ is $Q$. \\
(iii) A subdivision $\{(Q_i^\prime, A_i^\prime )\}$ of $(Q, A)$ is a
refinement of a subdivision $\{(Q_i, A_i)\}$ of $(Q, A)$ if the set
of all $Q_j^\prime \subset Q_i$ forms a subdivision of $(Q_i, A_i)$. \\
(iv) A subdivision $\{(Q_i, A_i)\}$ of $(Q, A)$ will be called a
triangulation if $A_i$ is an affinely independent subset of $\Z^n$
for all $n$.
\end{defn}

Observe that the definition of a triangulation implies that $Q_i$
are simplices for all $i$. In order to describe $\Sigma (A)$, one
needs the notion of coherent subdivisions. These subdivisions come
from examining the upper boundary of a piecewise linear concave
function. More precisely, consider an element $\psi \in \R^A$ as a
function from $A$ to $\R$. We can form the convex hull $$G_\psi =
Conv\{ (\alpha, y) : y \leq \psi (\alpha ), \alpha \in A , y \in \R
\} \subset \R^n \times \R$$ The upper boundary of $G_\psi$ will be
the graph of a piecewise linear function $g_\psi : Q \to \R$ where
$Q = Conv (A)$. Form a subdivision $S(\psi ) = \{(Q_i, A_i)\}$ of
$(Q, A)$ by letting $Q_i$ be the maximal domains of linearity for
$g_\psi$ and $A_i = \{\alpha \in A : \psi (\alpha ) = g_\psi (\alpha
)\}$. It is immediate that $S(\psi )$ actually forms a subdivision
and moreover, for generic $\psi$, $S(\psi )$ will give a
triangulation.

\begin{defn} A subdivision $\{(Q_i, A_i)\}$ of $(Q, A)$ is called
coherent if there exists a $\psi $ such that $S( \psi ) = \{(Q_i,
A_i)\}$.
\end{defn}

Now, given a triangulation $T = \{(Q_i, A_i)\}$ of $(Q, A)$ and any
$\psi \in \R^A$ we can construct a piecewise linear function
$g_{\psi, T} : Q \to \R$ which is defined by linearly interpolating
$\psi$ inside every individual $Q_i$. \begin{defn} Let $C(T) \subset
\R^A$ consist of those $\psi \in \R^A$ which satisfy \\ (i) The
function $g_{\psi, T}: Q \to \R$ is concave. \\ (ii) For any $\alpha
\in A$ which is not a vertex of any simplex from $T$, we have
$g_{\psi, T} (\alpha ) \geq \psi (\alpha )$.
\end{defn}

Taking all coherent triangulations, the cones $C(T)$ piece together
to form a fan called the secondary fan. It is known that this fan is
the normal fan of a polytope called the secondary polytope $\Sigma
(A)$ of $A$. Given a coherent triangulation $T$, one can define the
element $\varphi_T$ in $(\R^A)^*$ via $$\varphi_T (v_\alpha ) =
\Sigma_{ \alpha \in A_i} Vol (Q_i)$$ where the sum is over all
$(Q_i, A_i)$ appearing in $T$ with $\alpha \in A_i$ and $\{v_\alpha
: \alpha \in A \}$ is the canonical basis of $\R^A$. Define $\Sigma
(A) \subset (\R^A)^*$ to be the convex hull of $\{\varphi_T : T
\hbox{ a coherent triangulation of } (Q, A)\}$ and if $S$ is a
subdivision of $(Q, A)$ let $F(S)$ be the convex hull of all
$\varphi_T$ for which $T$ is a refinement of $S$. Dually, let $C(S)$
be the set of $\psi \in \R^A$ for which $S$ is a refinement of
$S(\psi )$.

\begin{thm}(\cite{GKZ}) The secondary fan is the normal fan of $\Sigma (A)$.
Furthermore, the faces of $\Sigma (A)$ are the polytopes $F(S)$ for
all coherent subdivisions of $(Q, A)$. The normal cone of $F(S)$ is
$C(S)$. \end{thm}

This theorem along with the fact that $Newt(E_A) = \Sigma(A)$ and
the product formula for $E_A$ are all we will need on the secondary
polytope in this paper.

\begin{defn} Let $(Q, A)$ be a marked polytope in $\R^n$ such that $0 \in A$ and $0 \in Int (Q)$. We
will say that a subdivision $\{(Q_0, A_0), (Q_1, A_1)\}$ is a
bisection of $(Q, A)$ if $0 \in Int (Q_0)$ and $A_0 \cup A_1 = A$.
\end{defn}

A bisection corresponds to breaking a given polytope into two
pieces, with one piece containing the origin in its interior. Given
such a bisection $S$, there is an element $\eta_S \in C(S)$ which is
uniquely defined by the properties that $\eta_S \in \Z^A$,
$\eta_S|_{A_0} = 0$ and $\int_Q \eta_S$ is maximal. Now, given a
potential $W = \Sigma_{\alpha \in A } c_\alpha z^\alpha $ and a
bisection $S$ of $(Q, A)$, we define a perturbation of $W$ via
$$W_t( z ) = \sum_{\alpha \in A } c_\alpha t^{- \eta_S
(\alpha)} z^\alpha $$ With such a perturbation of $W$, consider $E_A
(W_t( z )) $. From \cite{GKZ} we have that $E_A$ has Newton polytope
$\Sigma (A)$, so
\begin{displaymath} \begin{array}{ccc} E_A (W_t( z )) &  = & \sum_{\varphi \in \Sigma (A)} d_\varphi
\prod_{\alpha \in A} c_\alpha^{\varphi (\alpha )} t^{- \varphi
(\alpha )\eta_S (\alpha)} \\ & = & \sum_{\varphi \in \Sigma (A)}
\left( d_\varphi \prod_{\alpha \in A} c_\alpha^{\varphi (\alpha )}
\right) t^{- \sum_{\alpha \in A} \varphi (\alpha )\eta_S (\alpha)}
\end{array}
\end{displaymath}
As $\eta_S \in C(S)$ Theorem $5$ implies that $- \eta_S $ achieves
its minimum on $F(S)$; say $- \sum_{\alpha \in A} \varphi (\alpha
)\eta_S (\alpha) = k$ for $\varphi \in F(S)$. One can check that $k$
is positive. Then the above equation can be written
\begin{displaymath} \begin{array}{ccc} E_A (W_t( z ))  & = & t^{k} \sum_{\varphi \in F(S) } \left( d_\varphi \prod_{\alpha \in A}
c_\omega^{\varphi (\alpha )} \right)   + t^{k + 1}{q}(t, c_\alpha)
\end{array}
\end{displaymath}
Where $q(t, c_\alpha)$ is the part of the principal $A$-determinant
with support off of $F(S)$. Now let us recall the product formula
from \cite{GKZ} for principal A-determinants. In general if $S =
\{(Q_i, A_i)\}$ is a subdivision of a marked polytope $(Q, A)$ and
$f$ is a Laurent polynomial with support on $A$, this theorem
asserts that the part of $E_A (f)$ whose support is the face $F(S)$
is simply the product $\prod E_{A_i}(f|_{A_i})$. Here $f|_{A_i}$ is
the part of $f$ whose support lies in $A_i$. Using this, we can
further reduce the above equation to:
\begin{displaymath} \begin{array}{ccc} E_A (W_t( z ))  & = & t^{k} E_{A_0}(W|_{A_0})E_{A_1}(W|_{A_1})   + t^{k + 1}{q}(t, c_\alpha)
\end{array}
\end{displaymath}
Now let us set the constant $c_0$ equal the variable $s$, then we
can regard $E_A(W_t(z))$ as a polynomial $g(s, t)$ in two variables.
For a given value of $t$, assuming $W_t(z)$ is generic, the zeros of
$g(s, t)$ will then be the critical values of $W_t(z)$. One can see
this by investigating the secondary polytope which reveals that, for
all $t \not= 0$, $g(s, t)$ has degree $Vol (Q)$ as a polynomial in
$s$. To see this just choose a generic function $\psi$ on $A$
setting $\psi (0) = 0$ and all other elements negative. This will
yield a triangulation $T$ with $\varphi_T(0) = Vol (Q)$ which is the
degree of $s$ in $E_A(W_t(z))$ (again, assuming generic
coordinates). In particular, the number of zeros up to multiplicity
in the variable $s$ of $E_A(W_t(z))$ is $Vol(Q)$ for $t \not= 0$
which is in line with the number of critical points of $W_t(z)$ as
given by Kouchnirenko's Theorem \cite{Kouchnirenko}. Observe that
for $t = 0$ we have that $g(s, t) = 0$, so $W_0$ is not generic.
However, dividing $g(s, t)$ by $t^k$ we have
\begin{displaymath} \begin{array}{ccc} \frac{g(s, t)}{t^k} & = & E_{A_0}(W|_{A_0})E_{A_1}(W|_{A_1})   + t{q}(t, c_\alpha)
\end{array}
\end{displaymath}
Now, $E_{A_1}(W|_{A_1})$ does not depend on $s$ as $0 \notin A_1$,
so for a generic choice of coefficients this is a non-zero constant.
On the other hand $E_{A_0}(W|_{A_0})$ will be a polynomial in $s$
with zeros at the critical values of $W|_{A_0}$. So, in particular,
as $t$ tends to zero, $g(s, t)/t^k$ converges to a non-zero
polynomial in $s$ with zeros at the critical values of $W|_{A_0}$.
Applying a convergence argument, one sees that all other zeros of
$g(s, t)/t^k$ must tend to infinity as $t \to 0$. For a given $t
\not= 0$ let $\{\nu_1 (t), \ldots, \nu_r (t)\}$ be the zeros of
$g(s, t)$. The preceding argument shows that the critical values of
$W_t$ split into two groups as $t$ tends to zero, the first group
$\{\nu_1(t), \ldots, \nu_m (t)\}$ corresponding to those critical
points of $W|_{A_0}$ and the second $\{\nu_{m + 1}(t), \ldots, \nu_r
(t)\}$ tending out to infinity. One sees then that the number of
critical values in the second group $r - m $ is the number of
critical values of $W|_{A_1}$. Indeed, we would like to identify
these critical values to those of $W|_{A_1}$; this requires
reparameterizing the deformation.

Given an affine map $L: \Z^n \to \Z$, one can reparameterize
$W_t(z)$ to give the map $$W_{L, t}(z_1, \ldots, z_n) = t^{L(0)}
W_t(t^{L(e_1) - L(0)}z_1, \ldots, t^{L(e_n)- L(0)}z_n)$$ Using the
homogeneity properties of $E_A$, we have $$E_A(W_{L, t}) = t^{L_Q }
E_A(W_t)$$ Here $L_Q = L (\int_Q x dx ) + L(0) \cdot [(n + 1)Vol(Q)
- 1]$. Now, one can define another element $\tau_S \in C(S)$ by
setting $\tau_S = \eta_S - L$ where $L$ is the affine map from
$\Z^n$ to $\Z$ which is the restriction of $\eta_S$ to $A_1$. One
notes that $\tau_S$ is zero when restricted to $A_1$ and $\tau_S (0)
< 0$. The graphs of $-\eta_S$ and $-\tau_S$ for $A = \{-1, 0, 1,
2\}$, $A_0 = \{-1, 0, 1\}$ and $A_1 = \{1, 2\}$ are shown below.
\\[5pt]
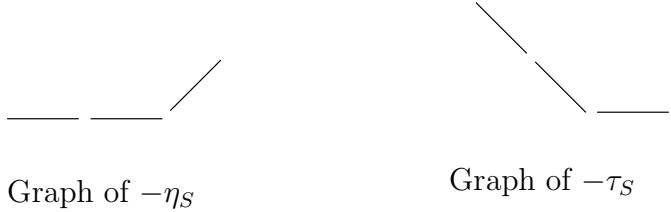
\begin{center}
\begin{picture}(0,0)%
\includegraphics{deform1.pstex}%
\end{picture}%
\setlength{\unitlength}{3947sp}%
\begingroup\makeatletter\ifx\SetFigFont\undefined%
\gdef\SetFigFont#1#2#3#4#5{%
  \reset@font\fontsize{#1}{#2pt}%
  \fontfamily{#3}\fontseries{#4}\fontshape{#5}%
  \selectfont}%
\fi\endgroup%
\begin{picture}(4180,1325)(2239,-1244)
\put(5026,-1111){\makebox(0,0)[lb]{\smash{{\SetFigFont{12}{14.4}{\rmdefault}{\mddefault}{\updefault}{Graph of $-\tau_S$}%
}}}}

\put(2251,-1186){\makebox(0,0)[lb]{\smash{{\SetFigFont{12}{14.4}{\rmdefault}{\mddefault}{\updefault}{Graph of $-\eta_S$}%
}}}}

\end{picture}%
\end{center}

If we define a deformation of $W$ using $\tau_S$ via
$${\tilde W}_t (z) = \sum_{\alpha \in A} c_\alpha t^{-\tau_S (\alpha
)} z^\alpha$$ then one can easily show that ${\tilde W}_t (z) =
W_{L, t}(z)$. Therefore $E_A({\tilde W}_t (z)) = t^{L_Q} g(s, t)$.
However, fixing a $t \not= 0$, the $s$ zeros of $g(s, t)$ will not
represent critical values of ${\tilde W}_t (z)$ as the constant term
is no longer $c_0 = s$. Having reparameterized, the new constant
term is $t^{L(0)}c_0 = t^{L(0)}s $, so the critical values of
${\tilde W}_t (z)$ for a fixed $t \not= 0$ are the points
$\{t^{L(0)}\nu_1 (t), \ldots, t^{L(0)}\nu_r (t)\}$. At $t = 0$ we
see that ${\tilde W}_t (z) = W|_{A_1}$ which has $r - m$ critical
values. As $t$ varies slightly, there will be $r - m$ critical
values an $\varepsilon$-distance away from these values. Assuming
that we chose coefficients of $W$ in a sufficiently generic way, we
may assume that all of the critical values of $W|_{A_1}$ are
distinct. Let us assume also that all of these critical values are
nonzero. Putting all of this together, we see that as $t$ tends
towards zero, the first $m$ critical values $\{t^{L(0)}\nu_1 (t),
\ldots, t^{L(0)}\nu_m (t)\}$ of ${\tilde W}_t(z)$ tend towards zero
since it was shown above that the values $\nu_i (t)$ for $1 \leq i
\leq m$ converge to finite values. On the other hand, there must be
$r - m$ values converging to the critical points of $W|_{A_1}$ and
as these values are distinct, this forces the remaining critical
values $\{t^{L(0)}\nu_{m + 1} (t), \ldots, t^{L(0)}\nu_r (t)\}$ to
converge to these finite values. We codify these results in the
following proposition.
\begin{prop} Given a generic potential $W$ with support $A$ and
Newton polytope $Q$ and an elementary subdivision $\{(Q_0, A_0),
(Q_1, A_1)\}$, one can define two deformations $W_t$ and ${\tilde
W}_t$ of $W$. There is a labeling $\{\nu_1( t), \ldots, \nu_r (t)\}$
of the critical values of $W_t$ for $t \not= 0$ such that
$\{t^{L(0)}\nu_1( t), \ldots, t^{L(0)}\nu_r (t)\}$ are the critical
values of ${\tilde W}_t$ where $L(0) > 0$. Furthermore, as $t$ tends
towards zero, there is an $m$ such that \\ (i) $\{\nu_1 (t) ,
\ldots, \nu_m (t)\}$ approaches the critical values of $W|_{A_0}$.
\\ (ii) $\{t^{L(0)}\nu_{m + 1}( t), \ldots, t^{L(0)}\nu_r (t)\}$
approaches the critical values of $W|_{A_1}$, each of which is
non-zero.
\end{prop}

We now turn to the main result of this subsection. Suppose $(W, U,
a)$ is a partial Lefschetz fibration. If $W_t$ is a perturbation of
the potential, one sees that for a sufficiently small $\varepsilon$
the openness of transversality ensures that for $t \in
D_\varepsilon$, $(W_t, U, a)$ is a partial Lefschetz fibration if
and only if no critical values of $W_t$ enter into $U$. This fact
combined with the previous proposition gives us the following
theorem.

\begin{thm} Suppose $W$ is a Laurent polynomial with support $A$ and
Newton polytope $Q$ and  $\{(Q_0, A_0), (Q_1, A_1)\}$ is a bisection
of $(Q, A)$. Assume $W$ is generic in the following sense: \\ (i)
$W$ is a Lefschetz fibration on $\gmc$
\\(ii) $W|_{A_0}$ is a Lefschetz fibration on $\gmc$ \\ (iii)
$(W|_{A_1}, \C - D_{\varepsilon}, a)$ is a partial Lefschetz
fibration for sufficiently large $a$ \\ (iv) $W_{A_1}$ has nonzero
critical values \\ Then, after choosing a distinguished basis of
paths $\{\delta_i \}$ for $(W|_{A_1}, \C - D_{\varepsilon}, a)$,
there is the following semi-orthogonal decomposition:
$$\D (Fuk (W)) = \left< \D (Fuk (W|_{A_0}), \D (Fuk
(W|_{A_1}, \C - D_\varepsilon , a, \{\delta_i\}) \right>$$
\end{thm}
\begin{proof} By Proposition 10 and the above observation we have that
${\tilde W}_t$ is a deformation of the partial Lefschetz fibration
$(W|_{A_1}, \C - D_{\varepsilon}, a)$ for $|t| < \epsilon$. But, by
assumption (i), $W = {\tilde W}_1$ is a Lefschetz fibration and from
Proposition 10 ${\tilde W}_t$ is a deformation of $W$ for $t \not=
0$. Thus
$$\D (Fuk (W )) \simeq \D (Fuk ({\tilde W_\epsilon})) = \left< \T ,
\D (Fuk (W|_{A_1}, \C - D_\varepsilon , a, \{\delta_i\}) \right>$$
where the last equality follows from Theorem $4$ and the first from
Corollary $3$. Now, the category $\T$ is generated by those
vanishing cycles coming from critical values of ${\tilde
W}_\epsilon$ in $D_\varepsilon$. Deforming to $W$ and then using the
deformation $W_t$ to deform to $W_0$, Corollary $3$ gives us that
$\T \simeq \D (Fuk (W_0, D_R , a^\prime, \{\gamma_i\}))$ for some
$R$. But as $W_0 = W|_{A_0}$ the latter category here is just $\D (
Fuk (W|_{A_0}))$.
\end{proof}

\section{Homological Mirror Symmetry for Weighted Projective
Blowups of Toric Surfaces}
\subsection{The potential $W_\A$ on $\C^2$} Recall that homological mirror symmetry for Fano toric
varieties asserts an equivalence between $\D (Coh (X_\tr))$ and $\D
(Fuk (W_\tr))$ where $W_\tr$ is a Laurent polynomial with Newton
polytope equal to the convex hull of $\tr(1)$, the primitives of the
one dimensional cones of the fan $\tr$. Suppose $\blo$ is a weighted
projective blowup of a smooth toric variety of a fan $\tr$ at a
$\gmc$ invariant point. In section $2$ we saw that there is an $n +
1$ cone $\sigma $ in $\tr $ generated by $\{v_0, \ldots, v_n\}$ such
that $\blo$ is the toric stack of a fan $\tilde \tr$ with ${\tilde
\tr}(1) = \tr (1) \cup \{a_0v_0 + \cdots + a_nv_n\}$. Thus, if
$W_\tr$ is the potential associated to $\tr$ then $W_{\tilde \tr} =
W_\tr + C z_0^{a_0} \cdots z_n^{a_n}$ (in fact, one could use other
monomials in this expansion, but we will see that this is not
necessary). The picture below shows the situation for the $\A = (2,
3)$ weighted projective blowup of $\p^2$.
\\[8pt]
\begin{picture}(0,0)%
\includegraphics{blowup2.pstex}%
\end{picture}%
\setlength{\unitlength}{3947sp}%
\begingroup\makeatletter\ifx\SetFigFont\undefined%
\gdef\SetFigFont#1#2#3#4#5{%
  \reset@font\fontsize{#1}{#2pt}%
  \fontfamily{#3}\fontseries{#4}\fontshape{#5}%
  \selectfont}%
\fi\endgroup%
\begin{picture}(5424,3214)(2389,-3419)
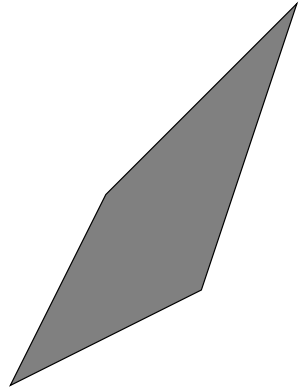
\put(6001,-3361){\makebox(0,0)[lb]{\smash{{\SetFigFont{12}{14.4}{\rmdefault}{\mddefault}{\updefault}{Newton Polytope for $W_{\tilde \tr}$}%
}}}}
\put(4201,-361){\makebox(0,0)[lb]{\smash{{\SetFigFont{12}{14.4}{\rmdefault}{\mddefault}{\updefault}{  $u = (2, 3)$}%
}}}}
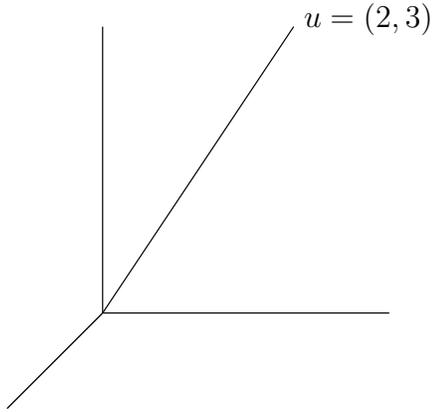
\put(2401,-3361){\makebox(0,0)[lb]{\smash{{\SetFigFont{12}{14.4}{\rmdefault}{\mddefault}{\updefault}{Fan for $Bl_{(2, 3)}(\p^2)$}%
}}}}
\end{picture}%
\\[6pt]
Now, in section $2$ we saw that given a weighted projective blowup
$\blo$ of $X$, the derived category can be expressed as a
semi-orthogonal decomposition $\left< \D (Coh (X)) , \T \right>$
where $\T$ has an exceptional collection Koszul dual to the
pushforwards of line bundles. On the other hand, in section $3$ we
showed that the derived Fukaya category of a potential $W$ on $\gmc$
also has a semi-orthogonal decomposition associated to a bisection.
For the potential $W_{\tilde \tr}$, one sees that in many situations
there is a clear bisection of the convex hull of the primitives of
${\tilde \tr}(1)$. Namely, take the convex hull of $\tr$ and the
convex hull $Q_\A$ of $A_\A = \{v_0, \ldots, v_n, a_0v_0 + \cdots
a_n v_n\}$. Note that in order for this to be a bisection one must
have that $Conv(\tr(1)) \cup Q_A = Conv({\tilde \tr}(1))$. This will
be assumed for the rest of this paper. We illustrate this
subdivision below for the potential associated to $Bl_{(2,
3)}(\p^2)$.
\\[8pt]
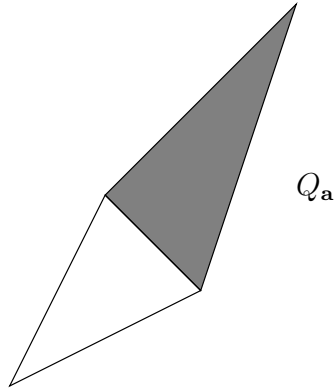
\begin{center}
\begin{picture}(0,0)%
\includegraphics{blowup3.pstex}%
\end{picture}%
\setlength{\unitlength}{3947sp}%
\begingroup\makeatletter\ifx\SetFigFont\undefined%
\gdef\SetFigFont#1#2#3#4#5{%
  \reset@font\fontsize{#1}{#2pt}%
  \fontfamily{#3}\fontseries{#4}\fontshape{#5}%
  \selectfont}%
\fi\endgroup%
\begin{picture}(1824,2424)(4189,-2773)
\put(6001,-1561){\makebox(0,0)[lb]{\smash{{\SetFigFont{12}{14.4}{\rmdefault}{\mddefault}{\updefault}{$Q_\A$}%
}}}}
\end{picture}%
\end{center}

Now, there are several elementary examples where it is known that
Homological mirror symmetry holds, e.g. weighted projective planes
and smooth Fano toric surfaces (\cite{AKO}, \cite{Ueda}). Thus our
strategy is to assume that homological mirror symmetry holds for
$\tr$ and show that the extra semi-orthogonal categories are
equivalent as well. Although this does not fully prove mirror
symmetry, as the interaction between the two semi-orthogonal
categories is neglected, it does indicate that the mirror categories
have equivalent pieces in natural decompositions. We will pursue
this strategy for weighted projective blowups of smooth toric
surfaces.

In general, the potential $W_\A$ associated to $(A_\A , Q_\A)$ can
be written $$W_\A (z_0, \ldots, z_n) = a_0 z_0 + \cdots + a_n z_n -
z_0^{a_0} \cdots z_n^{a_n}$$ It is the case that this potential is
generic enough to fit into the definition of a partial Lefschetz
fibration. It also has the advantage of having a good deal of
symmetry which helps in its investigation. Observe that $$\partial_i
W_\A(z_0, \ldots, z_n) = a_i - a_i z_0^{a_0} \cdots z_i^{a_i - 1}
\cdots z_n^{a_n}$$ so if $(c_0, \ldots, c_n)$ is a critical point of
$W_\A$ then $c_0^{a_0} \cdots c_i^{a_i - 1} \cdots c_n^{a_n} = 1 =
c_0^{a_0} \cdots c_j^{a_j - 1} \cdots c_n^{a_n}$ for every $i$ and
$j$. In particular, $c_i \not= 0$ for all $i$ and dividing the above
equation by $c_0^{a_0} \cdots c_i^{a_i - 1} \cdots c_j^{a_j - 1}
 \cdots c_n^{a_n}$ we have $c_i = c_j$. Thus a critical point of
 $W_\A$ can be written $(c, \ldots , c)$ where $c$ must satisfy the
 equation $c^{l - 1} = 1$ where $l = \sum_{i = 0}^n a_i$. Thus there are $l - 1$ critical points
 $\{p_0, \ldots, p_{l - 2}\}$ such that $p_i = (\zeta^i, \ldots, \zeta^i)$ where
 $\zeta$ is the $(l -1)$-th root of unity $e^{2\pi i /(l - 1)}$. One
 can compute that the corresponding critical value $q_i$ of $p_i$ as
 $q_i = W_\A (p_i) = (l - 1) \cdot \zeta^i$.

At first it appears that $W_\A$ has smooth fibers outside the
critical values $q_i$ and that parallel transport could be defined
for paths. This is indeed the case if we consider $W_\A$ as a
fibration on the partial compactification $\C^{n + 1}$ of $(\C^*)^{n
+ 1}$. However, there are two objections to viewing $W_\A$ as a
Lefschetz fibration on $\C^{n + 1}$. First, the standard K\"ahler
structure on $\C^{n +1}$ is not an extension of that on the complex
torus. This can easily be remedied by perturbing and using Theorem
$4$. The second and more serious objection is that the zero fiber of
$W_\A$ does not transversely intersect the normal crossing divisor
$\bigcup \{z_i = 0\}$, so the topology of this fiber in $(\C^*)^{n +
1}$ will abruptly change from those of neighboring fibers.
Nevertheless, it will be helpful to first regard $W_\A$ as a
potential on $\C^{n + 1}$ and obtain information on the vanishing
cycles at the zero fiber. We then work backwards and perturb to
neighboring fibers where such topological pathologies do not occur.
Pursuing this line we will use the method of matching paths which
rely on an auxiliary holomorphic function $f (z_0, \ldots, z_n) =
a_0z_0 + \cdots + a_n z_n$ (\cite{SeidVC}). To avoid confusion, we
will write the fiber at $q$ of $W_\A$ regarded as a map on $\C^{n +
1}$ as $F_q$, and when considering $W_\A$ as a map on $(\C^*)^{n +
1}$ we will use the notation $W^{-1}_\A(q)$.

For any fiber $F_q$ we can restrict $f$ to $F_q$ and find the
critical points of $f$. This will be the intersection of $F_q$ with
the discriminant variety $D \subset \C^{n + 1}$ which is defined as
the variety where the differential of $(W_\A, f) : \C^{n + 1} \to
\C^2$ has rank $< 2$. The variety $D$ can be seen to be the variety
defined by $(df \wedge dW)_{z} = 0$. One can compute
\begin{displaymath} (df \wedge dW)_{z} = \sum_{0
\leq i < j \leq n} a_i a_j (z_0^{a_0} \cdots z_j^{a_j - 1} \cdots
z_n^{a_n} - z_0^{a_0} \cdots z_i^{a_i - 1} \cdots z_n^{a_n}) dz_i
\wedge dz_j \end{displaymath} so that the discriminant variety is
simply
$$D = \{(c, \ldots, c)\} \cup (\cup_{a_i > 1} \{z_i = 0\})$$  Thus
the discriminant variety can be written as the union of two
components $C_1$ and $C_2$. The diagonal component $C_1$ contains
the critical points of $W_\A$, and the other component $C_2$ is a
union of divisors $\{z_i = 0\}$. Let us examine the critical values
of $f$ restricted to the fiber $F_q$ associated to the second
component of the discriminant variety. Since any element $p \in F_q
\cap C_2$ must have at least one $z_i = 0$ we have that $z_0^{a_0}
\cdots z_n^{a_n} = 0$ so that $q = W_\A(p) = f(p)$. Thus $f(F_q \cap
C_2) = q$; more generally, we have
$$f|_{F_q}^{-1}(q) = F_q \cap (\cup_{i = 0}^n \{z_i = 0\})$$  If
$C_2 = \emptyset$ (i.e. if $a_i = 1$ for all $i$) it will still be
convenient to keep track of this image. Doing so will allow us to
see which fibers have vanishing cycles intersecting the divisor
$\cup_{i = 0}^n \{z_i = 0\}$.

We now examine the the critical values of $f$ restricted to a fiber
$F_q$ associated to the diagonal component $C_1$. These will be
$f(c, \ldots, c) = l \cdot c$ where $c$ satisfies $W_\A (c, \ldots
c) = l c - c^l = q$. In other words, these critical values are the
roots of the polynomial $h_q(x) = x^l - l^l x + l^l q$. Either by
the theory of matching paths or by direct computation, one sees that
$h_q(x)$ has multiple roots only if $q$ is among the critical values
of $W_\A$. Furthermore, $q$ is a root of $h_q(x)$ only for $q = 0$
which implies $f(C_1 \cap F_q) \cap f(C_2 \cap F_q)$ is empty for
all $q \not= 0$. This implies that for a path $\gamma: [0, 1] \to
\C^*$, the vanishing cycle in $F_{\gamma(t)}$ can be consistently
isotoped so that it does not intersect the divisor $\cup_{i = 0}^n
\{z_i = 0\}$ and thus is contained in $(\C^*)^{n + 1}$.

 Now, observe that $\{0, l^{l / l -1}, l^{l / l - 1} \zeta ,
\ldots, l^{l /l - 1} \zeta^{l - 2}\}$ are the roots of $h_0(x)$ and
thus are the critical values of $f$ restricted to $F_0$. We would
like to see the movement of the critical values of $f|_{F_q}$ as $q$
moves from $0$ to $l - 1$. These will be the zeros of $h_q(x)$ and
the point $q$ itself. Let us focus on the real roots of $h_q(x)$ for
a moment. We observe that, as a real function, $h_0(x)$ has two
different graphs depending on whether $l$ is even or odd. These are
drawn below:
\\[8pt]
\begin{picture}(0,0)%
\includegraphics{graphs1.pstex}%
\end{picture}%
\setlength{\unitlength}{3947sp}%
\begingroup\makeatletter\ifx\SetFigFont\undefined%
\gdef\SetFigFont#1#2#3#4#5{%
  \reset@font\fontsize{#1}{#2pt}%
  \fontfamily{#3}\fontseries{#4}\fontshape{#5}%
  \selectfont}%
\fi\endgroup%
\begin{picture}(5754,2376)(1849,-2819)
\put(5476,-2761){\makebox(0,0)[lb]{\smash{{\SetFigFont{12}{14.4}{\rmdefault}{\mddefault}{\updefault}{$h_0(x)$ for $l$ odd}%
}}}}
\put(2401,-2761){\makebox(0,0)[lb]{\smash{{\SetFigFont{12}{14.4}{\rmdefault}{\mddefault}{\updefault}{$h_0(x)$ for $l$ even}%
}}}}
\end{picture}%

As $q$ tends from $0$ to $l - 1$, we see that the real roots of
$h_q(x)$ are simply the $x$-coordinates of the intersection of the
horizontal line at $y = - l^l q$ with the above graphs. Thus, the
two real roots $0$ and
 $l^{l/l - 1}$ contract together along the real axis without any
 other critical values associated to $C_1$ passing between them. Furthermore, the critical value
 associated to $C_2$ is simply the intersection of $y = - l^l q$ with the line $y = - l^l x$. Now, $h_0(x) - (-l^l x) = x^l$ is positive for
 all $x > 0$; so the real critical value associated to $C_2$ does not pass between the contracting roots either (although the $0$ root "hits" this value
 for $q = 0$). To put this in the language of
 matching paths, we can say that to the path ${\tilde \delta}_0 (t) = (l - 1) t$
 from the regular value $0$ to the critical value $l - 1$, we
 associate the matching path $\eta_0 (s) = l^{l/ l - 1} s$. These
 two paths occur on different lines, $\tilde \delta_0$ occurs in the
 range of $W_\A$ while $\eta_0$ occurs on the image of $F_0$ via $f$
 connecting the critical values. This is illustrated below:
 \\[8pt]
 \begin{center}
 \begin{picture}(0,0)%
\includegraphics{mpath1.pstex}%
\end{picture}%
\setlength{\unitlength}{3947sp}%
\begingroup\makeatletter\ifx\SetFigFont\undefined%
\gdef\SetFigFont#1#2#3#4#5{%
  \reset@font\fontsize{#1}{#2pt}%
  \fontfamily{#3}\fontseries{#4}\fontshape{#5}%
  \selectfont}%
\fi\endgroup%
\begin{picture}(5063,2530)(1414,-2744)
\put(4201,-2686){\makebox(0,0)[lb]{\smash{{\SetFigFont{12}{14.4}{\rmdefault}{\mddefault}{\updefault}{and matching path $\eta_0$}%
}}}}
\put(1426,-2461){\makebox(0,0)[lb]{\smash{{\SetFigFont{12}{14.4}{\rmdefault}{\mddefault}{\updefault}{Critical Values of $W_\A$}%
}}}}
\put(1426,-2686){\makebox(0,0)[lb]{\smash{{\SetFigFont{12}{14.4}{\rmdefault}{\mddefault}{\updefault}{and path $\tilde \delta_0$}%
}}}}
\put(4201,-2461){\makebox(0,0)[lb]{\smash{{\SetFigFont{12}{14.4}{\rmdefault}{\mddefault}{\updefault}{Critical Values of $f|_{F_0}$}%
}}}}
\end{picture}%
\end{center}

Observe that $W_\A(\zeta {\bf z}) = \zeta W_\A ({\bf z})$ and
$f(\zeta {\bf z}) = \zeta f ({\bf z})$. This symmetry implies that
the matching path associated to ${\tilde \delta_i}(t) = \zeta^i
{\tilde \delta_0}(t)$ is $\eta_i (s) = \zeta^i \eta_0 (s)$. Thus, in
the picture above, one can rotate both the path in the image of
$W_\A$ and its matching path in the image of $f|_{F_0}$ by the angle
$2 \pi \over (l  - 1)$ to obtain all other pairs $({\tilde
\delta_i}, \eta_i)$. Now, by the theory of matching paths, a
matching path $\eta_i$ for $\tilde \delta_i$ is isotopic to the
image of the vanishing cycle $L_i$ under the map $f$, and the
endpoints of the matching path are the critical values of a standard
height function on the vanishing cycle $L_i \simeq S^n$. This
machinery applies if the critical points of $f$ on $F_q$
corresponding to the endpoints are Morse. While this is the case for
all $q \not=0$ , once $q = 0$ we acquire a more complicated
singularity on $F_0$. We will ignore this fact as we are going to
perturb to a neighboring fiber in the next section, however, the
vanishing cycles we will describe in $F_0$ will necessarily have
corners at their common intersection which is in $F_0 \cap
f^{-1}(0)$. Thus we will identify the vanishing cycle $L_i$
associated to ${\tilde \delta}_i$ with the component of the
pre-image $f^{-1}(\eta_i)$ which contains the critical point $p_i$.
After doing so, we see that all of the intersection points of
vanishing cycles lie in the subspace $F_0 \cap f^{-1}(0)$.

At this point we restrict our attention to the case of $n = 1$. It
will be convenient to adopt the convention that $a_0 \leq a_1$ which
we will assume throughout. Observe that $F_0 \cap f^{-1}(0) = \{(0,
0)\}$ so that all vanishing cycles $L_i$ in $F_0$ intersect in
precisely one point, namely $\{(0, 0)\}$. In what follows, we will
represent $F_0$ as a $(l - 1)$-fold branched covering of $\p^1 -\{0,
\infty\}$ and describe the vanishing cycles as lifts of a single
curve in $\p^1 - \{0, \infty\}$. To pursue this aim, we quotient
$\C^2$ and $\C$ by the diagonal action of $G = \{\zeta^j\ : 0 \leq j
\leq l - 2 \}$. We have the diagram:
\begin{displaymath}
\begin{array}{ccc}
\C^2 & {\buildrel W_\A \over \rightarrow} & \C \\
{\psi \scriptstyle \downarrow} & & {\downarrow \scriptstyle z^{l - 1}} \\
\C^2 / G & {\buildrel {\tilde W_\A} \over \hookrightarrow} & \C
\end{array}
\end{displaymath}
Here we take $$\psi (z_0, z_1) = (z_0^{a_0 - 1} z_1^{a_1}, z_0^{a_0}
z_1^{a_1 - 1}) = (u_0, u_1)$$ for $(z_0, z_1) \in (\C^*)^2$ and find
that
$${\tilde W_\A}(u_0, u_1) = u_0^{a_0} u_1^{a_1}
\left(\frac{a_0}{u_0} + \frac{a_1}{u_1} - 1 \right)^{l - 1}$$ for
$u_0 \not= 0 \not= u_1$. Thus the closure of ${\tilde W_\A}^{-1}(0)$
is the variety $V = \{(u_0, u_1): a_0u_1 + a_1 u_0 - u_1 u_0 = 0\}$.
One observes that $V \cap \{u_0 = 0\} = \{(0, 0)\} = V \cap \{u_1 =
0\}$ which implies that $\psi: F_0 \to V$ is a branched covering
with branch point $\{(0, 0)\}$ whose ramification is $l - 1$.
Furthermore, the defining equation for $V$ is a quadratic which
implies $V$ is $\p^1$ minus the two points intersecting the divisor
at infinity. In other words, one can view the compactification of
$V$ in $\p^2$ as the hypersurface $Z = Zero (a_0u_1 u_2 + a_1 u_0
u_2 - u_1 u_0)$. Then we have $V = Z - \{[0:1:0], [1:0:0]\}$.

Now, the fiber of $\psi$ over any point $p \in V -\{(0, 0)\}$ is a
torsor over $\Z / (l - 1) \Z \simeq G$ where the action is given by
multiplication by $\zeta$. We will describe the monodromy of a loop
around the points $[0:1:0]$ and $[1:0:0]$ to give a complete
description of the topology of $F_0$. For this we give a coordinate
map $j: \C \to \C^2 \subset \p^2$ for $Z$ where $j(w) = (w, a_1 w /
w - a_0)$. Observe that in $\p^2$, $j(a_0) = [0: 1: 0]$. Now, in
$(\C^*)^2$ there is a global multiple valued inverse of $\psi$ which
can be written $$\psi^{-1}(u_0, u_1) = \left( \sqrt[l -
1]{\frac{u_1^{a_1}}{u_0^{a_1 - 1}}} \hbox{ },\sqrt[l -
1]{\frac{u_0^{a_0}}{u_1^{a_0 - 1}}} \hbox{ } \right) $$ In the
coordinate $w$ this gives $$\psi^{-1} (w) = \left(\sqrt[l -
1]{\frac{a_1^{a_1}w}{(w - a_0)^{a_1}}} \hbox{ },\sqrt[l - 1]{\frac{w
(w - a_0)^{a_0 - 1}}{a_1^{a_0 - 1}}} \hbox{ } \right) $$ From this
one can see that the monodromy around $a_0$ is multiplication by
$\zeta^{-a_1}$, or, from the torsor point of view, it is addition by
$-a_1$. To find the monodromy around the point $[1: 0: 0]$ we simply
add the monodromies around $[0:0:1]$ and $[0:1:0]$ yielding $a_0$ as
$1 - a_1 \equiv a_0 (mod (l - 1))$. This gives a complete picture of
the topology of $F_0$.

Now we would like to investigate the vanishing cycles $L_i \subset
F_0$ (with corners). These are maps of $S^1$ in $F_0$, which are
embeddings except at the point $(0,0)$ which may have a corner. Each
image contains the point $(0,0)$ which is the only intersection
point of any two; furthermore we have $L_i = \zeta^i L_0$. Thus
$\psi ( L_i) = \psi (L_j)$ for all $i$ and $j$, $(0, 0) \in \psi
(L_0)$ and $\psi (L_0 )$ has no self intersection points in $V$.
Now, $V$ is topologically a cylinder and the only closed curves with
no intersection points in $V$ up to isotopy are therefore the
homotopically trivial curve or the curve wrapping around the
cylinder. From standard theory on vanishing cycles, the middle
dimensional homology of the fiber is generated by the vanishing
cycles, so we can rule out the trivial curve. We summarize these
results in the following proposition and picture.
\begin{prop} The zero fiber $F_0 \subset \C^2$ admits a $(l -
1)$-fold branched covering $\psi$ to $Z - \{[1:0:0], [0:1:0]\}$ with
one branch point over $[0:0:1]$ of ramification $l - 1$. The
monodromy around the points $[1:0:0]$ and $[0:1:0]$ are $a_0$ and
$-a_1$ respectively. All vanishing cycles $L_i$ have the same image
$\psi (L_i)$ which is isotopic to a closed curve containing the
point $[0:0:1]$, generating the homotopy of $Z - \{[1:0:0],
[0:1:0]\}$ and containing no self intersections. \end{prop}

\begin{center}
\begin{picture}(0,0)%
\includegraphics{base2.pstex}%
\end{picture}%
\setlength{\unitlength}{3947sp}%
\begingroup\makeatletter\ifx\SetFigFont\undefined%
\gdef\SetFigFont#1#2#3#4#5{%
  \reset@font\fontsize{#1}{#2pt}%
  \fontfamily{#3}\fontseries{#4}\fontshape{#5}%
  \selectfont}%
\fi\endgroup%
\begin{picture}(5701,3364)(3301,-3344)
\put(6301,-2761){\makebox(0,0)[lb]{\smash{{\SetFigFont{12}{14.4}{\rmdefault}{\mddefault}{\updefault}{$[0:0:1]$}%
}}}}
\put(3976,-1111){\makebox(0,0)[lb]{\smash{{\SetFigFont{12}{14.4}{\rmdefault}{\mddefault}{\updefault}{$+a_0$}%
}}}}
\put(7801,-1036){\makebox(0,0)[lb]{\smash{{\SetFigFont{12}{14.4}{\rmdefault}{\mddefault}{\updefault}{$-a_1$}%
}}}}
\put(8326,-136){\makebox(0,0)[lb]{\smash{{\SetFigFont{12}{14.4}{\rmdefault}{\mddefault}{\updefault}{$[0:1:0]$}%
}}}}
\put(3301,-136){\makebox(0,0)[lb]{\smash{{\SetFigFont{12}{14.4}{\rmdefault}{\mddefault}{\updefault}{$[1:0:0]$}%
}}}}
\put(5176,-2011){\makebox(0,0)[lb]{\smash{{\SetFigFont{12}{14.4}{\rmdefault}{\mddefault}{\updefault}{$+1$}%
}}}}
\put(6451,-961){\makebox(0,0)[lb]{\smash{{\SetFigFont{12}{14.4}{\rmdefault}{\mddefault}{\updefault}{$\psi (L_i)$}%
}}}}
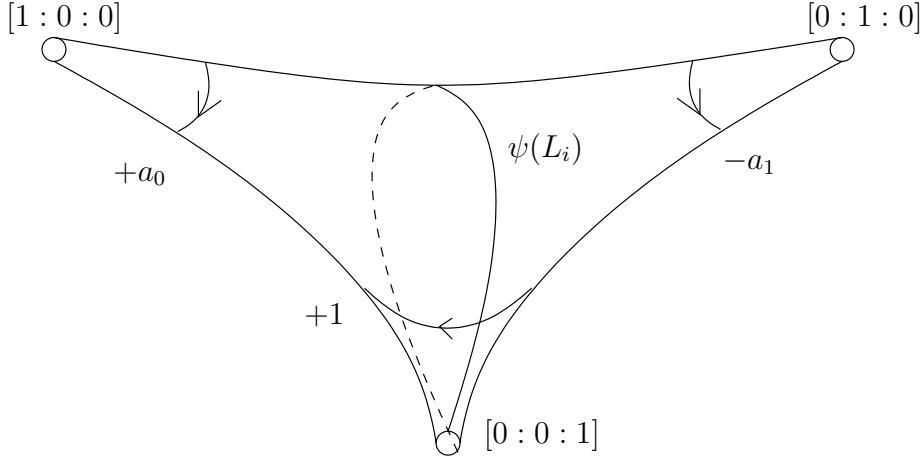
\put(5251,-3286){\makebox(0,0)[lb]{\smash{{\SetFigFont{12}{14.4}{\rmdefault}{\mddefault}{\updefault}{The base $V$ of $\psi$}%
}}}}
\end{picture}%

\end{center}
One can see from the above picture that any holomorphic discs
connecting intersection points of $L_i$ in a nearby fiber must map
to a neighborhood of $[0:0:1]$. Thus we ought to examine the
placement of the vanishing cycles in a neighborhood of $(0, 0) \in
F_0$. If we take an $\delta$-disc neighborhood $D$ of $(0, 0)$ in
$F_0$, then $\psi (z) = z^{l - 1}$ for a local chart of $V$. Now,
taking slits from $[0:0:1]$ to $[0:1:0]$ and from $[0:1:0]$ to
$[1:0:0]$ gives $l -1$ domains gluing together to form
$W_\A^{-1}(0)$. We fix the first slit to give the line segments
$\{r\zeta^k : 0 \leq r \leq \delta \}$ in the neighborhood $D$.
Then, after an isotopy, we have that the vanishing cycles in $D$ are
approximately the lines $\ell_k = \{r e^{(2 k + 1) \pi i / 2(l -
1)}: 0 \leq r \leq \delta \}$ for $0 \leq k \leq 2l - 3$ (in fact
these would give the vanishing cycles with corners at $0$).
Multiplying this picture by a power of $\zeta$ if necessary, we can
assume that $\ell_0 \subset L_0$. Then, following the curve $L_0$
around $V$ we see that it passes through the second slit from
$[0:1:0]$ to $[0:0:1]$ so that it gains $-a_0$ monodromy.
Additionally, the angle that $L_0$ returns with to $(0, 0)$ acquires
a $\zeta^{1/2}$ factor. Thus if $\ell_0$ is the outgoing part of
$L_0$ in $D$ then $\zeta^{{1/2} - a_0}\ell_0 = \ell_{1 - 2a_0}$ is
the incoming part of $L_0$ where the subscript in the second line
segment should be taken modulo $2(l - 1)$. Utilizing the fact that
$\zeta^{i} L_0 = L_i$ we have that
$$L_i \cap D = \ell_{2i} \cup \ell_{1 - 2a_0 + 2i}$$ We illustrate
this for $\A = (1,4)$ and $\A = (2, 3)$ below.

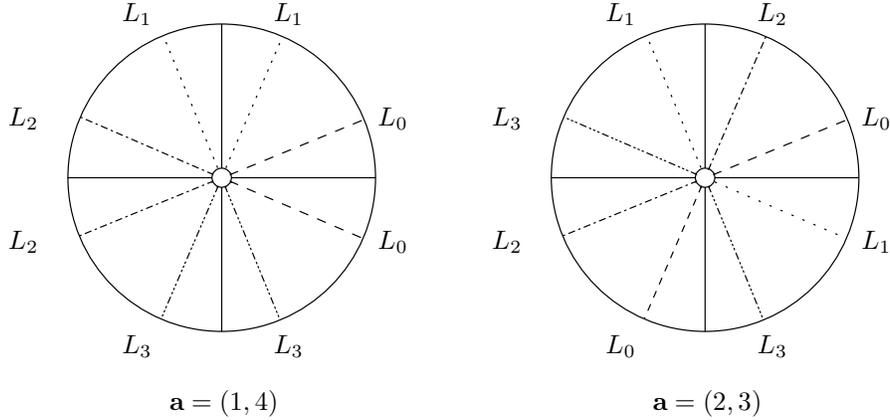
\begin{center}
\begin{picture}(0,0)%
\includegraphics{circ3b.pstex}%
\end{picture}%
\setlength{\unitlength}{3947sp}%
\begingroup\makeatletter\ifx\SetFigFont\undefined%
\gdef\SetFigFont#1#2#3#4#5{%
  \reset@font\fontsize{#1}{#2pt}%
  \fontfamily{#3}\fontseries{#4}\fontshape{#5}%
  \selectfont}%
\fi\endgroup%
\begin{picture}(5736,2609)(3901,-2815)
\put(7949,-2769){\makebox(0,0)[lb]{\smash{{\SetFigFont{10}{12.0}{\rmdefault}{\mddefault}{\updefault}{$\A = (2, 3)$}%
}}}}
\put(6223,-984){\makebox(0,0)[lb]{\smash{{\SetFigFont{10}{12.0}{\rmdefault}{\mddefault}{\updefault}{$L_0$}%
}}}}
\put(6223,-1758){\makebox(0,0)[lb]{\smash{{\SetFigFont{10}{12.0}{\rmdefault}{\mddefault}{\updefault}{$L_0$}%
}}}}
\put(5567,-330){\makebox(0,0)[lb]{\smash{{\SetFigFont{10}{12.0}{\rmdefault}{\mddefault}{\updefault}{$L_1$}%
}}}}
\put(4616,-330){\makebox(0,0)[lb]{\smash{{\SetFigFont{10}{12.0}{\rmdefault}{\mddefault}{\updefault}{$L_1$}%
}}}}
\put(3901,-984){\makebox(0,0)[lb]{\smash{{\SetFigFont{10}{12.0}{\rmdefault}{\mddefault}{\updefault}{$L_2$}%
}}}}
\put(3901,-1758){\makebox(0,0)[lb]{\smash{{\SetFigFont{10}{12.0}{\rmdefault}{\mddefault}{\updefault}{$L_2$}%
}}}}
\put(4616,-2413){\makebox(0,0)[lb]{\smash{{\SetFigFont{10}{12.0}{\rmdefault}{\mddefault}{\updefault}{$L_3$}%
}}}}
\put(5567,-2413){\makebox(0,0)[lb]{\smash{{\SetFigFont{10}{12.0}{\rmdefault}{\mddefault}{\updefault}{$L_3$}%
}}}}
\put(9259,-984){\makebox(0,0)[lb]{\smash{{\SetFigFont{10}{12.0}{\rmdefault}{\mddefault}{\updefault}{$L_0$}%
}}}}
\put(7652,-330){\makebox(0,0)[lb]{\smash{{\SetFigFont{10}{12.0}{\rmdefault}{\mddefault}{\updefault}{$L_1$}%
}}}}
\put(6938,-1758){\makebox(0,0)[lb]{\smash{{\SetFigFont{10}{12.0}{\rmdefault}{\mddefault}{\updefault}{$L_2$}%
}}}}
\put(8604,-2413){\makebox(0,0)[lb]{\smash{{\SetFigFont{10}{12.0}{\rmdefault}{\mddefault}{\updefault}{$L_3$}%
}}}}
\put(7652,-2413){\makebox(0,0)[lb]{\smash{{\SetFigFont{10}{12.0}{\rmdefault}{\mddefault}{\updefault}{$L_0$}%
}}}}
\put(9259,-1758){\makebox(0,0)[lb]{\smash{{\SetFigFont{10}{12.0}{\rmdefault}{\mddefault}{\updefault}{$L_1$}%
}}}}
\put(8604,-330){\makebox(0,0)[lb]{\smash{{\SetFigFont{10}{12.0}{\rmdefault}{\mddefault}{\updefault}{$L_2$}%
}}}}
\put(6938,-984){\makebox(0,0)[lb]{\smash{{\SetFigFont{10}{12.0}{\rmdefault}{\mddefault}{\updefault}{$L_3$}%
}}}}
\put(4913,-2769){\makebox(0,0)[lb]{\smash{{\SetFigFont{10}{12.0}{\rmdefault}{\mddefault}{\updefault}{$\A = (1, 4)$}%
}}}}
\end{picture}%

\end{center}
\subsection{The potential $W_\A$ on $(\C^*)^2$} We will now perturb
the picture presented in the previous subsection to a fiber
$W_\A^{-1}(\varepsilon)$ in $(\C^*)^2$. The first observation we
make is that the point $(0, 0) \in F_0$ no longer exists when $W_\A$
is considered as a potential on $(\C^*)^2$. Furthermore, this end of
the fiber splits into two ends after perturbing slightly. More
precisely, for a given $q \in \C$,  we have that $F_q \cap (\{z_0 =
0\} \cup \{z_1 = 0\}) = \{(0, q /a_1), (q / a_0, 0)\}$. In the
$(\C^*)^2$ picture these two points are not included in the fiber.
Furthermore, if one parallel transports the fiber of $W_\A$ around a
small $\varepsilon$ circle about $0$, the resulting monodromy map
will be a full Dehn twist of the two ends. In the $\C^2$ picture,
this monodromy map is isotopic to the identity; however, after
removing the two points, the map becomes non-trivial. In what
follows, we will give an explicit map which is symplectically
isotopic to the parallel transport map for an arc about the origin
in $\C$.

Fix $\varepsilon$ to be a sufficiently small real number and let
$q_\theta = \varepsilon e^{i\theta}$. Given any fiber
$W_\A^{-1}(q_\theta)$, we can form a local chart $U_\theta = \{w :
|w| < 4 \varepsilon, w \not= 0, w \not= q_\theta\}$  centered at one
of the ends, say $(0, \varepsilon e^{i \theta} /a_1)$ such that the
other end is $q_\theta$. Indeed, from the previous subsection we can
identify $U_\theta = D - \{0, q_\theta \}$ where $D$ was a
neighborhood of $(0,0)$ in $F_0$. This gives an identification of
$W_\A^{-1}(q_\theta)$ with $F_0 - \{0, q_\theta\}$. Utilizing such
an identification, we will explicitly write a map $M_\theta :
W_\A^{-1}(q_0) \to W_\A^{-1}(q_\theta)$ symplectically isotopic to
the parallel transport map along the arc $\gamma (t) =
 \varepsilon e^{i t\theta } = q_{t\theta}$. To define $M_\theta$, let $h:
\R_{\geq 0} \to \R_{\geq 0}$ be a smooth non-increasing function
which satisfies
\begin{displaymath} h(x) = \left\{
\begin{array}{cc} 1 & for \hbox{ } x < 2 \varepsilon \\ 0 & for \hbox{ } x >
3 \varepsilon
\end{array} \right. \end{displaymath}
Then we define \begin{displaymath} M_\theta (w) = \left\{
\begin{array}{cc} e^{h(|w|)\theta i} w & for \hbox{ } w \in U_\theta
\\ w & otherwise \end{array} \right. \end{displaymath}

It is easy to see that $M_{2\pi}: W^{-1}_\A(\varepsilon ) \to
W^{-1}_\A ( \varepsilon)$ gives a full Dehn twist about a circle
encompassing the two ends and yields the isotopy class of the
identity if these ends are included.

On the other hand, we have that multiplying a fiber $W^{-1}_\A (q)$
by $\zeta^k$ will give a map to the fiber $W^{-1}_\A (\zeta^k q)$ as
well. We arrange the parameterization of $U_\theta$ in such a way
that this map is symplectically equivalent to multiplication by
$\zeta^k$. This, together with Proposition 11 yields the following
proposition.

\begin{prop} For $\varepsilon$ sufficiently small and $q = \varepsilon e^{i \theta}$, the fiber $\wq$ is symplectically isotopic to the curve
$F_0 - \{ 0, q \}$ where $0$ and $q$ are regarded as elements of the
disc neighborhood $D$. Letting $U_\theta = D - \{ 0, q\}$, the
monodromy map along an arc from $q$ to $e^{i \varphi}q$ is
symplectically isotopic $M_\varphi: \wq \to W^{-1}_\A(e^{i\varphi
}q)$ and the map from $\wq$ to $W^{-1}_\A(\zeta^k q)$ induced by
multiplication is given by the action of $\zeta^k$ on the fiber
$F_0$ (i.e. the action on the fibers of the branched cover).
\end{prop}

The motivation for describing this parallel transport map is to
obtain a description for the vanishing cycles of $W_\A$ on
$(\C^*)^2$. Unlike the situation in $\C^2$, and indeed for all
situations involving partial Lefschetz fibrations, we need to use
more caution in choosing our distinguished basis of paths in the
image of $W_\A$. We saw in the last subsection that if one parallel
transports the vanishing cycle to the origin, then the vanishing
cycle will "fall off" the subspace $W_\A^{-1} (0)$ (i.e. it will hit
the point $(0, 0)$). Thus, any such path is not allowed. Therefore,
instead of taking the paths $\tilde \delta_k$ we will take paths
$\delta_k$ homotopic in $\C - \{0\}$ to
\begin{displaymath} g_k(t) = \left\{ \begin{array}{cc} \varepsilon
e^{4\pi k t i / (l - 1)} & for \hbox{ } 0 \leq t \leq 1/2 \\
\left( 2(1 - t)\varepsilon + (2t - 1)(l - 1) \right) e^{2\pi i /(l -
1)} & for \hbox{ } 1/2 \leq t \leq 1 \end{array} \right.
\end{displaymath} An example of such a basis is illustrated below.
\\[6pt]
\begin{center}
\begin{picture}(0,0)%
\includegraphics{mpath3.pstex}%
\end{picture}%
\setlength{\unitlength}{3947sp}%
\begingroup\makeatletter\ifx\SetFigFont\undefined%
\gdef\SetFigFont#1#2#3#4#5{%
  \reset@font\fontsize{#1}{#2pt}%
  \fontfamily{#3}\fontseries{#4}\fontshape{#5}%
  \selectfont}%
\fi\endgroup%
\begin{picture}(5683,2375)(2389,-2444)
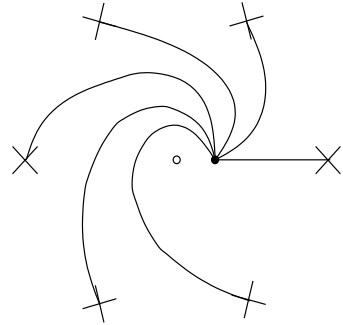
\put(6076,-2386){\makebox(0,0)[lb]{\smash{{\SetFigFont{12}{14.4}{\rmdefault}{\mddefault}{\updefault}{The distinguished basis $\{\delta_k\}$}%
}}}}
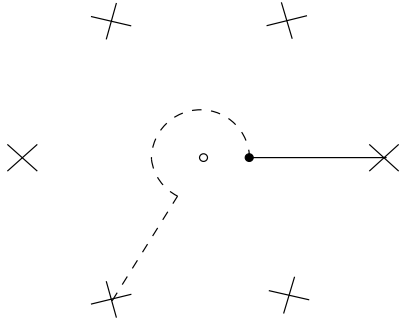
\put(2551,-2386){\makebox(0,0)[lb]{\smash{{\SetFigFont{12}{14.4}{\rmdefault}{\mddefault}{\updefault}{The curves $g_0(t)$ and $g_4(t)$}%
}}}}
\end{picture}%

\end{center}

By the previous proposition and the construction of the
distinguished basis, we can describe the vanishing cycles $L_k$
associated to $\delta_k$ in terms of the vanishing cycle $L_0$.

\begin{cor}  The vanishing cycle $L_k$ associated to $\delta_k$ is isotopic to
$M_{2  \pi k / (l - 1)}^{-1}(\zeta^k \cdot L_0 )$
\end{cor}
\begin{proof} In Proposition 12, we saw that multiplying the fiber $\we$ pointwise by $\zeta^k$
will give the fiber $W_\A^{-1}(\zeta^k \varepsilon)$ . But
multiplication by $\zeta^k$ also yields the vanishing cycle
associated to $\zeta^k \delta_0$ in $W_\A^{-1}(\zeta^k
\varepsilon)$. Parallel transporting this vanishing cycle around the
arc from $\zeta^k \varepsilon$ to $\varepsilon$ gives the monodromy
map $M_{2 k \pi / (l - 1)}^{-1}$. The concatenation of $\zeta^k
\delta_0$ and the arc from $\zeta^k \varepsilon$ to $\varepsilon$ is
just the path $g_k$ described above which is, by definition,
isotopic to $\delta_k$. This implies the result. \end{proof}

From this corollary we see that the problem of the placements of
vanishing cycles in $\we$ is reduced to understanding the placement
of the vanishing cycle $L_0$ associated to $\delta_0$. In order to
understand this subspace, we will return to a matching paths
argument. In the previous subsection we saw that the matching path
of the
 curve $\tilde \delta_0$ from $0$ to $l - 1$ was a line segment from $0$ to $l^{l - 1}$.
Utilizing the arguments given to show this, one can see that if
instead one takes the path $\delta_0$ from $\varepsilon$ to $l - 1$
then the matching path will again be a straight line segment
connecting the two positive real roots of $h_\varepsilon (x)$.
Unlike the previous situation, the critical points of $f$ restricted
to $W^{-1}(\varepsilon)$ are all Morse. Thus, the component of the
inverse via $f|_{W_\A^{-1}(\varepsilon )}$ of this segment which
contains the critical points of $f$ will give the vanishing cycle
$L_0$ in $W_\A^{-1}(\varepsilon )$. From Proposition 12 and the
above observations, we see that the topology of the fiber
$W_\A^{-1}(\varepsilon)$ will be the same as $F_0$ minus two points.
These two points appear as $0$ and $q_0$ in the coordinate chart
$U_0$ detailed above. As $\varepsilon$ approaches zero, these points
contract and the vanishing cycle $L_0$ associated to $\tilde
\delta_0$ admits the description given at the end of the previous
subsection. The following proposition describes the position of the
perturbed vanishing cycle $L_0$ associated to $\delta_0$ in
$W_\A^{-1}(\varepsilon)$.

\begin{prop} Let $\vec s$ be the line segment in $U_0$ connecting
$0$ to $q_0$. Then, after a reparametrization of $U_0$ if necessary,
the vanishing cycle $L_0$ obtained from the path $\delta_0$
intersects $\vec s$ in exactly one point.
\end{prop}

\begin{proof} To prove this proposition, we return to the
matching paths argument. In the fiber $W^{-1}_\A(\varepsilon)$ we
have seen above that the matching path of $\delta_0$ is the line
segment connecting the two positive real roots of $h_\varepsilon
(x)$. Let us label the roots of $h_\varepsilon (x)$,  $\{ r_0, r_1 ,
\ldots, r_{l - 1}\}$ where $r_0$ and $r_1$ are the two positive real
roots with $r_0 < r_1$. We recall that these are the critical values
of the map $f$ restricted to the fiber $W^{-1}_\A(\varepsilon )$. By
arguments provided in the last subsection, we have that $\varepsilon
< r_0$ and $f(W^{-1}_\A (\varepsilon )) = \C - \{\varepsilon\}$. By
the continuity of the roots of $h_q (x)$ with respect to $q$, we can
assume that $r_i$ is close to $l^{l/(l - 1)}\zeta^{i - 1}$ for $1
\leq i \leq l - 1$ and $r_0$ is close to $\varepsilon$. We start by
examining the monodromy of $f|_{W^{-1}_\A (\varepsilon ) }$ about
the critical values $r_i$ and the missing point $\varepsilon$.

Recall that $f(z_0, z_1) = a_0z_0 + a_1 z_1$ so if $f(z_0, z_1) = s$
then $z_0 = (s - a_1 z_1 ) / a_0$. Thus for any $q$ and $s \in \C$
$$f|_{F_q}^{-1} (s) = \{((s - a_1 z_1 ) / a_0 , z_1) \in F_q : C z_1^{a_1} (s - a_1
z_1)^{a_0} = s - q\}$$ In particular, for a regular value $s$ we
have the fiber of $f$ restricted to $\wq$ contains $l$ points so
that $f|_\wq$ is an $l$-fold branched cover. As was mentioned
before, the critical values of $f$ restricted to a fiber in
$(\C^*)^2$ are non-degenerate except for critical fibers. Thus $f$
restricted to $\we$ is an $l$-fold branched cover of $\C -
\{\varepsilon\}$ with ordinary double points over the critical
values $r_i$. To find the monodromy around the point $\varepsilon$
we recall that there are two points $(0, \varepsilon /a_1)$,
$(\varepsilon /a_0, 0)$ in the partial compactification
$F_\varepsilon$ which have image $\varepsilon$ under $f$. From the
above characterization of the fibers we see that
$$f|_{F_\varepsilon }^{-1} ( \varepsilon ) = \{((\varepsilon - a_1 z_1 ) / a_0, z_1) \in F_q : C z_1^{a_1} (\varepsilon - a_1
z_1)^{a_0} = 0\}$$ Therefore the point $(0, \varepsilon /a_1)$ has
ramification $a_0$ and the point $(\varepsilon /a_0, 0)$ has
ramification $a_1$. So, as a permutation on the fibers, the
monodromy around $\varepsilon$ can be written as a product of
disjoint cycles $\tau_0$ and $ \tau_1$ of order $a_0$ and $a_1$
respectively.

Now let us examine the monodromy around a curve $\gamma$ which
encloses $r_0$ and $\varepsilon$ and no other $r_i$. Perturbing the
fiber $\we$ to the degenerate fiber $W^{-1}_\A (0)$, we see that as
a permutation, the monodromy will be the same as that around the
value $0$ for the degenerate fiber. But $f|_{F_0}^{-1}(0) = (0, 0)$
so the ramification of $f|_{F_0}$ at $(0, 0)$ must be $l$. This
implies that the monodromy of $\gamma$ is a cycle of order $l$. Now,
$\gamma$ is also the composition of the monodromy around
$\varepsilon$ and the monodromy around $r_0$ which is a
transposition $\sigma$ (since $r_0$ is an ordinary double point). So
if $s$ is a regular value with $\varepsilon < s < r_0$, and
$f|_\we^{-1} (s) = \{p_1, \ldots, p_l\}$ labels the fiber in such a
way that $\tau_0 = (p_1 p_2 \cdots p_{a_0})$ and $\tau_1 = (p_{a_0 +
1} \cdots p_l)$ then we see that, as a permutation, the monodromy
around $r_0$ must be $\sigma = (p_j p_k)$ where $1 \leq j \leq a_0$
and $a_0 + 1 \leq k \leq l$.

As $0$ and $q_0$ tend towards each other in $U_0 \subset
W_\A^{-1}(\varepsilon)$ as $\varepsilon$ tends to zero, for
sufficiently small $\varepsilon$, we have that $f(\vec s )$ is
contained in a small neighborhood of $\varepsilon$, say the interior
$U$ of $\gamma$ (recall $\gamma$ is a curve going around
$\varepsilon $ and $r_0$ in the image of $f|_{\we}$). In particular,
$f(\vec s)$ stays sufficiently far away from $r_i$ for $i \not= 0$
as these values converge to the non-zero roots of $h_0(x)$. Now, in
the domain $U_0$ we can find a curve $C$ connecting $(0, \varepsilon
/a_1)$ to $(\varepsilon /a_0, 0)$. This is simply the curve which
exits $(0, \varepsilon /a_1)$ along the fiber $p_j$, tends towards
the double point over $r_0$ and loops back along the $p_k$ fiber
towards $(\varepsilon /a_0, 0)$. It is clear from the description of
the monodromy permutations above that any other curve connecting
$(0, \varepsilon /a_1)$ to $(\varepsilon /a_0, 0)$ in $U_0$ without
self intersection must be isotopic to $C$ relative boundary points.
As $\vec s$ connects $(0, \varepsilon /a_1)$ to $(\varepsilon /a_0,
0)$ and is contained in $U_0$ we see that $\vec s$ is isotopic to
$C$. Thus, reparameterizing $U_0$ if necessary we can assume that
$\vec s = C$. Finally, recall that the vanishing cycle $L_0$
associated to $\delta_0$ has $f(L_0) = \overrightarrow{r_0 r_1}$ and
must therefore contain the double point over $r_0$. Also, by the
above description, $f(C) = \overrightarrow{\varepsilon r_0 }$. This
implies $L_0 \cap \vec s$ is precisely this double point yielding
the result.
\end{proof}

The picture below shows the perturbed vanishing cycle $L_0$ in
contrast to the degenerate case.
\\[6pt]
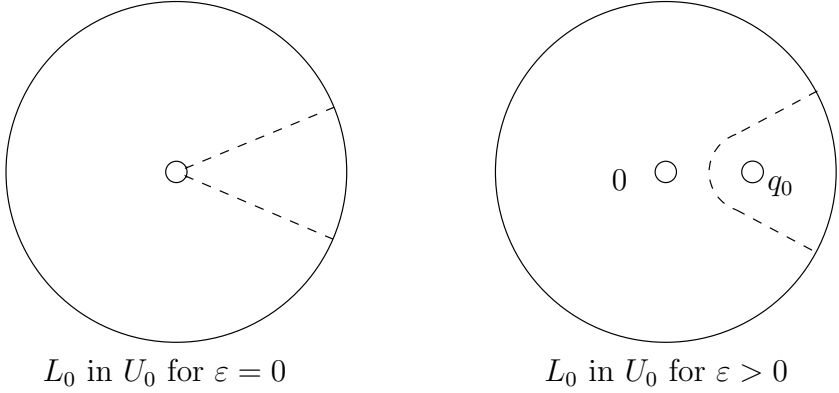
\begin{picture}(0,0)%
\includegraphics{deform2.pstex}%
\end{picture}%
\setlength{\unitlength}{3947sp}%
\begingroup\makeatletter\ifx\SetFigFont\undefined%
\gdef\SetFigFont#1#2#3#4#5{%
  \reset@font\fontsize{#1}{#2pt}%
  \fontfamily{#3}\fontseries{#4}\fontshape{#5}%
  \selectfont}%
\fi\endgroup%
\begin{picture}(5411,2449)(2026,-2819)
\put(5576,-2761){\makebox(0,0)[lb]{\smash{{\SetFigFont{12}{14.4}{\rmdefault}{\mddefault}{\updefault}{$L_0$ in $U_0$ for $\varepsilon > 0$}%
}}}}
\put(6001,-1561){\makebox(0,0)[lb]{\smash{{\SetFigFont{12}{14.4}{\rmdefault}{\mddefault}{\updefault}{$0$}%
}}}}
\put(6976,-1561){\makebox(0,0)[lb]{\smash{{\SetFigFont{12}{14.4}{\rmdefault}{\mddefault}{\updefault}{$q_0$}%
}}}}
\put(2426,-2761){\makebox(0,0)[lb]{\smash{{\SetFigFont{12}{14.4}{\rmdefault}{\mddefault}{\updefault}{$L_0$ in $U_0$ for $\varepsilon = 0$}%
}}}}
\end{picture}%
\\[6pt]

Combining this proposition and the corollary, we obtain a complete
picture of the vanishing cycles $L_k$ in $\we$. Let us summarize
this picture. The fiber $\we$ is symplectically equivalent $F_0 -
\{0, \varepsilon\}$ which admits an $(l - 1)$-fold cover to $\p^1 -
\{0, 1 , \infty\}$ except at the fiber containing $\varepsilon$
where there is a point missing. Outside the neighborhood $U_0$ of $0
\in F_0$, the description of the vanishing cycles $L_k$ stays the
same as in the previous subsection, i.e. they are $l - 1$ lifts of
the loop around $\p^1 - \{1, \infty\}$. Inside $U_0$ we can describe
$L_0$ as a curve which enters $U_0$ at the same points as it had
entered $D$ but now passes once through the curve connecting $0$ to
$q$ (one may note that there are several such curves, but utilizing
the smoothness of parallel transport of $W_\A$ in $\C^2$, we are
restricted to the picture above). Now all other vanishing cycles
$L_k \subset \we$ are obtained in $U_0$ by multiplying $L_0$ by
$\zeta^k$ and then acting on $U_{2 \pi k /(l - 1)}$ by the monodromy
map $M_{2 \pi k / (l - 1)}^{-1}$. The picture below shows the
vanishing cycles $L_0$ and $L_3$ in $U_0$ for the cases $\A = (1,
4)$ and $\A = (2, 3)$.

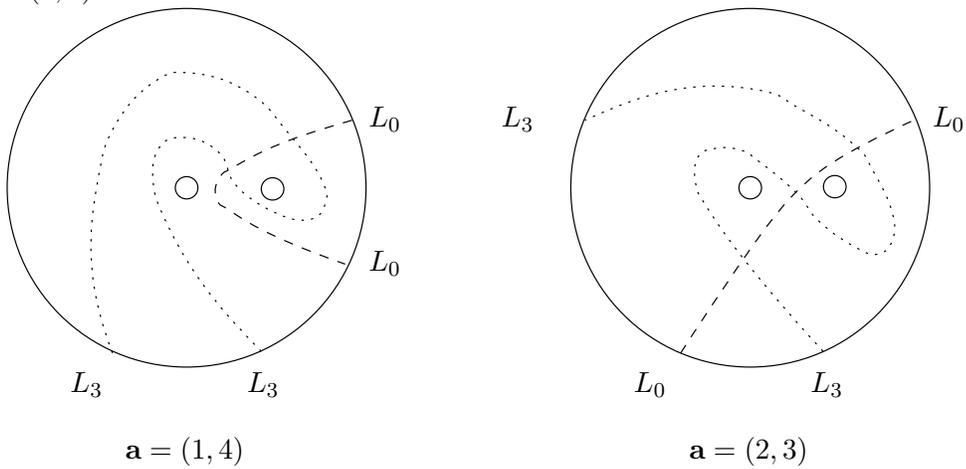
\begin{picture}(0,0)%
\includegraphics{deform3.pstex}%
\end{picture}%
\setlength{\unitlength}{3947sp}%
\begingroup\makeatletter\ifx\SetFigFont\undefined%
\gdef\SetFigFont#1#2#3#4#5{%
  \reset@font\fontsize{#1}{#2pt}%
  \fontfamily{#3}\fontseries{#4}\fontshape{#5}%
  \selectfont}%
\fi\endgroup%
\begin{picture}(6287,2887)(4145,-3268)
\put(7259,-1133){\makebox(0,0)[lb]{\smash{{\SetFigFont{11}{13.2}{\rmdefault}{\mddefault}{\updefault}{$L_3$}%
}}}}
\put(4897,-3216){\makebox(0,0)[lb]{\smash{{\SetFigFont{11}{13.2}{\rmdefault}{\mddefault}{\updefault}{$\A = (1, 4)$}%
}}}}
\put(8438,-3216){\makebox(0,0)[lb]{\smash{{\SetFigFont{11}{13.2}{\rmdefault}{\mddefault}{\updefault}{$\A = (2, 3)$}%
}}}}
\put(6424,-1133){\makebox(0,0)[lb]{\smash{{\SetFigFont{11}{13.2}{\rmdefault}{\mddefault}{\updefault}{$L_0$}%
}}}}
\put(6424,-2037){\makebox(0,0)[lb]{\smash{{\SetFigFont{11}{13.2}{\rmdefault}{\mddefault}{\updefault}{$L_0$}%
}}}}
\put(4550,-2801){\makebox(0,0)[lb]{\smash{{\SetFigFont{11}{13.2}{\rmdefault}{\mddefault}{\updefault}{$L_3$}%
}}}}
\put(5660,-2801){\makebox(0,0)[lb]{\smash{{\SetFigFont{11}{13.2}{\rmdefault}{\mddefault}{\updefault}{$L_3$}%
}}}}
\put(9966,-1133){\makebox(0,0)[lb]{\smash{{\SetFigFont{11}{13.2}{\rmdefault}{\mddefault}{\updefault}{$L_0$}%
}}}}
\put(9201,-2801){\makebox(0,0)[lb]{\smash{{\SetFigFont{11}{13.2}{\rmdefault}{\mddefault}{\updefault}{$L_3$}%
}}}}
\put(8091,-2801){\makebox(0,0)[lb]{\smash{{\SetFigFont{11}{13.2}{\rmdefault}{\mddefault}{\updefault}{$L_0$}%
}}}}
\end{picture}%

In what follows, we will give an alternative representation of the
vanishing cycles that yields the equivalence of $\D (Fuk (W_\A, \C -
\D_\epsilon, a , \{\delta_k\})$ with the category $\T \subset \D
(Coh (\blo ))$.

\subsection{The derived category of the partial Lefschetz fibration
$W_\A$}

From the previous subsection, we see that all intersection points of
the vanishing cycles $L_k$ occur in $U_0$. It is also clear that any
holomorphic disc occurring in the moduli spaces defining products in
the Fukaya category are also contained in $U_0$; otherwise, upon
projection to $\p^1 - U_0 - \{[1:0:0] , [0:1:0]\}$ we would have
that the image of $L_0$ is contractible. So in order to obtain the
Fukaya category it suffices to consider the subspace $U_0$ and its
intersections with the vanishing cycles $L_k$. Now, all of the data
which gives the directed Fukaya category in our situation is
invariant up to scaling and isotopy of the vanishing cycles. So
instead of using $U_0$ we can take the space $D_e - \{0, 1\}$ by
scaling by $\varepsilon^{-1}$ and fixing the intersection points of
the vanishing cycles with the boundary. In the representation we
will give we will take the logarithm $4(l - 1) / 2\pi  \log (z)$ of
this setup and examine the images of the vanishing cycles in $\C - 4
(l - 1) i \Z$. Recall that at the end of subsection 4.1 we described
the ordering of the vanishing cycles in $D$ in terms of the lines
$\ell_k = \{r e^{(2 k + 1) \pi i / 2(l - 1)}: 0 \leq r \leq
\epsilon\}$ and that $L_k \cap D = \ell_{2k} \cup \ell_{1 - 2a_0 +
2k}$. Thus the intersection of $L_k$ with the boundary of $D_e$ will
be $e \zeta^{(2k + 1)/4}$ and $e \zeta^{(4k - 4a_0 + 3)/4}$. We will
define a set of points in $\C - 4(l - 1)i \Z$ as follows:
\begin{displaymath} \begin{array}{ccc} P_{k - } & = & (2k + 1 - 2( l - 1))i \\ P_{k + } & = & (2k + 1)i \\
 Q_{k - } & = & 1 + (4k + 1 - 2(l - 1))i \\ Q_{k + }& = &Q_{k -} +
 (4a_0 - 2)i \end{array} \end{displaymath}

Now let $\vec s_{k \pm}$ be the class of line segments connecting
$P_{k \pm} + 4(l - 1)a i$ to $Q_{k \pm} + 4 (l - 1) a i$ for all
integers $a$. Let $C_k = \{z: |z - P_{k -} + (l - 1)i| = l - 1
\hbox{ }, \hbox{ } Im (z) \leq 0\} + 4(l - 1)a i$ be the set of half
circles connecting $P_{k + } + 4(l - 1) a i$ to $P_{k -} + 4(l - 1)a
i$ in the negative real half plane. Finally, take
$$\lag_k = \vec s_{k +} \cup \vec s_{k -} \cup C_k$$ Then after
re-examining Corollary 4, the comments of the last subsection and
the definition of the monodromy map, we have the following
proposition.

\begin{prop} The piecewise smooth curves $\exp(\frac{2\pi }{ 4(l -
1)} \lag_k)$ are isotopic to the vanishing cycles $L_k \cap U_0$.
\end{prop}

By isotopic, we will mean that we have smoothed the corners on the
imaginary line and that we have rotated the image to agree with the
boundary conditions. We illustrate this set-up for our two examples.
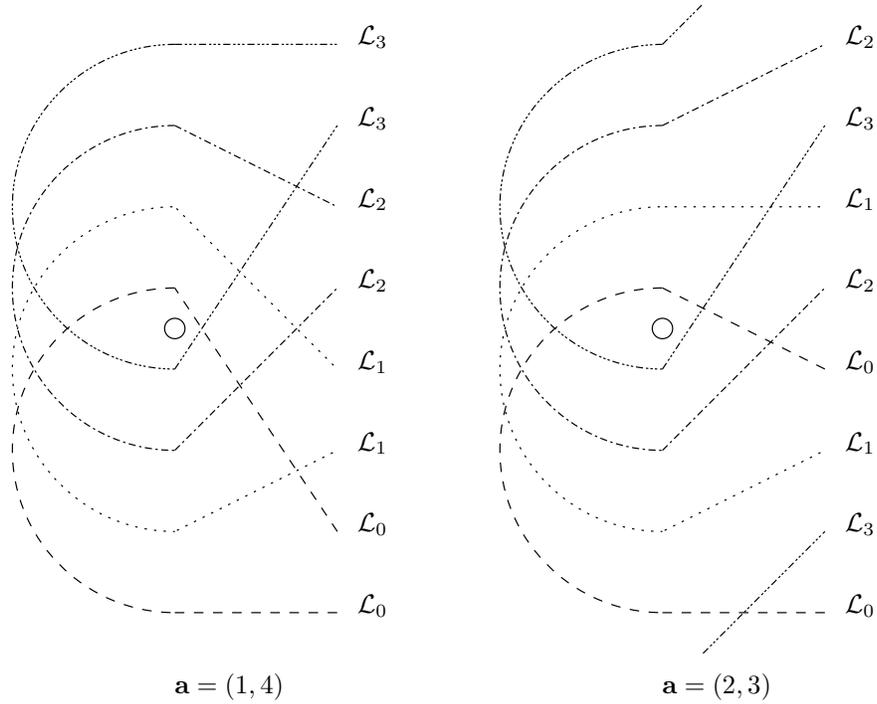
\begin{center}
\begin{picture}(0,0)%
\includegraphics{wpfuk1.pstex}%
\end{picture}%
\setlength{\unitlength}{3947sp}%
\begingroup\makeatletter\ifx\SetFigFont\undefined%
\gdef\SetFigFont#1#2#3#4#5{%
  \reset@font\fontsize{#1}{#2pt}%
  \fontfamily{#3}\fontseries{#4}\fontshape{#5}%
  \selectfont}%
\fi\endgroup%
\begin{picture}(5241,4399)(2369,-4748)
\put(6462,-4699){\makebox(0,0)[lb]{\smash{{\SetFigFont{10}{12.0}{\rmdefault}{\mddefault}{\updefault}{$\A = (2, 3)$}%
}}}}
\put(4548,-616){\makebox(0,0)[lb]{\smash{{\SetFigFont{10}{12.0}{\rmdefault}{\mddefault}{\updefault}{$\lag_3$}%
}}}}
\put(4548,-1637){\makebox(0,0)[lb]{\smash{{\SetFigFont{10}{12.0}{\rmdefault}{\mddefault}{\updefault}{$\lag_2$}%
}}}}
\put(4548,-2147){\makebox(0,0)[lb]{\smash{{\SetFigFont{10}{12.0}{\rmdefault}{\mddefault}{\updefault}{$\lag_2$}%
}}}}
\put(4548,-2657){\makebox(0,0)[lb]{\smash{{\SetFigFont{10}{12.0}{\rmdefault}{\mddefault}{\updefault}{$\lag_1$}%
}}}}
\put(4548,-3168){\makebox(0,0)[lb]{\smash{{\SetFigFont{10}{12.0}{\rmdefault}{\mddefault}{\updefault}{$\lag_1$}%
}}}}
\put(4548,-3678){\makebox(0,0)[lb]{\smash{{\SetFigFont{10}{12.0}{\rmdefault}{\mddefault}{\updefault}{$\lag_0$}%
}}}}
\put(4548,-4188){\makebox(0,0)[lb]{\smash{{\SetFigFont{10}{12.0}{\rmdefault}{\mddefault}{\updefault}{$\lag_0$}%
}}}}
\put(4548,-1126){\makebox(0,0)[lb]{\smash{{\SetFigFont{10}{12.0}{\rmdefault}{\mddefault}{\updefault}{$\lag_3$}%
}}}}
\put(7610,-2147){\makebox(0,0)[lb]{\smash{{\SetFigFont{10}{12.0}{\rmdefault}{\mddefault}{\updefault}{$\lag_2$}%
}}}}
\put(7610,-3168){\makebox(0,0)[lb]{\smash{{\SetFigFont{10}{12.0}{\rmdefault}{\mddefault}{\updefault}{$\lag_1$}%
}}}}
\put(7610,-4188){\makebox(0,0)[lb]{\smash{{\SetFigFont{10}{12.0}{\rmdefault}{\mddefault}{\updefault}{$\lag_0$}%
}}}}
\put(7610,-1126){\makebox(0,0)[lb]{\smash{{\SetFigFont{10}{12.0}{\rmdefault}{\mddefault}{\updefault}{$\lag_3$}%
}}}}
\put(7610,-2657){\makebox(0,0)[lb]{\smash{{\SetFigFont{10}{12.0}{\rmdefault}{\mddefault}{\updefault}{$\lag_0$}%
}}}}
\put(7610,-1637){\makebox(0,0)[lb]{\smash{{\SetFigFont{10}{12.0}{\rmdefault}{\mddefault}{\updefault}{$\lag_1$}%
}}}}
\put(7610,-616){\makebox(0,0)[lb]{\smash{{\SetFigFont{10}{12.0}{\rmdefault}{\mddefault}{\updefault}{$\lag_2$}%
}}}}
\put(7610,-3678){\makebox(0,0)[lb]{\smash{{\SetFigFont{10}{12.0}{\rmdefault}{\mddefault}{\updefault}{$\lag_3$}%
}}}}
\put(3399,-4699){\makebox(0,0)[lb]{\smash{{\SetFigFont{10}{12.0}{\rmdefault}{\mddefault}{\updefault}{$\A = (1, 4)$}%
}}}}
\end{picture}%
\end{center}
The advantage of considering this representation of the vanishing
cycles is that the morphisms $(L_i, L_k)$ have a natural
decomposition which exhibits the isomorphism with the mirror
category. We will now state this as a theorem modulo issues related
to grading.

\begin{thm} Let $V$ be a graded vector space generated by $e_0$ and $e_1$ with weights $a_0$ and $a_1$.
Define $\psi_{j,k} : (L_j, L_k) \to \oplus_{t \leq k - j}
\bigwedge^{\bullet t} V$ as
\begin{displaymath}  \begin{array}{ccc} \psi_{j, k} (\exp(\frac{2\pi }{ 4(l -
1)}C_j \cap C_k)) & = & 1
\\ \psi_{j,k} (\exp(\frac{2\pi }{ 4(l -
1)}\vec s_{j +} \cap \vec s_{k -})) & = & e_0
\\ \psi_{j, k} (\exp(\frac{2\pi }{ 4(l -
1)}\vec s_{k +} \cap \vec s_{j -})) & = & e_1 \end{array}
\end{displaymath}
Then $\psi_{j, k}$ is an isomorphism that commutes with composition
in $Fuk(W_\A, \C - D_\epsilon, a, \{\delta_i\})$. Furthermore, all
products $m_i$ vanish for $i \not= 2$.
\end{thm}
The proof of this theorem will occupy the rest of this subsection.
To see why this theorem is technically difficult, one can examine
the figures above and find many embedded polygons with four or more
edges that a priori could contribute to higher products. One point
that will be proven is that all such polygons do not have the
appropriate ordering of intersection points along the boundary to
contribute to $m_i$. Although this proof is lengthy, the techniques
are completely elementary.

\begin{proof} We will prove
this theorem in several steps.
\\[6pt]
{\it Step 1}: The map $\psi_{j, k}$ is a well defined isomorphism of vector spaces. \\[6pt]

Recall that $\bigwedge^{\bullet t}V$ consists of the elements of the
exterior algebra with weighted degree equal to $t$. So $\oplus_{t
\leq k - j} \bigwedge^{\bullet t} V$ will always contain the
identity element, it will contain the generator $e_0$ iff $a_0 \leq
k - j$ and it will contain the generator $e_1$ iff $a_1 \leq k - j$.
Note that $e_0 \wedge e_1$ will not be contained in this space as $k
- j \leq l - 2$. Now, it is clear that $C_j \cap C_k$ will always
contain exactly one intersection point modulo $4(l - 1)i \Z$. So to
complete this step we must show that, modulo $4(l - 1)i \Z$, $\vec
s_{j +} \cap \vec s_{k -}$ contains $1$ element iff $a_0 \leq k - j$
and is empty otherwise and $ \vec s_{k +} \cap \vec s_{j -}$
contains $1$ element iff $a_1 \leq k - j$ and is empty otherwise. We
use brute force and calculate the equations for the lines containing
the segments $\vec s_{k \pm}$. Using cartesian coordinates we write
the equations as follows
\begin{displaymath} \begin{array}{ccccc} \vec s_{k -} & \hbox{ } & y
& = & 2k x + 2k + 1 - 2(l - 1) + 4(l - 1)b \\ \vec s_{j + } & \hbox{
} & y & = & (2j - 2l + 4a_0) x + (2k + 1) + 4 (l - 1)c \end{array}
\end{displaymath} where $b$ and $c$ are integers. Then $(x, y)
\in \vec s_{j +} \cap \vec s_{k -}$ iff $0 < x < 1$ and there exists
an integer $d = b - c$ such that
\begin{equation} (k - j + l - 2a_0) x  =  (j - k ) + (l - 1)(2d
+ 1)
\end{equation}

where we have used the above equations to solve for $x$.\\[4pt] {\it
Case 1}: $\vec s_{j +} \cap \vec s_{k -}$ \\[4pt] We have that $k - j > 0$ and $l - 2a_0
= a_1 + a_0 - 2a_0 = a_1 - a_0 \geq 0$ by our convention that $a_0
\leq a_1$. Thus $k - j + l - 2a_0 > 0$ and dividing equation ($3$)
by this term, the inequality $0 < x < 1$ gives $0  <  (j - k) + (l -
1)(2d + 1)  < (k - j + l - 2a_0)$ which in turn gives
\begin{displaymath}
\begin{array}{ccccc} k - j & < & (l - 1)(2d + 1) & < & 2(k - j) + l - 2a_0
\end{array}
\end{displaymath}
Now, we have $0 < k - j < l - 1$ and $l - 2a_0 = a_1 - a_0 < l - 1$
which implies \begin{displaymath} \begin{array}{ccccc} 0 & < & (l -
1)(2d + 1) & < & 3(l - 1)
\end{array} \end{displaymath} and therefore $d = 0$. So for $k > j$ there is an intersection point
of $\vec s_{j +} \cap \vec s_{j -}$ iff the inequality $k - j < l -
1 < 2(k - j) + l - 2a_0$ is satisfied. The left hand side of this
inequality is always true, while the right hand side gives $2a_0 - 1
< 2(k - j)$ which is the case iff $a_0 \leq k - j$. The fact that
this intersection point is unique mod $4(l - 1)i \Z$ follows from
the fact that the integer $d$ was determined.
\\[4pt] {\it Case 2}: $\vec s_{k +} \cap \vec s_{j -}$ \\[4pt]
For this case we have more options to consider. Our equation ($3$)
is now:
\begin{equation} (j - k + l - 2a_0) x  =  (k - j ) + (l - 1)(2d
+ 1)
\end{equation}

We start with the possibility that $j - k + l - 2a_0 > 0$. This will
imply that $a_1 > a_1 - a_0 = l - 2a_0 > k - j$ so that we should
expect no intersection points. Indeed, proceeding as in case $1$ we
obtain the inequality:
\begin{displaymath}
\begin{array}{ccccc} j - k & < & (l - 1)(2d + 1) & < & 2(j - k) + l - 2a_0
\end{array} \end{displaymath}
Here $- (l - 1) < j - k < 0$ and $l - 2a_0 < l - 1$ which implies $-
(l - 1)  <  (l - 1)(2d + 1)  <   l - 1$ and therefore $- 1 < d < 0$.
But $d$ is required to be integral so there exists no intersection
points as expected.

The second option to consider is $j - k + l - 2a_0 = 0$. This will
occur iff the slope of $\vec s_{k +}$ equals that of $\vec s_{j -}$
so that no intersection points occur. But again we have $a_1 > a_1 -
a_0 = l - 2a_0 = k - j$ so we expect no intersection points for this
case.

The final option for this case is that $j - k + 1 - 2a_0 < 0$. Here,
there is an intersection point in $ \vec s_{k +} \cap \vec s_{j -}$
iff there exists an integer $d$ satisfying the inequality
\begin{displaymath}
\begin{array}{ccccc} j - k & > & (l - 1)(2d + 1) & > & 2(j - k) + l - 2a_0
\end{array} \end{displaymath}
Again we have that $0 > j - k > - (l - 1)$ which implies $2(j - k) +
l - 2a_0 > -2(l - 1)$ so that $0 > 2d + 1 > - 2$ which implies $d =
- 1$. Putting this into the inequality above we have
\begin{displaymath}
\begin{array}{ccccc} k - j & < & (l - 1) & < & 2(k - j) - l + 2a_0
\end{array} \end{displaymath} The left hand side is always satisfied
while the right hand side gives $2l - 2a_0 - 1 < 2(k - j)$ or $2a_1
- 1 < 2(k - j)$. This is satisfied iff $a_1 \leq k - j$ as desired.
Again, the uniqueness of such an intersection point mod $4(l - 1)i
\Z$ follows from the uniqueness of $d$.
\\[4pt] {\it Step 2}: Holomorphic discs which intersect the negative
real half plane.
\\[4pt] For what follows, we will give an orientation to each
$\lag_i$ by assuming it flows from $Q_{i +}$ to $Q_{i -}$. Given a
holomorphic disc $u \in \M_0 (p_{i_0}, \ldots, p_{i_r})$ we observe
that the image of $u$ in $\C$ makes an angle at $ p_{i_k}$ less than
$180$ degrees. This is due to the fact that $u \in \M_0 (p_{i_0},
\ldots, p_{i_r})$ instead of simply in $\M (p_{i_0}, \ldots,
p_{i_r})$. We will represent the disc $u$ in the following way: Let
$C_{i}(\pm)$ denote the circle $C_i$ with the given orientation or
the opposite orientation depending on the $\pm$. Similarly let $s_{i
\pm}(\pm)$ denote the oriented line segments. If $G_{k}$ is the
oriented curve on $\lag_{i_k}$ which is the piece of the boundary of
$u$ connecting $p_{i_{k-1}}$ to $p_{i_k}$ we represent $G_k$ as the
word $X_{k_1} \cdots X_{k_t}$ where $X_{k_j}$ is one of the letters
$C_{i_k}(\pm)$ or $s_{i_k}(\pm)$, each $X_{k_j}$ intersects $G_{k}$
and the string of $X_{k_j}$'s is ordered so that their concatenation
gives the minimal curve in $\lag_{i_k}$ that intersects $G_k$. We
will label this word $Y_k$ and represent $u$ as the total word $Y_u
= \prod_{j = 0}^r Y_j$. An important observation is that the
subscripts of the letters in the word $Y_u$ are non-decreasing from
left to right which follows from the definition of $\M_0 (p_{i_0},
\ldots, p_{i_r})$.

Let us start by examining the case where the combination $C_i (\pm)
C_j (\pm ) C_k (\pm )$ occurs in the word $Y_u$. One can see that
the orientation of $u$ forces the first part of this combination to
be $C_i(\pm ) C_j(\mp)$. Now, as one follows $C_j$ along its
orientation, the intersection points $C_n \cap C_j$ along $C_j$ are
decreasing (i.e. $n$ is decreasing). Thus if the first part is
$C_i(+) C_j(-)$ then $k > i$. On the other hand, if the first part
is $C_i(-) C_j(+)$ then $k < i$, but such a situation is not
allowed, as the ordering of the vanishing cycles in the word must be
increasing. Thus any such combination occurring in $Y_u$ must be of
the form $C_i(+) C_j(-) C_k (+ )$. One can see by the same argument
that there is no $C_s( \pm)$ for which $C_i(+) C_j(-) C_k (+ )
C_s(\pm)$ or $C_s(\pm ) C_i(-) C_j(+) C_k (- )$ occurs in $Y_u$. But
observe that there is a holomorphic triangle $\tilde u \in \M_0 (C_i
\cap C_j, C_j \cap C_k , C_k \cap C_i)$ which has the word $C_i(-)
C_j(+) C_k (- )$. As no other intersection points for $u$
consecutively occur before or after $C_i \cap C_j$ and $C_j \cap
C_k$ in the negative real half plane we have that $im({\tilde u})
\subset im (u)$. But the boundary of $\tilde u$ is also contained in
the boundary of $u$ which implies that $u = {\tilde u}$. It is clear
that such $\tilde u$ defines the product $m_2$ for those
intersection points sent to $1 \in \bigwedge^{\bullet \bullet} V$
and their composition is compatible.

Now let us assume that $C_i(\pm)C_j(\mp)$ occurs in $Y_u$ and $u$ is
not contained in the negative real half plane. Let us take the case
$C_i(+) C_j(-)$ as the alternative case has a similar argument. In
this case we have that $P_{i +}, P_{j +} \in im (u)$ so that $im
(u)$ intersects the positive real half plane $H_+$. Now, the
component $K$ of $im (u) \cap H_+$ containing $P_{i +}$ is a convex
set as the boundary consist of line segments with intersection
angles $< \pi$. Therefore, no $P_{k -}$  or $P_{t +}$ for $t < i $
or $t > j$ occurs in $K$ as these points are not contained in the
convex set lying above $\vec s_{i+}$ and below $\vec s_{j +}$.
Furthermore, if $P_{k +}$ occurs in the boundary of $u$ with $i < k
< j$ then the symbol $C_k(\pm)$ must occur somewhere in the word
$Y_u$. But as the indices of the letters in this word are
increasing, this is not possible. Therefore, $K$ contains $P_{i +}$
and $P_{j+}$. Taking the union $\tilde K$ of $K$ and the triangle
$P_{i +}(C_i \cap C_j)P_{j +}$ gives a holomorphic disc contained in
the image of $u$ whose boundary is also contained in the boundary of
the image of $u$. Thus, $\tilde K = im (u)$ so that any $u$ that is
not a holomorphic triangle but contains two consecutive letters
$C_i(\pm ) C_j(\mp )$ in $Y_u$ contains no other $C_k(\pm )$.

Now suppose that the first letter in $Y_u$ is $C_i (\pm )$ and the
last letter is $C_j(\pm)$ and that $im (u)$ is not contained in the
negative real half plane. Then we can see that the signs must be
$C_i(+)$ and $C_j(-)$. From the above paragraph we know that
$C_i(+)s_{i +}(+)$ is the beginning of $Y_u$ and $s_{j -}(-) C_j(-)$
is the end of $Y_u$. So $P_{i +}$ and $P_{j -}$ are elements of the
boundary of $u$. If we take $K$ to be the convex set $im (u) \cap
H_+$ containing $P_{j -}$ then by the orientation of $u$, $K$ is
contained in the set of elements above $\vec s_{j -}$. Now, $K$
cannot contain $0$ as this would imply that $0 \in im (u)$ so $K
\cap \R i$ must be a line segment from $P_{j -}$ to $P_{k -}$ with
$k > j$ and $P_{k -}$ in the boundary of $u$. But this implies that
the letter $C_k(\pm)$ occurs in $Y_u$ and as $C_j(-)$ is the last
letter, it must occur before $C_j(-)$. This violates the increasing
order on the subscripts of letters in $Y_u$. Therefore, there is no
curve $u$ for which the first and last letters are semi-circles.

Our last case to consider for this step is that $C_i$ occurs in
$Y_u$ without any other consecutive letters being semi-circles
(consecutive in the cyclic sense). This will imply that there are no
intersection points lying in the negative real half plane. But then
the interior of the region bounded by $C_i$ and $\R i$ is contained
in $im (u)$ and, in particular, $0 \in im (u)$ which is not
possible. So this case does not occur.

Summing up these results we have: \\[6pt]
\indent Given $u \in \M_0 (p_{i_0}, \ldots, p_{i_r})$ there are
three possibilities for $Y_u$ \\ \indent \indent (1) $u$ is a
holomorphic triangle contained in the negative real half plane \\
\indent \indent (2) The word $Y_u$ contains the expression
$C_i(\pm)C_j(\mp)$ and all other \\ \indent \indent \indent letters
are of the type $s_{k \pm}(\pm )$ \\ \indent \indent (3) The word
$Y_u$ contains no letters of the type $C_i(\pm )$
\\[6pt]
{\it Step 3}: Holomorphic discs which intersect the positive real
half plane.
\\[4pt]

Before pursuing this case in detail, let us note that as one moves
along $\vec s_{i +}$ in the positively oriented direction, $j$
increases where $\{p_j\} = \vec s_{i +} \cap \vec s_{j -}$.
Similarly, as one moves along $\vec s_{i -}$ in the positively
oriented direction, $j$ increases where $\{p_j\} = \vec s_{i -} \cap
\vec s_{j +}$. Both of these facts can be seen from the equations
defining $\vec s_{k \pm}$. Now, if $i < j$ step $1$ gives us
$\{p_j\} = \vec s_{i -} \cap \vec s_{j +}$ iff $j - i \geq a_1$.
Now, if $j - i \geq a_1$ then there exists no $0 \leq k \leq l - 2$
for which $i - k \geq a_0$ for then $j - k \geq l$. Thus, if $l -
a_1 - 2 \geq i$ then the only intersection points $p_j$ along $\vec
s_{i -}$ must satisfy $j \geq i + a_1$. A similar argument shows
that if $i \geq a_0$ then the only
 intersection points $p_j$ along $\vec s_{i -}$ must satisfy $j
+ a_0 \leq i$. We summarize these and analogous results in the
following table:  \begin{center} \begin{tabular}{ r l} $l -2 - a_1
\geq i$ \hbox{  implies} & $j \geq i + a_1$ for $p_j \in \vec s_{i -}$ \\
$a_0 \leq i$ \hbox{  implies} & $j + a_0 \leq i$ for $p_j \in \vec
s_{i -}$ \\ $l -2 - a_0 \geq i$ \hbox{  implies} & $j \geq i + a_0$
for $p_j \in \vec s_{i +}$
\\ $a_1 \leq i$ \hbox{  implies} & $ j + a_1 \leq i$ for $p_j \in \vec s_{i +}$
\end{tabular} \end{center}

We will use this table in what follows. Finally recall that $\vec
s_{i +} \cap \vec s_{j +} = \emptyset = \vec s_{i -} \cap \vec s_{j
- }$ for all $i \not= j$.

Our first case to consider for this step is that $Y_u$ contains
three consecutive letters of the type $s_{i\pm }(\pm)$, say $s_{i
\pm}(\pm)s_{j \pm}(\pm ) s_{k \pm}(\pm)$. Utilizing the clockwise
orientation on $u$, one sees that the possible combinations of two
consecutive letter of this type are
\begin{displaymath}
\begin{array}{ccc} s_{i -}(+)s_{j +}(-) & & s_{i +}(+) s_{j -} (+)
\\ s_{i -}(-)s_{j +}(+) & & s_{i +}(-) s_{j -} (- ) \end{array}
\end{displaymath}

Let us take examine the cases where either $s_{i -} (+) s_{j +} (-)$
or $s_{i +} (-) s_{j -} (-)$ occurs. Then, as we are moving from the
intersection point $\vec s_{i \pm} \cap \vec s_{j \pm}$ along $\vec
s_{j \pm}$ against the orientation, we see that $k < i$. But as the
indices in $Y_u$ must be increasing, this is not possible, so these
cases does not occur.

Now suppose $s_{i +} (+) s_{j -} (+)$ occurs. If $a_1 \leq i$ then
by the above table, $j < i$ which contradicts increasing ordering of
indices, thus $l - 2 - a_0 \geq i$ and $j \geq i + a_0 \geq a_0$.
But again by the above table this implies that $k + a_0 \leq j$ or $
k < j$ which contradicts the increasing ordering of indices in $s_{i
+} (+) s_{j -} (+)s_{k +}(-)$. Now suppose that $s_{i -}(-)s_{j+}
(+)$ occurs in $Y_u$. If $a_0 \leq i$ then the table shows that $j <
i$ contradicting increased ordering of indices; so $l - 2 - a_1 \geq
i$ and $j \geq i + a_1 \geq a_1$. But then the table shows that $k +
a_1 \leq j$ or $k < j$ again contradicting the increasing ordering
of indices of $Y_u$. As these are the only options, we see that
there are never three consecutive letters of the type $s_{i
\pm}(\pm)$ in $Y_u$. In particular, there is no acceptable
holomorphic disc contained in the positive real half plane.

The above result, along with step 1, implies that the only
acceptable holomorphic discs intersecting $H_+$ will contain either
one point in $H_-$ and one in $H_+$, one point in $H_-$ and two in
$H_+$ or one point in $H_-$ and three in $H_+$. In each case $Y_u$
must contain the word $s_{i \pm}(\pm)C_i(\pm )C_j(\mp) s_{j
\pm}(\mp) $ where the signs are all forced by the orientation. Let
us eliminate the first and third options. For the first option we
have $Y_u = s_{i \pm}(\pm)C_i(\pm )C_j(\mp) s_{j \pm}(\mp) $ . One
quickly observes that the subscript sign for $s_{i \pm}$ is the same
as that for $s_{j \pm}$. But then $\vec s_{i \pm} \cap \vec s_{j
\pm} = \emptyset$ implying that the boundary of $u$ contains only
one intersection point. This is not allowed, so no such $u$ exists.
The third option must have a word of the type $$Y_u = s_{t
\mp}(\sigma_t)s_{i \pm}(\pm)C_i(\pm )C_j(\mp) s_{j \pm}(\mp)s_{k
\mp}(\sigma_k )$$ where $\sigma_t$ and $\sigma_k$ are the
appropriate signs. But again we have that the signs of the
subscripts for the beginning and end agree, so $\vec s_{t \mp} \cap
\vec s_{k \mp} = \emptyset$ and the curve $u$ does not close.

At this point we can affirm that $m_i = 0$ for $i \not= 2$. Indeed,
above we have shown that no $u$ has image in $H_+$, so every $u$
must have image intersecting $H_-$. By step 1 we saw that any $u$
with image intersecting $H_-$ that is not a triangle contains
precisely one point in $H_-$. By the above paragraph, we have ruled
out the case where $u$ contributes to $m_1$ and the case where $u$
contributes to $m_3$. As no three $s_{i \pm }(\pm )$ occur
consecutively in $Y_u$, we have that the only remaining case consist
of triangles with one point in $H_-$ and two in $H_+$.

Let us start by assuming $s_{i +}(+)C_i(+) C_j(-)s_{j +}(-)$ and
$s_{k -}(\pm)$ occurs in $Y_u$. If $k < i$ then $s_{k -}(\pm )$ must
occur before the word $s_{i +}(+)C_i(+) C_j(-)s_{j +}(-)$ and if $k
> j$ it must occur after. Assume $k < i$ so that $$Y_u =s_{k -}(-) s_{i +}(+)C_i(+) C_j(-)s_{j
+}(-)$$ The points $p_0, p_1, p_2$ occurring in the triangle satisfy
$\{p_0\} = \vec s_{k - } \cap \vec s_{i +}$, $\{p_1\} = C_i \cap
C_j$ and $\{p_2\} = \vec s_{j +} \cap \vec s_{k -}$. By the
definition of $\psi_{\bullet, \bullet}$ we have $\psi_{k, i}(p_0) =
e_1$, $\psi_{i, j} (p_1) = 1$ and $\psi_{k, j} (p_2) = e_1$. Thus,
if this triangle exists,  it defines composition compatible with
that in $\bigwedge^{\bullet \bullet} V$. One can see the existence
of such a triangle by taking the convex hull of $P_{i +}, P_{j +},
p_0$ and $p_2$ and taking the union with the region contained in the
triangle $p_1P_{i + } P_{j +}$.

Now suppose $k > j$, then we have $$Y_u =s_{i +}(+)C_i(+) C_j(-)s_{j
+}(-)s_{k -}(-) $$ and $\{p_0\} = C_i \cap C_j$, $\{p_1 \} = \vec
s_{j +} \cap \vec s_{k - }$ and $\{p_2\} = \vec s_{i +} \cap \vec
s_{k -}$. Again by the definition of $\psi_{\bullet, \bullet}$ we
have $\psi_{i, j}(p_0) = 1$, $\psi_{j, k} (p_1) = e_0$ and $\psi_{i,
k} (p_2) = e_0$. So again we see that the part of $m_2$ that this
triangle defines agrees via $\psi$ with multiplication in
$\bigwedge^{\bullet \bullet} V$. The existence of such triangles can
be shown in a similar way as above.

Following the procedure in the previous two paragraphs, one can take
care of the case where $s_{i -}(-)C_i(-) C_j(+)s_{j -}(+)$ occurs in
$Y_u$ to see that $\psi_{\bullet \bullet}$ commutes with
multiplication in $Fuk(W_\A , \C - D_\epsilon, a, \{\delta_i\})$.
This completes the proof.
\end{proof}

\subsection{Grading in $Fuk(W_\A , \C - D_\epsilon, a,
\{\delta_i\})$} The theorem in the preceding section showed an
equivalence between $Fuk(W_\A , \C - D_\epsilon, a, \{\delta_i\})$
and the category $\T$ associated to our exceptional collection
$\{S_k\}$ assuming we neglected grading. This subsection will finish
the proof of this equivalence by explicitly grading the vanishing
cycles $L_i$. The procedure for constructing these gradings is
fairly straightforward; however, calculating the Maslov indices
poses a technical challenge and is the most difficult part of the
argument. There are two obstructions to giving vanishing cycles
$\Z$-gradings. The first is the obstruction to lifting the
lagrangian Grassmanian $\Lambda$ of the tangent bundle to a
fiberwise universal covering $\tilde \Lambda$. Assuming this
obstruction vanishes and we have such a lift, the second obstruction
occurs as an obstruction of lifting the canonical section $\phi_i:
L_i \to \Lambda$ to ${\tilde \phi_i}: L_i \to {\tilde \Lambda}$.

Now suppose we have a complex valued $1$-form $\eta$ on the fiber
$W^{-1}(\varepsilon)$ that yields an isomorphism
$\eta_p:T_p(W_\A^{-1}(\varepsilon)) \to \C$ for every $p \in
W^{-1}(\varepsilon)$. In most applications, $\eta$ is a holomorphic
$1$-form for which the above requirement is automatic. Observe that
$\eta$ yields a fiberwise trivialization of $\Lambda$, $\eta_p:
\Lambda_p \to S^1$ via $\eta_p (\left< v \right> ) =  \arg (\eta_p
(v))$ where $v \in T_p(W_\A^{-1}(\varepsilon))$ is a generator for
the lagrangian line and $S^1$ is identified with $\R / \pi \Z$. Thus
we can form the pullback:
\begin{displaymath}
\begin{array}{ccc} {\tilde \Lambda} & {\buildrel \upsilon \over
\longrightarrow}  & \R \\ {{\scriptstyle \pi} \downarrow} & &
{\downarrow }
\\ \Lambda & {\buildrel  \arg \eta \over
\longrightarrow} &
 S^1 \end{array} \end{displaymath}

 So given any such $\eta$, the first obstruction will vanish (note
 however that different choices of $\eta$ may yield different
 bundles $\tilde \Lambda$). Let us assume for the moment that such a $1$-form has been constructed. Let $P_\varepsilon$ be the
 topological category of immersed paths in $W_\A^{-1}(\varepsilon)$.
 We will take an object of $P_\varepsilon$ to be a point $p \in
 W_\A^{-1}(\varepsilon)$ along with a tangent vector $ v
 \in
 T_p(W_\A^{-1}(\varepsilon))$ and morphisms to be immersed paths (we will systematically ignore the fact that this is only
 an $A^\infty$-category without units). Given a path $\alpha : [0, 1]
 \to W_\A^{-1}(\varepsilon)$ one has a canonical map $\phi_\alpha: [0, 1] \to
\Lambda$ and if ${\tilde \phi_\alpha}: [0, 1] \to {\tilde \Lambda}$
is any lift then define $$\xi (\alpha ) = \upsilon({\tilde
\phi_\alpha}(1)) - \upsilon ({\tilde \phi_\alpha}(0))$$ where
$\upsilon$ is the map defined in the above pullback diagram. It is
not hard to show that $\xi$ is a functor from $P_\varepsilon$ to
$\R$ where $\R$ is the additive category with one object. One can
also check that $\xi$ is constant on connected components of
morphisms (i.e. invariant under isotopies of immersed paths with
fixed tangent vectors at the endpoints). Furthermore, $\xi$ is also
constant on free isotopy classes of closed loops. Explicitly, we
have that $\xi (\alpha ) = \arg (\eta (\alpha^\prime(1))) - \arg
(\eta (\alpha^\prime (0)))$ where $\arg (\eta (\alpha^\prime (t)))
\in \R$ is continuous on $[0, 1]$.
 Using this language one can see that $L_i$
has a $\Z$-grading if and only if $\xi (L_i) = 0$ (more generally,
$L_i$ has a $\Z / n \Z$-grading iff $\xi (L_i )\equiv 0 (mod n)$).

This language is also helpful in calculating the Maslov index of an
intersection point of two vanishing cycles. Recall that a
$\Z$-graded vanishing cycle is a vanishing cycle $L_i$ along with a
lift ${\tilde \phi_i}:  L_i \to {\tilde \Lambda}$ of the canonical
section $\phi_i$. Giving a orientation to each $L_i$, we can take
${\tilde \phi_i}(p) = \arg (\eta (v_p)) \in \R$ where $v_p$ is the
positively oriented ray in $T_p L_i$ and we have chosen a particular
lift from $\R / \pi \Z$ to $\R$. If $p \in L_i \cap L_j$ with $j
> i$ then one can define the Maslov index at $p$ as $\mu (p) = -
\lfloor \frac{1}{ \pi} ( {\tilde \phi_j}(p) - {\tilde \phi_i}(p))
\rfloor$ where $\lfloor a \rfloor$ is the least integer function.
One should note that in other expositions this definition is
actually $1 - \mu (p)$, but as we have taken a counter-clockwise
orientation of our distinguished basis, we must take Maslov indices
relative the conjugate holomorphic structure.

Now suppose we have oriented each $L_i$ as in the previous
subsection. Given $p \in L_i \cap L_j$ with $j > i$ let $v_j$ and
$v_i$ be the tangent vectors at $p$ to $L_j$ and $L_i$ respectively.
Let $f_i$ and $f_j$ be oriented parameterizations of $L_i$ and $L_j$
respectively such that $f_i (0) = p = f_j(0)$. Now, give
$W^{-1}(\varepsilon)$ the orientation induced by the complex
conjugate structure (i.e. clockwise orientation). Define $or(p)$ to
be zero if $v_j \wedge v_i$ is oriented and $1$ otherwise. Given any
sufficiently small $\kappa$, we will define a path $\alpha_\kappa$
which starts along $L_j$ and ends along $L_i$ and depends on
$or(p)$. Define $\alpha_\kappa$ to be a smoothly embedded path
contained in a sufficiently small neighborhood of $p$ such that
$(\alpha_{\kappa}(0), \alpha_\kappa^\prime (0)) =
(f_j((-1)^{or(p)+1}\kappa), f_j^\prime ( (-1)^{or(p) + 1}\kappa))$,
$(\alpha_{\kappa}(1), \alpha_\kappa^\prime (1)) =
(f_i((-1)^{or(p)}\kappa), f_i^\prime ( (-1)^{or(p)} \kappa))$ and
$im (\alpha)$ intersects the vanishing cycles in only these two
points. We illustrate
$\alpha_\kappa$ below for the two cases. \\[4pt]

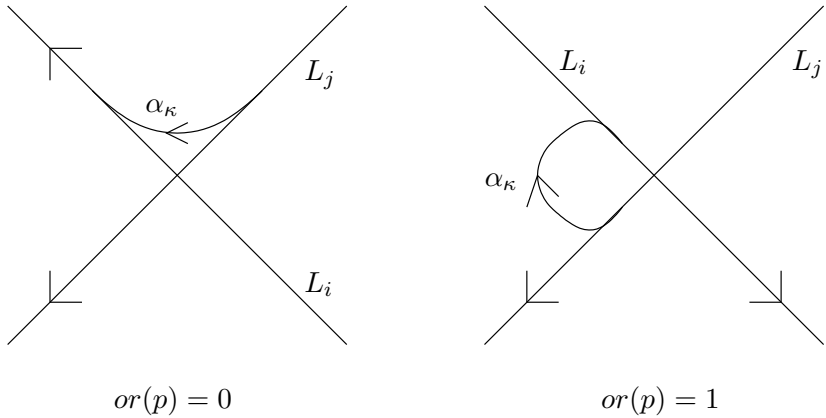
\begin{picture}(0,0)%
\includegraphics{alpha1.pstex}%
\end{picture}%
\setlength{\unitlength}{3947sp}%
\begingroup\makeatletter\ifx\SetFigFont\undefined%
\gdef\SetFigFont#1#2#3#4#5{%
  \reset@font\fontsize{#1}{#2pt}%
  \fontfamily{#3}\fontseries{#4}\fontshape{#5}%
  \selectfont}%
\fi\endgroup%
\begin{picture}(5320,2593)(2389,-2942)
\put(3067,-2890){\makebox(0,0)[lb]{\smash{{\SetFigFont{11}{13.2}{\rmdefault}{\mddefault}{\updefault}{$or(p) = 0$}%
}}}}
\put(4265,-828){\makebox(0,0)[lb]{\smash{{\SetFigFont{11}{13.2}{\rmdefault}{\mddefault}{\updefault}{$L_j$}%
}}}}
\put(4265,-2157){\makebox(0,0)[lb]{\smash{{\SetFigFont{11}{13.2}{\rmdefault}{\mddefault}{\updefault}{$L_i$}%
}}}}
\put(7327,-760){\makebox(0,0)[lb]{\smash{{\SetFigFont{11}{13.2}{\rmdefault}{\mddefault}{\updefault}{$L_j$}%
}}}}
\put(5862,-760){\makebox(0,0)[lb]{\smash{{\SetFigFont{11}{13.2}{\rmdefault}{\mddefault}{\updefault}{$L_i$}%
}}}}
\put(3267,-1026){\makebox(0,0)[lb]{\smash{{\SetFigFont{11}{13.2}{\rmdefault}{\mddefault}{\updefault}{$\alpha_\kappa$}%
}}}}
\put(5396,-1492){\makebox(0,0)[lb]{\smash{{\SetFigFont{11}{13.2}{\rmdefault}{\mddefault}{\updefault}{$\alpha_\kappa$}%
}}}}
\put(6129,-2890){\makebox(0,0)[lb]{\smash{{\SetFigFont{11}{13.2}{\rmdefault}{\mddefault}{\updefault}{$or(p) = 1$}%
}}}}
\end{picture}%
\\[4pt]

One sees that from this construction that $- \pi \leq \xi
(\alpha_\kappa) + or(p) \pi \leq 0$ and as $\kappa$ approaches zero,
$\xi (\alpha_\kappa) + or(p) \pi$ approaches $\arg (\eta (v_i)) -
\arg( \eta (v_j))$ mod $\pi \Z$. From this, one easily shows that
$$- \lim_{\kappa \to 0} [{1 \over \pi}({\tilde \phi_j}(p) - {\tilde
\phi_i(p)} + \xi (\alpha_\kappa) + or(p) \pi ) ] = -\lfloor {1 \over
\pi} ({\tilde \phi_j}(p) - {\tilde \phi_i(p)}) \rfloor = \mu (p)$$

Now suppose $p_i \in L_i$ and $p_j \in L_j$ and take the oriented
paths $\beta_j$ from $p_j$ to $f_j ((-1)^{or (p) + 1} \kappa)$ and
$\beta_i$ from $f_i((-1)^{or (p)} \kappa)$ to $p_i$. Let $\beta_{ji
\kappa}$ be the concatenation $\beta_i \circ \alpha_\kappa \circ
\beta_j$. Using the functorial properties of $\xi$ we have
\begin{displaymath}
\begin{array}{ccc} \xi (\beta_{ji \kappa})+ or(p) \pi  & = & \xi
(\beta_j)  + \xi (\beta_i)  + \xi (\alpha_{ji}) + or(p) \pi  \\ & =
& {\tilde \phi_j}(f_j ((-1)^{or (p) + 1} \kappa)) - {\tilde
\phi_j}(p_j) +  {\tilde \phi_i}(p_i) \\ & & - {\tilde
\phi_i}(f_i((-1)^{or (p)} \kappa)) + \xi(\alpha_\kappa)  + or(p) \pi
\\ & = & [{\tilde \phi_j}(f_j ((-1)^{or (p) + 1} \kappa)) - {\tilde
\phi_i}(f_i((-1)^{or (p)} \kappa)) \\ & & +
\xi(\alpha_\kappa) + or ( p) \pi] +  {\tilde \phi_i}(p_i) -  {\tilde \phi_j}(p_j)  \\
& = & -\pi \mu (p) +  {\tilde \phi_i}(p_i) -  {\tilde \phi_j}(p_j)
\end{array}
\end{displaymath}
We see that the last equality follows from the invariance of $\xi$
under isotopy. In other words, as $\kappa$ tends to zero we saw that
the second to last line approaches the last line, however, as these
paths are isotopic, their value under $\xi$ is constant, so we must
have equality. This equation gives \begin{equation} \mu (p) = -
\frac{1}{\pi}[\xi (\beta_{ji \kappa})+ or(p) \pi + {\tilde
\phi_j}(p_j) - {\tilde \phi_i}(p_i)] \end{equation}  Although this
approach may seem cumbersome at first, its utility lies in the fact
that we can calculate Maslov indices simply by knowing the path
$\beta_{ji}$, the monodromy data associated to $\eta$ and the values
of $\tilde \phi_i$ at a single point.

With these preliminaries in mind, we return to the fiber $\we$ and
prove our final theorem.

 \begin{thm} The map $\psi_{jk}$ induces an
isomorphism of triangulated categories from $\D (Fuk (W_\A , \C^* -
D_\epsilon , a, \{\delta_i \})) $ to $\T$.
\end{thm}

\begin{proof} The $1$-form $\eta$ is defined as the restriction to $\we$ of
a $1$-form $\tau$ on $(\C^*)^2$ which satisfies $${\tau} \wedge dW =
\frac{dz_1 \wedge dz_2}{z_1 z_2} := \Omega$$ While $\tau$ is not
defined globally on $(\C^*)^2$, it is well defined on the regular
fibers $\we$. One can see that $\Omega$ yields a non-zero phase
function on the Lefschetz thimbles $D_i$ and the restriction of
$\Omega$ to $L_i$ is precisely $\eta (L_i)$ (assuming one has
oriented $L_i$ properly). As the thimbles are simply connected, on
sees immediately that the vanishing cycles can be given
$\Z$-gradings via $\eta$.

Now, as $F_\varepsilon$ transversely intersects the divisor $\{z_i =
0\}$, we have that near such an intersection point, $\eta$ looks
like $\frac{dz}{z}$. Recall from the previous subsection that we
have a neighborhood $U_0 \subset \we$ where $U_0 = \{w : |w| < e, w
\not= 0, w \not= 1 \}$. As $w$ tends to zero in $U_0$, $z_1$ tends
to zero, while as $w$ tends to $1$, $z_0$ tends to zero. This
implies that on $U_0$,
$$\eta = h(w) \frac{dw}{w(w - 1)}$$ where $h(w)$ is a
nowhere vanishing holomorphic function. Of course, one can multiply
$\eta$ by $1/h(w)$ and it will not effect the Maslov indices.
Furthermore, perturbing $\eta$ near the boundary of $U_0$ we can
assume that $\eta = \frac{dw}{w^2} = - dw^{-1}$ near the boundary.
Although $\eta$ is no longer holomorphic, this perturbation will not
effect Maslov indices. From the previous subsection, we have a
description of the incoming and outgoing points of $L_i$ in $U_0$ as
$R_{i \pm} = \exp(\frac{2\pi}{ 4(l - 1)} Q_{i \pm})$. We can assume
that the tangent vector $v_{i \pm}$ at $R_{i \pm }$ is normal to the
circular boundary of $U_0$, pointing inward for $R_{i +}$ and
outward for $R_{i -}$. Thus, up to a positive real factor, $v_{i
\pm} = \mp e^{Im(\frac{2\pi}{ 4(l - 1)} Q_{i \pm})}
\partial_w$ and $\eta(v_{i \pm}) = \mp e^{\frac{2\pi}{ 4(l -
1)}(Im(Q_{i \pm}) - 2 Q_{i \pm})}$ so that $$\arg (\eta(v_{i +})) =
\pi - \frac{\pi}{ 2(l - 1)} Im (Q_{i +})$$ and $$\arg (\eta(v_{i
-})) = - \frac{\pi}{ 2(l - 1)} Im (Q_{i -})$$ We will give $L_i$ a
grading by letting
$${\tilde \phi_i}(R_{i +}) = \pi - \frac{\pi}{ 2(l - 1)} Im (Q_{i
+})$$ We would like then to find the real number ${\tilde
\phi_i}(R_{i -})$. Indeed, while we know that ${\tilde \phi_i}(R_{i
-}) \equiv  - \frac{\pi}{ 2(l - 1)} Im (Q_{i -})(mod 2\pi)$, we must
find the real number that is determined by our choice of ${\tilde
\phi_i}(R_{i +})$. In order to do this we will find $\xi (S_r)$ for
a certain path $S_r$ which we now describe. Define
$$S_r(t) = 1 + r e^{\frac{(2t - 1)\pi i}{2}}$$ to be the oriented
path of a semi-circle of radius $r$. We can pullback $\eta$ via the
map $e^{\frac{2\pi z}{ 4(l - 1)}}$. For the half space $Re(z) \geq
1$ we obtain the form
$${\tilde \eta}_z = \frac{\pi }{2(l - 1)} e^{-\frac{z \pi }{2(l -
1)}} dz$$ by extending $\eta$ to $-d(w^{-1}) $ outside of $U_0$ in
the complex plane. Thus we can transfer all calculations to the
covering space with parameter $z$ using the form $\tilde \eta$. One
should note that around points $4(l - 1)ai$, we have that ${\tilde
\eta}_z = k(z) \frac{dz}{z}$ with a non-zero holomorphic function
$k(z)$. Now let $G_r(t) = 1 + i(2rt - r)$ be the straight line
segment. One can see by perturbing the curve $S_r$ that $\xi ( S_r)
= \xi (G_r) + \pi$. Under the covering map, $G_r$ goes to a circular
arc. Calculating, we have
$${\tilde \eta}(G^\prime_r(t)) = {\tilde
\eta}_{G_r(t)}(2r\partial_y) = \frac{\pi r i}{(l - 1)} e^{-\frac{(1
+ i(2rt - r)) \pi }{2(l - 1)}}$$ So that $$\arg ({\tilde
\eta}(G^\prime_r(t)) = Im (\log ( {\tilde \eta}(G^\prime_r(t)))) =
\frac{\pi}{2} - \frac{\pi (2rt - r)}{2(l - 1)}$$ Implying $$\xi
(G_r) = - \frac{r\pi}{(l-1)}$$ Thus $\xi (S_r) = \pi - r\pi/ (l -
1)$. Also note that translating the curve $S_r$ along the line
$Re(z) = 1$ will not effect $\xi (S_r)$ as $\tilde \eta$ only picks
up a constant phase for such a translation. Using this calculation,
we can find the values ${\tilde \phi}_i (R_{i -})$. To do this just
assume that we have perturbed the corners of $\lag_i \subset \C$ so
that the tangent vectors at $Q_{i \pm}$ are normal to $Re(z) = 1$.
Then
$$F_i = \lag_i \cup (S_{Im[(Q_{i +} - Q_{i -})/2]} + (Q_{i +} + Q_{i
-})/2)$$ gives an oriented embedded closed curve which loops around
$0$ precisely once. This curve is freely homotopic to a small curve
around $0$ and as ${\tilde \eta} = k(z) \frac{dz}{z}$ in a
neighborhood of $0$ we have that $\xi (F_i ) = 0$. Using the
additivity of $\xi $ we then have \begin{displaymath}
\begin{array}{ccc}0 &  = & \xi (\lag_i) + \xi (S_{(Q_{i +} - Q_{i -}/2} +
(Q_{i +} + Q_{i -})/2)) \\& = & {\tilde \phi_i}(R_{i -}) - {\tilde
\phi_i}(R_{i +}) + \pi - \frac{\pi}{2 (l - 1)}Im(Q_{i +}) +
\frac{\pi}{2 (l - 1)}Im(Q_{i -})
 \\ & = & {\tilde \phi_i}(R_{i -}) + \frac{\pi}{2 (l - 1)}Im(Q_{i
 -}) \end{array} \end{displaymath}
So that $${\tilde \phi_i}(R_{i -}) = - \frac{\pi}{2 (l - 1)}Im(Q_{i
 -}) $$
We can now compute the Maslov indices of the various intersection
points. Let us start with the intersection points of the type $p \in
C_j \cap C_i$, i.e. those that are sent to the identity in $\wedge
V$. One can easily see either by computation or by examining the
figures in the previous subsection that $or(p) = 0$ for any such
point. Now, choose $p_j = Q_{j +}$ and $p_i = Q_{i -}$ as in our
setup for equation ($5$). Then the concatenated curve $\beta_{ji
\kappa}$ will be isotopic to an embedded curve $\beta_1$ starting at
$Q_{j +}$, ending at $Q_{i -}$ and going around $0$. As in the
derivation of ${\tilde \phi}_i (R_{i -})$ we see that $$\xi
(\beta_{ji \kappa}) = - \xi (S_{Im(Q_{j +} - Q_{i -})/2}) =
\frac{\pi}{2(l - 1)}Im(Q_{j +}) - \frac{\pi}{2(l - 1)}Im(Q_{i -}) -
\pi$$ implying by equation ($5$)
\begin{displaymath} \begin{array}{ccc} \mu (p)&  = & -\frac{1}{\pi}[ \xi (\beta_1) + {\tilde \phi}_i (R_{j +}) - {\tilde \phi}_i (R_{i -})] \\
&= & -\frac{1}{\pi}[\frac{\pi}{2(l - 1)}Im(Q_{j +}) - \frac{\pi}{2(l
- 1)}Im(Q_{i -}) - \pi + {\tilde \phi}_i (R_{j +}) - {\tilde \phi}_i
(R_{i -})] \\ & = & 0 \end{array}\end{displaymath} as expected.

For all other intersection points we see that $or(p) = 1$. Let $p
\in L_j \cap L_i$ be such a point that is sent to $e_0$ under the
isomorphism $\psi_{ji}$. We can describe the curve $\beta_{ji
\kappa}$ as an immersed curve starting at $Q_{j +}$, ending at $Q_{i
-}$, wrapping around $0$ twice and containing an "extra" loop as
drawn below.

\begin{picture}(0,0)%
\includegraphics{maslov1.pstex}%
\end{picture}%
\setlength{\unitlength}{3947sp}%
\begingroup\makeatletter\ifx\SetFigFont\undefined%
\gdef\SetFigFont#1#2#3#4#5{%
  \reset@font\fontsize{#1}{#2pt}%
  \fontfamily{#3}\fontseries{#4}\fontshape{#5}%
  \selectfont}%
\fi\endgroup%
\begin{picture}(5511,4071)(1195,-4579)
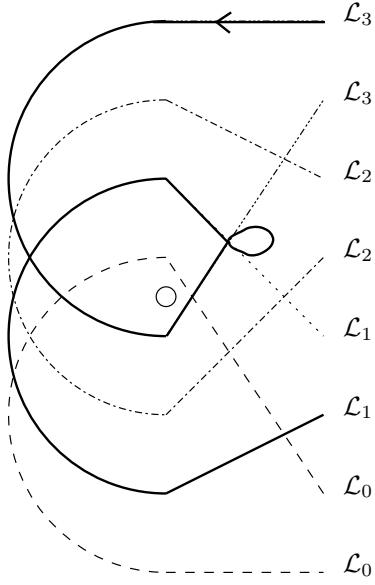
\put(1478,-4536){\makebox(0,0)[lb]{\smash{{\SetFigFont{10}{12.0}{\rmdefault}{\mddefault}{\updefault}{Curve approximating $\beta_{31\kappa}$}%
}}}}
\put(5108,-1622){\makebox(0,0)[lb]{\smash{{\SetFigFont{10}{12.0}{\rmdefault}{\mddefault}{\updefault}{$\beta_2$}%
}}}}
\put(6294,-2020){\makebox(0,0)[lb]{\smash{{\SetFigFont{10}{12.0}{\rmdefault}{\mddefault}{\updefault}{$\beta_3$}%
}}}}
\put(4689,-1212){\makebox(0,0)[lb]{\smash{{\SetFigFont{10}{12.0}{\rmdefault}{\mddefault}{\updefault}{$\beta_1$}%
}}}}
\put(3314,-636){\makebox(0,0)[lb]{\smash{{\SetFigFont{10}{12.0}{\rmdefault}{\mddefault}{\updefault}{$\lag_3$}%
}}}}
\put(3314,-1626){\makebox(0,0)[lb]{\smash{{\SetFigFont{10}{12.0}{\rmdefault}{\mddefault}{\updefault}{$\lag_2$}%
}}}}
\put(3314,-2121){\makebox(0,0)[lb]{\smash{{\SetFigFont{10}{12.0}{\rmdefault}{\mddefault}{\updefault}{$\lag_2$}%
}}}}
\put(3314,-2616){\makebox(0,0)[lb]{\smash{{\SetFigFont{10}{12.0}{\rmdefault}{\mddefault}{\updefault}{$\lag_1$}%
}}}}
\put(3314,-3110){\makebox(0,0)[lb]{\smash{{\SetFigFont{10}{12.0}{\rmdefault}{\mddefault}{\updefault}{$\lag_1$}%
}}}}
\put(3314,-3605){\makebox(0,0)[lb]{\smash{{\SetFigFont{10}{12.0}{\rmdefault}{\mddefault}{\updefault}{$\lag_0$}%
}}}}
\put(3314,-4101){\makebox(0,0)[lb]{\smash{{\SetFigFont{10}{12.0}{\rmdefault}{\mddefault}{\updefault}{$\lag_0$}%
}}}}
\put(3314,-1131){\makebox(0,0)[lb]{\smash{{\SetFigFont{10}{12.0}{\rmdefault}{\mddefault}{\updefault}{$\lag_3$}%
}}}}
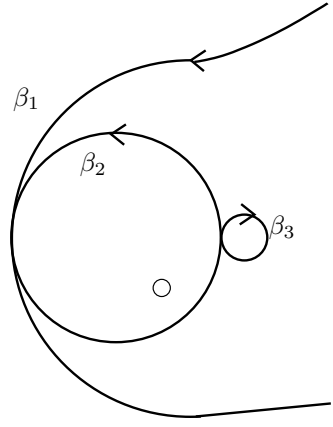
\put(4758,-4536){\makebox(0,0)[lb]{\smash{{\SetFigFont{10}{12.0}{\rmdefault}{\mddefault}{\updefault}{Curve isotopic to $\beta_{31\kappa}$}%
}}}}
\end{picture}%

 Thus, by perturbing $\beta_{ji \kappa}$,
it can be broken up into a curve $\beta_1$ from $Q_{j +}$ to $Q_{i
-}$ around $0$ as above, a loop $\beta_2$ around  $0$ and an extra
contractible loop $\beta_3$ oriented clockwise. Then using the
additivity of $\xi$ we have $\xi(\beta_{ji \kappa}) = \xi (\beta_1)
+ \xi (\beta_2) + \xi (\beta_3) = \xi (\beta_1) - 2\pi$. Thus, using
the above computation we have $\mu (p) = -\frac{1}{\pi}[\xi
(\beta_1) - 2\pi + \pi + {\tilde \phi}_i (R_{j +}) - {\tilde \phi}_i
(R_{i -})] = -\frac{1}{\pi}[\xi (\beta_1)  + {\tilde \phi}_i (R_{j
+}) - {\tilde \phi}_i (R_{i -})] + 1 = 1$ as expected.

Finally, assume $p \in L_j \cap L_i$ is sent to $e_1$ via
$\psi_{ji}$. The curve $\beta_{ji \kappa}$ is then isotopic to a
curve $\alpha_1$ that starts at $Q_{j +} - 4(l - 1)i$, ends at $Q_{i
-}$ and a contractible loop $\alpha_2$ oriented clockwise as
illustrated below. \\[2pt]

\begin{picture}(0,0)%
\includegraphics{maslov3.pstex}%
\end{picture}%
\setlength{\unitlength}{3947sp}%
\begingroup\makeatletter\ifx\SetFigFont\undefined%
\gdef\SetFigFont#1#2#3#4#5{%
  \reset@font\fontsize{#1}{#2pt}%
  \fontfamily{#3}\fontseries{#4}\fontshape{#5}%
  \selectfont}%
\fi\endgroup%
\begin{picture}(5909,1904)(1191,-2332)
\put(3428,-1189){\makebox(0,0)[lb]{\smash{{\SetFigFont{10}{12.0}{\rmdefault}{\mddefault}{\updefault}{$\lag_3$}%
}}}}
\put(5839,-1740){\makebox(0,0)[lb]{\smash{{\SetFigFont{10}{12.0}{\rmdefault}{\mddefault}{\updefault}{$Q_{0-}$}%
}}}}
\put(5839,-1259){\makebox(0,0)[lb]{\smash{{\SetFigFont{10}{12.0}{\rmdefault}{\mddefault}{\updefault}{$Q_{3+} - 4(l - 1)i$}%
}}}}
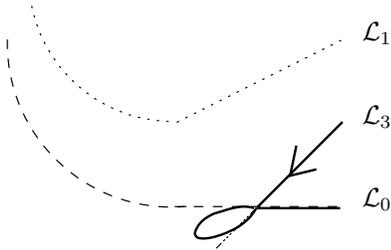
\put(1457,-2265){\makebox(0,0)[lb]{\smash{{\SetFigFont{10}{12.0}{\rmdefault}{\mddefault}{\updefault}{Curve approximating $\beta_{30\kappa}$}%
}}}}
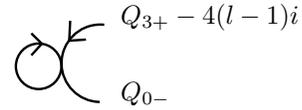
\put(4726,-2283){\makebox(0,0)[lb]{\smash{{\SetFigFont{10}{12.0}{\rmdefault}{\mddefault}{\updefault}{Curve isotopic to $\beta_{30\kappa}$}%
}}}}
\put(3428,-667){\makebox(0,0)[lb]{\smash{{\SetFigFont{10}{12.0}{\rmdefault}{\mddefault}{\updefault}{$\lag_1$}%
}}}}
\put(3428,-1712){\makebox(0,0)[lb]{\smash{{\SetFigFont{10}{12.0}{\rmdefault}{\mddefault}{\updefault}{$\lag_0$}%
}}}}
\end{picture}%

Now $\alpha_1 \circ S_{Im (Q_{j +} - Q_{i -})/2 - 2(l - 1)}$ forms a
counter clockwise contractible loop so that
\begin{displaymath} \begin{array}{ccc} 2 \pi & = & \xi (\alpha_1 \circ S_{Im (Q_{j +} -
Q_{i -})/2 - 2(l - 1)}) \\ & = & \xi (\alpha_1) + \xi(S_{Im (Q_{j +}
- Q_{i -})/2 - 2(l - 1)}) \\ & = & \xi (\alpha_1) + \pi -
\frac{\pi}{2(l - 1)}Im(Q_{j +}) + \frac{\pi}{2(l - 1)}Im(Q_{i -}) +
2 \pi \\ & = & 2 \pi + \xi ( \alpha_1 ) + {\tilde \phi_j}(R_{j +}) -
{\tilde \phi_i}(R_{i -})
\end{array}\end{displaymath}
 Therefore $$\xi (\alpha_1) ={\tilde \phi_i}(R_{i -}) - {\tilde \phi_j}(R_{j
 +})$$ and \begin{displaymath} \begin{array}{ccc} \mu (p) & = &
 -\frac{1}{\pi}[ \xi (\alpha_1) + \xi (\alpha_2) + \pi + {\tilde \phi}_i (R_{j +}) - {\tilde \phi}_i (R_{i
 -})] \\ & = & -\frac{1}{\pi}[- 2\pi + \pi] = 1
\end{array}\end{displaymath}
as expected. \end{proof}

\end{document}